\documentclass[11pt,oneside]{article}
\usepackage[english]{babel}

\usepackage{mabliautoref}

\usepackage[curve,matrix,arrow,cmtip]{xy}
\NoComputerModernTips

\newdir^{ (}{{}*!/-3pt/\dir^{(}}    
\newdir_{ (}{{}*!/-3pt/\dir_{(}}

\usepackage[ansinew]{inputenc}
\usepackage[dvipsnames,usenames]{xcolor}

\usepackage[scale=0.7,hmarginratio=1:1,vmarginratio=2:3]{geometry}
\usepackage{makecell}

\newcounter{themargin}
\usepackage{latexsym}
\usepackage{amsmath}
\usepackage{amstext}
\usepackage{amssymb}
\usepackage{amsbsy}
\usepackage{paralist}
\usepackage[active]{srcltx}
\usepackage[titles]{tocloft}

\usepackage{tikz}
\usepackage{verbatim}

\usetikzlibrary{arrows,positioning} 
\usetikzlibrary{shapes.multipart}
\tikzset{
    >=stealth',
    punkt/.style={
           rectangle,
           rounded corners,
           draw=black,
           text width=16em,
           minimum height=2em,
           text centered},
    pil/.style={
           ->,
           shorten <=2pt,
           shorten >=2pt,}
}

\newbox\dottobox
\setbox\dottobox=\hbox{$\longrightarrow$}
\setbox\dottobox=\hbox to\wd\dottobox{\hss$
\UseComputerModernTips\xymatrix@C=.5cm{\ar@{.>}[r]&\\}
                                      $\hss}

\newbox\leftdottobox
\setbox\leftdottobox=\hbox{$\longleftarrow$}
\setbox\leftdottobox=\hbox to\wd\leftdottobox{\hss$
\UseComputerModernTips\xymatrix@C=.5cm{\ar@{<.}[r]&\\}
                                      $\hss}

\newbox\dotintobox
\setbox\dotintobox=\hbox{$\longrightarrow$}
\setbox\dotintobox=\hbox to\wd\dotintobox{\hss$
\UseComputerModernTips\xymatrix@C=.5cm{\ar@{^{ (}.>}[r]&\\}
                                      $\hss}

\usepackage{GBplainAutF}

\usepackage{graphicx}

\renewcommand{\thesubsection}{\arabic{section}.\arabic{subsection}}

\theoremstyle{plain}   
\maketheorem{Thm}{Theorem}{theorem}
\maketheorem{Prop}{Proposition}{theorem}
\maketheorem{Prop-Def}{Proposition-Definition}{theorem}
\maketheorem{Cor}{Corollary}{theorem}
\maketheorem{Lem}{Lemma}{theorem}
\maketheorem{SubLem}{Sublemma}{theorem}
\maketheorem{Conj}{Conjecture}{theorem}
\newtheorem*{ThmInt}{Theorem}

\theoremstyle{definition}    %------------text not italic style------------------

\maketheorem{Def}{Definition}{theorem}
\maketheorem{DefRem}{Definition-Remark}{theorem}
\maketheorem{Claim}{Claim}{theorem}
\maketheorem{Alg}{Algorithm}{theorem}

\theoremstyle{remark}    %----------------also text not italic ------------------

\maketheorem{Rem}{Remark}{theorem}
\maketheorem{Ex}{Example}{theorem}
\maketheorem{Ques}{Question}{theorem}
\maketheorem{QuesRem}{Question-Remark}{theorem}
\maketheorem{Ass}{Assumption}{theorem}
\maketheorem{Con}{Convention}{theorem}

\hypersetup{
bookmarks,
bookmarksdepth=3,
bookmarksopen,
bookmarksnumbered,
pdfstartview=FitH,
colorlinks,backref,hyperindex,
linkcolor=RawSienna,
anchorcolor=BurntOrange,
citecolor=OliveGreen,
filecolor=BlueViolet,
menucolor=Yellow,
urlcolor=OliveGreen
}

\DeclareMathOperator{\crs}{\mathrm{crs}}
\newcommand{\ost}{\overline\st}
\newcommand{\ust}{\mathrm{ust}}
\newcommand{\ad}{\mathrm{ad}}

\newcommand{\ori}{\mathrm{or}}
\newcommand{\St}{\mathrm{St}}
\newcommand{\har}{\mathrm{har}}

\DeclareMathOperator{\Iw}{{\mathcal I}}

\DeclareMathOperator{\SB}{SB}

\newcommand{\Tempnewpage}%{}%
{\newpage}

\newcommand{\longto}{\longrightarrow}
\newcommand{\into}{\hookrightarrow}

\def\SMat#1#2#3#4{{{\renewcommand{\arraystretch}{.6} 
\addtolength{\arraycolsep}{-.4\arraycolsep}
\big( \begin{array}{cc} \scriptstyle #1 & \scriptstyle #2 \\ \scriptstyle #3 & \scriptstyle #4   \end{array}  \big)}}}

\title{ A Hecke-equivariant decomposition of spaces of Drinfeld cusp forms via representation theory, and an investigation of its subfactors}
\author{Gebhard B\"ockle, Peter Mathias Gr\"af and Rudolph Perkins
\\
}
\date{\today}
\begin{document}
\maketitle

\begin{abstract}
There are various reasons why a naive analog of the Maeda conjecture has to fail for Drinfeld cusp forms. Focussing on double cusp forms and using the link found by Teitelbaum between Drinfeld cusp forms and certain harmonic cochains, we observed a while ago that all obvious counterexamples disappear for certain Hecke-invariant subquotients of spaces of Drinfeld cusp forms of fixed weight, which can be defined naturally via representation theory. The present work extends Teitelbaum's isomorphism to an adelic setting and to arbitrary levels, it makes precise the impact of representation theory, it relates certain intertwining maps to hyperderivatives of Bosser-Pellarin, and it begins an investigation into dimension formulas for the subquotients mentioned above. We end with some numerical data for $A=\BF_3[t]$ that displays a new obstruction to an analog of a Maeda conjecture by discovering a conjecturally infinite supply of $\BF_3(t)$-rational eigenforms with combinatorially given (conjectural) Hecke eigenvalues at the prime~$t$.
\end{abstract}
\tableofcontents

\setcounter{secnumdepth}{3}

\setcounter{section}{0}

\section*{Introduction}

The conjecture of Maeda, which is now supported by much computational evidence but no theoretical insight, asserts that the spaces $S_k^{\mathrm{cl}}(\SL_2(\BZ))$ of classical modular forms of weight $k$ and level one consist of a single Hecke orbit under the natural action of $\Gal(\overline \BQ/\BQ)$ on systems of Hecke eigenvalues. This suggests that $S_k^{\mathrm{cl}}(\SL_2(\BZ))$ possesses no non-trivial decompositions into Hecke stable subspaces. Moreover the Maeda conjecture asserts that for fixed $k$ the Galois group associated to the Hecke field  is the symmetric group $S_m$ on $m=\dim S_k^{\mathrm{cl}}(\SL_2(\BZ))$ letters. Recently variants of the conjecture for levels $\Gamma$ other than one have been suggested; see \cite{Dieulefait-Pacetti-Tsaknias}. Once the subspace spanned by CM forms has been removed from $S_k^{\mathrm{cl}}(\Gamma)$, the number of Hecke orbits seems to be related to the number of inertial types of conductor given by the level -- but this does not fully explain what is observed.

For Drinfeld modular forms, the situation is different. Let $A$ be the coordinate ring of a smooth projective curve over a finite field of characteristic $p$ minus one point $\infty$ with quotient field $F$ and let $\Gamma\subset \GL_2(F)$ denote a congruence subgroup (throughout this introduction). Denote by $S_{k,l}(\Gamma)$ the space of Drinfeld cusp forms of weight $k$, type $l$ and level $\Gamma$. The work \cite{Bosser-Pellarin-DQMF} of Bosser and Pellarin for $A=\BF_q[t]$ and $\Gamma=\GL_2(A)$ gives maps $S_{k',l'}(\Gamma)\to S_{k,l}(\Gamma)$ under certain numerical conditions on $(k,l)$ and $(k',l')$, that have to do with the vanishing of certain binomial coefficients mod $p$ and require $k-k'$ to be even. The maps are given in terms of hyperderivatives and are Hecke equivariant up to a twist by a character, which for a prime $\Fp$ of $A$ is given by $\Fp^{(k-k')/2}$. Since these maps are in general neither trivial nor surjective, the spaces  $S_{k,l}(\Gamma)$ can contain non-trivial Hecke stable subspaces. In particular, a direct analog of a Maeda type conjecture is not possible for Drinfeld modular forms. An even simpler argument to dispute such an analog is to use the $p$-power map $S_{k,l}(\Gamma)\to S_{pk,pl}(\Gamma)$. It is injective and, up to a Frobenius-twist, Hecke equivariant; but clearly the target has in general a strictly larger dimension than the domain. It has also been known from the beginning that there are cuspidal Hecke eigenforms that are not doubly-cuspidal. However they could play a role similar to CM forms in the classical case.

Our starting point towards a systematic study of Maeda-style conjectures in the Drinfeld setting is the Hecke-equivariant residue isomorphism $S_{k,l}(\Gamma)\cong C_\har(V_{k,l})^\Gamma$ of Teitelbaum from \cite{Teitelbaum-Poisson} between spaces of cusp forms of a given weight $k$, type $l$ and level $\Gamma$ and spaces of $\Gamma$-invariant harmonic cochains with values in a certain $\GL_2$-representation $V_{k,l}$ that is a finite-dimensional $F$-vector space. 
The results in \cite{Teitelbaum-Poisson} are only fully developed for groups $\Gamma$ that have no prime-to-$p$ torsion, and the results are not developed in an adelic setting. In \cite{Boeckle-EiSh} the first author gave such an adelic setting, but only for {\em small} adelic level groups $\CK$.
As the group $\SL_2(A)$, which has prime-to-$p$ torsion, is the most natural candidate for studying Maeda-style behaviors in our setup, and as one needs a suitable Hecke theory, \autoref{Section2} of this article is devoted to removing this restriction and to proving the following theorem which generalizes both \cite{Teitelbaum-Poisson} and  \cite{Boeckle-EiSh}. Here, we denote by $S_{k,l}(\CK)$ the space of adelic Drinfeld cusp forms of weight $k$, type $l$ and level $\CK$, by $C^\ad_\har(V_{k,l},\CK)$ the corresponding space of adelic harmonic cochains and by $\St_\CK$ the adelic Steinberg module.  All three carry an action of a naturally defined Hecke algebra $\CH_\CK$.

\begin{ThmInt}[\autoref{thm:Sk-CHar-Adel}, \autoref{thm:CHar-St-global}]
\label{Thm-ThmInt}
Let $\CK$ be any compact open subgroup of $\GL_2(\BA_F^\infty)$.
\begin{compactenum}
\item 
There are Hecke-equivariant isomorphisms
\[
S_{k,l}(\CK)\stackrel\simeq\longto C^\ad_\har(V_{k,l},\CK)^{\GL_2(F)} \otimes_F\BC_\infty \quad \text{and}\quad
C^\ad_\har(V_{k,l},\CK)^{\GL_2(F)} \stackrel\simeq\longto \St_{\CK}\otimes_{\GL_2(F)}V_{k,l}.\]
\item 
The assignment $N\mapsto C^\ad_\har(N,\CK)^{\GL_2(F)}$ defines an exact functor from the category of $F[\GL_2(F)]$-modules of finite $F$-dimension to the category of $\CH_\CK$-modules of finite $F$-dimension.
\end{compactenum}
\end{ThmInt}
Let us mention that in the course of proving the above theorem, we provide an alternative description of the Hecke action in the adelic setting not given \cite{Boeckle-EiSh}, and we indicate a proof why the two actions agree. The action here is much simpler than that in \cite{Boeckle-EiSh}. The latter was motivated by translating the Hecke action on (adelic) Drinfeld cusp forms via the residue isomorphism. The remaining part of the proof is based on general cohomological considerations presented in \autoref{Section3}. A main observation is that the Steinberg modules $\St$ and $\St_\CK$ recalled in formulas \eqref{eq:Def-Steinberg} and \eqref{eqn:StAd-Defi}, respectively, have nice cohomological properties for all congruence subgroups $\Gamma$ or all compact open subgroups $\CK$, and not only $\Gamma$ that are $p'$-torsion free, or $\CK$ that are small.

Part (b) of the above theorem turns any composition series of $V_{k,l}$ into a Hecke-stable composition series of $S_{k,l}(\CK)$. Therefore in \autoref{Section4} we shall study in detail the representation theory of $F[\GL_2(F)]$-representations of finite $F$-dimension. Building on \cite{Bonnafe}, we classify the irreducible representation of $\GL_2(F)$ and explain how to algorithmically determine the simple factors of $V_{k,l}$; observe that this has no counterpart in characteristic zero, where the analogous $\GL_2(\BQ)$-representations are irreducible. Moreover in \autoref{Prop-DsAsMapOnReps} we provide non-trivial maps between representations $V_{k,l}$ and $V_{k',l'}$ under certain conditions on $(k,l)$ and $(k',l')$, to which we were  led to by $\cite{Bosser-Pellarin-DQMF}$. The kernels and images of these maps contribute to composition series of the $V_{k,l}$, but in general do not give a maximal composition series. For arbitrary $(k,l)$, we cannot describe a composition series of the $V_{k,l}$, not even algorithmically; we can only describe as a quotient its Jordan-H\"older factor of highest~weight. 

\autoref{Section5} explains how the maps between different representations $V_{k,l}$ of \autoref{Section4} turn into the hyperderivatives of $\cite{Bosser-Pellarin-DQMF}$ as well as the Frobenius map on Drinfeld cusp forms under the functor and the isomorphism of the above theorem, see \autoref{Cor-BP} and \autoref{Prop-FrDrin}. As a byproduct we provide natural extensions of the maps of $\cite{Bosser-Pellarin-DQMF}$ to any ring $A$ and any level subgroup $\CK$. Taking the existence of single-cuspidal eigenforms into account, our representation-theoretic approach covers all previously known obstructions to the irreducibility of $S_{k,l}(\SL_2(A))$ as a Hecke-module. 

Part (b) of the above theorem also makes it clear that from a representation theoretic viewpoint it is essential to understand the Hecke-modules $C^\ad_\har(N,\CK)^{\GL_2(F)}$ for simple $N$. A first step is to determine their dimension. If $\CK$ is small, or if $\Gamma$ has no prime-to-$p$ torsion, such formulas are implicit in \cite{Teitelbaum-Poisson}; for instance, for any such $\Gamma$ there exists a constant $c_\Gamma>0$ such that $\dim_F C_\har(V)^\Gamma=c_\Gamma \dim_F V$ for any $F[\GL_2(F)]$-module $V$ of finite $F$-dimension. This does not apply, however, to groups like $\SL_2(A)$, for which one should first study analogs of a Maeda conjecture. Therefore in  \autoref{Section6} we investigate the dimension of $C_\har(L_k)^{\SL_2(A)}$ for the irreducible $\SL_2(F)$-modules $L_k$, $k\ge0$, introduced in \autoref{Section4}, and for  $A=\BF_q[t]$. If $q$ is odd and $k$ is odd, the dimension is zero. Otherwise we expect 
\[\dim_F C_\har(L_k)^{\SL_2(A)} \approx \frac{\gcd(2,q^2-1)}{q^2-1}\dim_F L_k.\]
For $q=2,3,5$, we present closed formulas for all $k$, which support this expectation, see Propositions \ref{prop:DimLkQ3}, \ref{prop:DimLkQ2} and \ref{prop:DimLkQ5}.

Having understood the impact of representation theory, in \autoref{MaedaSection} we explore the possibility of there being a Maeda-style conjecture. We take $A=\BF_q[t]$ and $\Gamma=\SL_2(A)$, and we consider a natural increasing sequence of weights  $k_n$ for which the $\GL_2(F)$-representation $V_{k_n,l}$ is irreducible (for $\SL_2$ the types are irrelevant). Making use of computer algebra systems, we explicitly determine the characteristic polynomials of the Hecke operator $T_t$ on $C_\har(V_{k_n,l})^\Gamma$ at the prime $(t)\subset A$, disregarding single-cuspidal Hecke eigenforms. The outcome  of our experiments in \autoref{Table1} for $p=3$ caught us by surprise. In the weights $k_n$ searched, we discovered an increasing number (in $n$) of $F$-rational Hecke eigenforms of multiplicity one, that seem to be new to the literature. We have conjectural recipes for the number of such forms in weight $k_n$ and for the occurring $T_t$-eigenvalues, see \autoref{Table2} and \autoref{Conj:RatEigen}. The existence of these forms beyond the range of our data is open. Putting these eigenforms and the single cuspidal eigenforms aside, our data still gives no Maeda type conjecture. In our computations we encountered one or two Hecke orbits; in all cases their associated Galois group over $F$ is a symmetric group, as conjectured in the classical case. But we feel that more computations are needed -- which is difficult, since $k_n$ grows exponentially. Our search for a conceptual explanation for the existence of these special $F$-rational eigenforms has failed so far, but we plan to further investigate this question in future work.

\subsubsection*{Notation and conventions}
\label{Sec-Notation}
\begin{compactitem}
\item $F$ will denote a global function field with constant field $\BF_q$ of cardinality $q$ and characteristic~$p$, and we set $e:=\log_pq$, so that $q=p^e$.
\item We fix a place $\infty $ of $F$ and define $A$ as the subring of $F$ of functions regular away from $\infty$. This is the coordinate ring of a smooth affine curve with finite unit group and constant field $k$. 
\item We let  $\hat A$ be the profinite completion $\invlim A/\Fn$, the limit being over all non-zero ideals $\Fn$ of~$A$.
\item For any non-zero prime ideal $\Fp$ of $A$, we let $A_\Fp$ be the completion of $A$ at $\Fp$ and we let $F_\Fp$ be the fraction field of $A_\Fp$, i.e., the completion of $F$ at~$\Fp$.
\item We write $\BA_F^\infty=\hat A\otimes_A F$ for the ring of adeles away from $\infty$.
\item The completion of $F$ at $\infty$ will be $F_\infty$, its valuation ring $\CO_\infty\subset F_\infty$, and a uniformizer is $\pi\in F_\infty$.
\item The completion of an algebraic closure of $F_\infty$ is denoted by $\BC_\infty$.
\end{compactitem}

{\bf Acknowledgements:} We are grateful to A. Petrov for carrying out many computations at an~early stage of our investigations. We would also like to thank the anonymous referee for the valuable comments. G.B. and P.G. received support from the DFG within the  FG 1920 and the SPP~1489.

\section{Harmonic cochains and the Steinberg module}
\label{Section2}

In this section we recall results on harmonic cochains and the Steinberg module, and we briefly describe the link to Drinfeld cusp forms for $\GL_2$. 

\subsection{The local theory}

Let $\CT$ be the Bruhat-Tits-tree for $\PGL_2$ for $F_\infty$. We consider it as a graph with vertex set $\CT_0$ and oriented edge set $\CT_1^{\ori}$. For any oriented edge $e$, write $e^*$ for the edge opposite to $e$, and write $o(e)$ and $t(e)$ for the origin and the terminus of~$e$. The tree comes with a left-action of $\PGL_2(F_\infty)$, i.e., we have such an action on $\CT_0$ and on $\CT_1^\ori$ and the maps $o(\cdot)$ and $t(\cdot)$ are equivariant with respect to these actions. For a precise description of the tree, see \cite[II.1]{Serre-Trees}.

Let $N$ be an $F[\GL_2(F)]$-module. Define the space of $N$-valued harmonic cochains on $\CT$ as
\[C_\har(N)=\Big\{c\colon \CT^\ori_1\to N\mid c(e^*)=-c(e) \hbox{ for all }e\in \CT^\ori_1, \sum_{t(e)=v}c(e)=0\hbox{ for all }v\in \CT_0\Big\}.\]
The $F$-vector space $C_\har(N)$ carries an action of $\GL_2(F)$ by defining for $\gamma\in\GL_2(F)$ and $c\in C_\har(N)$ the map
\[ \gamma\circ c\colon \CT_1^\ori\to N , e \mapsto \gamma(c(\gamma^{-1} e)).\]

We will be interested in invariants of $C_\har(N)$ under certain subgroups of $\GL_2(F)$. To define them, let $P$ be a rank $2$ projective $A$-submodule of $F^2$, so that we have a group monomorphism $\Aut_A(P)\to\GL_2(F)$. For a non-zero ideal $\Fn$ of $A$ define $\Aut_A(P,\Fn)$ as the subgroup of elements of $\Aut_A(P)$ whose induced action on $P/\Fn P$ is trivial.
\begin{Def}
One calls a subgroup $\Gamma$ of $\GL_2(F)$ a congruence subgroup if there are $P$ and $\Fn$ as above such that $\Aut_A(P,\Fn)\subset\Gamma\subset\Aut_A(P)$.
\end{Def}
If $N$ has finite $F$-dimension, it is shown in \cite{Teitelbaum-Poisson}, building on \cite{Serre-Trees}, that for any congruence subgroup $\Gamma$ of $\GL_2(F)$ the space $C_\har(N)^\Gamma$ of $\Gamma$-invariant $N$-valued harmonic cochains is a finite-dimensional $F$-vector space.

Denote by $F^2$ the tautological $\GL_2(F)$-representation and write $\det$ for the representation on $F$ via the determinant $\GL_2(F)\to\GL_1(F)$. For any $F[\GL_2(F)]$-module $N$ we denote by $N^*$ the $F$-dual $\Hom_F(N,F)$; it is again an $F[\GL_2(F)]$-module. Of central interest to us are the $F[\GL_2(F)]$-modules 
\begin{equation}
V_{k,l}:= \Big(\det^{l-1}\otimes_F \Sym^{k-2} ((F^2)^*)\Big)^*.
\end{equation}

Let $S_{k,l}(\Gamma)$ denote the space of cusp forms of weight $k$ and type $l$ for $\Gamma$ from \cite{Gekeler}, for $A=\BF_q[t]$, or from \cite{Boeckle-EiSh}, in general. The following result is from \cite{Teitelbaum-Poisson}, though proofs in full generality are only given in \cite{Boeckle-EiSh}.
\begin{Thm}
For any congruence subgroup $\Gamma\subset \GL_2(F)$ there is a natural isomorphism of $\BC_\infty$-vector spaces.
\[S_{k,l}(\Gamma)\stackrel\simeq\longto C_\har(V_{k,l})^\Gamma\otimes_F\BC_\infty.\]
\end{Thm}
For $A=\BF_q[t]$ the isomorphism above is also one of Hecke modules. As it is more natural to discuss the Hecke action in an adelic context, we shall postpone this for the moment.

The description of Drinfeld cusp forms via harmonic cochains is a first combinatorial description of such forms. Another one is via the Steinberg module, which we recall next. We consider the projective space $\BP^1(F)$ over $F$ as a left $\GL_2(F)$-set via $\gamma\cdot (a:b):=(a:b)\gamma^{-1}$ for $\gamma\in\GL_2(F)$ and $(a:b)\in\BP^1(F)$. This induces a left $\BZ[\GL_2(F)]$-module structure on $\BZ[\BP^1(F)]$. The degree map $\deg\colon \BZ[\BP^1(F)]\to\BZ,\sum_i n_i P_i\mapsto \sum n_i$ is $\GL_2(F)$-equivariant. The projective line $\BP^1(F_\infty)$ can be interpreted as the boundary of $\CT$, its subset $\BP^1(F)$ is related to the pair $(\CT,\Gamma)$. The Steinberg module for $\GL_2(F)$ is defined as the left $\BZ[\GL_2(F)]$-module
\begin{equation}\label{eq:Def-Steinberg}
\St:=\BZ[\BP^1(F)]_0:=\kernel\Big(\deg\colon \BZ[\BP^1(F)]\to\BZ \Big).
\end{equation}

Recall that a group is called $p'$-torsion free, if all its torsion elements have order a power of $p$. If $\Gamma$ is a congruence subgroup of $\GL_2(F)$ that is $p'$-torsion free, then it is shown in \cite[II.2.9]{Serre-Trees} that $\St$ is a finitely generated projective $\BZ[\Gamma]$-module. We recall this in some detail, because it is relevant when we later compare the Hecke action defined below with the Hecke action from \cite{Boeckle-EiSh}. 

Following Serre  \cite{Serre-Trees}, we say that a simplex $s$ of $\CT$ is $\Gamma$-stable if $\Stab_\Gamma(s)$ is trivial. We write $\CT_1^{\ori,\Gamma\dash\st}$ for the set of $\Gamma$-stable oriented edges, and $\CT_0^{\Gamma\dash\st}$ for the set of $\Gamma$-stable vertices. By their definition, $\Gamma$ acts freely on $\CT_1^{\ori,\Gamma\dash\st}$ and on $\CT_0^{\Gamma\dash\st}$. Hence $\BZ[\CT_1^{\ori,\Gamma\dash\st}]$ and $\BZ[\CT_0^{\Gamma\dash\st}]$ are free $\BZ[\Gamma]$-modules. There is a natural boundary map $\partial_\Gamma\colon \BZ[\CT_1^{\ori,\Gamma\dash\st}]\to \BZ[\CT_0^{\Gamma\dash\st}]$, defined by mapping a $\Gamma$-stable edge $e$ to $[t(e)]^{\Gamma\dash\st}-[o(e)]^{\Gamma\dash\st}$ with the convention that for a simplex $s$ of $\CT$ we set $[s]^{\Gamma\dash\st}=[s]$, if $s$ is $\Gamma$-stable, and $[s]^{\Gamma\dash\st}=0$, otherwise; so indeed $\partial_\Gamma$ depends on $\Gamma$. The following is from \cite[II.2.9]{Serre-Trees}; the formulation follows \cite[\S~5.3 and Proposition 5.16]{Boeckle-EiSh}.
\begin{Lem}\label{lem:FactsOnSt}
Suppose that $\Gamma\subset \GL_2(F)$ is a $p'$-torsion free congruence subgroup. Then we have:
\begin{compactenum}
\item $\partial_\Gamma\colon \BZ[\CT_1^{\ori,\Gamma\dash\st}]\to \BZ[\CT_0^{\Gamma\dash\st}]$ is a surjective $\BZ[\Gamma]$-module homomorphism.
\item The action $[e]\mapsto -[e^*]$ induces an action of $\BZ/(2)$ on $\BZ[\CT_1^{\ori,\Gamma\dash\st}]$, and the resulting space of coinvariants $\BZ[\overline\CT_1^{\Gamma\dash\st}] :=\BZ[\CT_1^{\ori,\Gamma\dash\st}]/\{e+e^*\mid e\in\CT_1^{\ori,\Gamma\dash\st}\}$ is again a free $\Gamma$-module.
\item There is an induced surjective $\BZ[\Gamma]$-module homomorphism $\bar\partial_\Gamma\colon \BZ[\overline\CT_1^{\Gamma\dash\st}]\to \BZ[\CT_0^{\Gamma\dash\st}]$.
\item There is a natural isomorphism $ b_\Gamma\colon \St\to \kernel(\bar\partial_\Gamma)$ of $\BZ[\Gamma]$-modules.
\item The number of $\Gamma$-orbits of  $\CT_1^{\ori,\Gamma\dash\st}$ and of $\CT_0^{\Gamma\dash\st}$ are finite.
\item 
\label{lem:FactsOnStf}
$\St$ is a finitely generated projective $\BZ[\Gamma]$-module.
\end{compactenum}
\end{Lem}

\begin{Rem}\label{rem:NormalSt}
Note that if $\Gamma\subset \wt\Gamma \subset \GL_2(F)$ are congruence subgroups such that $\Gamma$ is $p'$-torsion free and normal in $\wt\Gamma$, then $\BZ[\overline\CT_1^{\Gamma\dash\st}]$ and $\BZ[\CT_0^{\Gamma\dash\st}]$ carry natural structures as $\BZ[\wt\Gamma]$-modules and the map $\bar\partial_{\Gamma}$ and the isomorphism $b_\Gamma$ are $\BZ[\wt\Gamma]$-module homomorphisms. To see this, one verifies that all objects and morphisms of the diagram in the proof of \cite[Proposition 5.16]{Boeckle-EiSh} carry a natural $\wt\Gamma$-action.
\end{Rem}

Let us give an explicit description of $b_\Gamma$ from \autoref{lem:FactsOnSt}(d), as it will be used in the proof of the following result.\footnote{We do assume throughout this paragraph that the reader is familiar with \cite[II.2.9]{Serre-Trees}. We implicitly recall $b_\Gamma^{-1}$ in the second paragraph above \autoref{lem:FactsOnStAd}.}
 As a $\BZ$-module, $\St$ is free, and a basis is given by the elements $(a:1)-(1:0)$ with $a\in F$. Thus it will suffice to describe $b_\Gamma(P'-P)$ for any pair of distinct points $P,P'$ of $\BP^1(F)$. Let $\wp_{P\to P'}$ be the geodesic in $\CT$ from $P$ to $P'$. We think of $\wp_{P\to P'}$ as a sequence $(e_i)_{i\in\BZ}$ of edges such that $t(e_{i-1})=o(e_i)$ for $i\in\BZ$. The geodesic $\wp_{P\to P'}$ is then characterized by requiring that this path has no back-tracking, and that the half lines $(e_i)_{i\le0}$ and $(e_i)_{i\ge0}$ represent the boundary points $P$ and $P'$, respectively. From \cite[II.2.9]{Serre-Trees} one deduces: 
\begin{equation}\label{eq-bGamma-1}
 0\neq b_\Gamma(P'-P)=\sum_{i\in\BZ} [e_i]^{\Gamma\dash\st}  \in \kernel(\bar\partial_\Gamma)\subset \BZ[\overline\CT_1^{\Gamma\dash\st}].
\end{equation}
 The sum is finite: there are infinite half lines contained in $(e_i)_{i\le0}$ and in $(e_i)_{i\ge0}$, that end in $P$ and $P'$, respectively, on which $\Stab_\Gamma(e_i)$ is non-trivial (and growing in $|i|$) and contained in $\Stab_\Gamma(P)$ or $\Stab_\Gamma(P')$, respectively. The non-vanishing is clear since $b_\Gamma$ is an isomorphism and $P'-P$ is non-zero in $\St$.

The following lemma will be used to prove the comparison result \autoref{prop-TwoIGammaN}.
\begin{Lem}\label{Lem-ForTwoIGammaN}
Suppose $\Gamma'\subset\Gamma$ are $p'$-torsion free congruence subgroups of $\GL_2(F)$. Define maps 
\[p_{0,\Gamma'\subset\Gamma}\colon \BZ[\CT_0^{\Gamma'\dash\st}] \to \BZ[\CT_0^{\Gamma\dash\st}],[v]^{\Gamma'\dash\st} \mapsto [v]^{\Gamma\dash\st} \ \hbox{ and } \ \overline p_{1,\Gamma'\subset\Gamma}\colon \BZ[\overline\CT_1^{\Gamma'\dash\st}] \to \BZ[\overline\CT_1^{\Gamma\dash\st}],[e]^{\Gamma'\dash\st} \mapsto [e]^{\Gamma\dash\st} .\] 
Then the following diagram is commutative:
\[
\xymatrix{
0\ar[r] & \St \ar[r]^-{b_{\Gamma'}} \ar@{=}[d]&  \BZ[\overline\CT_1^{\Gamma'\dash\st}] \ar[r]^-{\bar\partial_{\Gamma'}} \ar[d]^-{\overline p_{1,\Gamma'\subset\Gamma}}& \BZ[\CT_0^{\Gamma'\dash\st}] \ar[r] \ar[d]^-{p_{0,\Gamma'\subset\Gamma}}&0 \\
0\ar[r] & \St \ar[r]^-{b_{\Gamma}} &  \BZ[\overline\CT_1^{\Gamma\dash\st}] \ar[r]^-{\bar\partial_\Gamma} & \BZ[\CT_0^{\Gamma\dash\st}] \ar[r] &0 \\
}
\]
\end{Lem}
\begin{proof}
The maps $p_{0,\Gamma'\subset\Gamma}$ and $\overline p_{1,\Gamma'\subset\Gamma}$ are clearly well-defined, since we describe them on bases. The commutativity of the diagram is an immediate consequence of the description of $b_\Gamma$ in \eqref{eq-bGamma-1} -- and our definitions of the maps $\bar\partial_\Gamma$,  $p_{0,\Gamma'\subset\Gamma}$ and $\overline p_{1,\Gamma'\subset\Gamma}$ and $t\mapsto [t]^{\Gamma\dash\st}$.
\end{proof}

Now Teitelbaum in \cite{Teitelbaum-Poisson} defines a map $\iota_{N,\Gamma}$ of $F$-vector spaces as follows.\footnote{ As introduced, $\CT_0$, $\CT_1$ and $\St$ carry a left $\GL_2(F)$-action. We regard them as right modules by defining $(\cdot)\gamma:=\gamma^{-1}(\cdot)$.} He considers
\begin{equation}\label{eq:Def-CharToSt}
C_\har(N)^{\Gamma}\to   \BZ[\overline\CT_1^{\Gamma\dash\st}]\otimes_{\BZ[\Gamma]} N,c\mapsto \sum_{e\in \Gamma\backslash\CT_1^{\ori,\Gamma\dash\st}/\{\pm 1\}} e\otimes_{\Gamma} c(e),
\end{equation}
and he verifies that that this map is well-defined and if composed with $\bar\partial_{\Gamma}\otimes_{\Gamma} \id_N$ is zero. If now $\Gamma$ is $p'$-torsion free, then $\St\otimes_{\Gamma} N\cong\kernel(\bar\partial_{\Gamma}\otimes_{\Gamma} \id_N) $, because $ \BZ[\CT_0^{\Gamma\dash\st}]$ is a projective $\BZ[\Gamma]$-module. One obtains an induced map $\iota_{N,\Gamma}\colon C_\har(N)^{\Gamma}\to \St\otimes_{\Gamma} N$.
\begin{Thm}[{\cite[Prop.~21]{Teitelbaum-Poisson}}\footnote{ Strictly speaking, Teitelbaum proves the result only for $V_{k,l}$. But his proof on \cite[pp.~504--506]{Teitelbaum-Poisson} carries over with no changes to all $N$.}]
\label{thm:CHar-St-local}
Suppose that $\Gamma$ is a $p'$-torsion free congruence subgroup of $\GL_2(F)$. Then the map  $\iota_{\Gamma,N}$ is an $F$-vector space isomorphism
\[ C_\har(N)^\Gamma\to \St\otimes_{\Gamma} N. \]
In particular $\iota_{\Gamma,V_{k,l}}\colon C_\har(V_{k,l})^\Gamma\to \St\otimes_{\Gamma} V_{k,l}$ is an isomorphism.
\end{Thm}
We shall remove the constraint that $\Gamma$ be $p'$-torsion free in \autoref{thm:CHar-St-local-Gen} and show further that the functor $\St\otimes_\Gamma\ublank$ is exact on $F[\GL_2(F)]$-modules.

\subsection{The global adelic theory}

To define a Hecke action, for general $F$ and $\infty$ with class number possibly different from $1$, we adelize the above situation. We fix a compact open subgroup $\CK\subset\GL_2(\BA_F^\infty)$. As before, we let $N$ be an $F[\GL_2(F)]$-module. 

\medskip

Let $e\in\CT_1^\ori$, $v\in\CT_0$ and $g\in\GL_2(\BA_F^\infty)$. We shall write  $[e,g]_\CK$ and $[v,g]_\CK$ for respective classes in  $ \CT_{1,\CK}^\ori:= \CT^\ori_1\times\GL_2(\BA_F^\infty)/\CK$ and  $\CT_{0,\CK}:=\CT_0\times\GL_2(\BA_F^\infty)/\CK$. Both $ \CT_{1,\CK}^\ori$ and $\CT_{0,\CK}$ carry a left $\GL_2(F)$-operation by $\gamma[s,g]_\CK:=[\gamma s,\gamma g]_\CK$ for $s$ a simplex and $\gamma\in\GL_2(F)$. The space $C^\ad_\har(N,\CK)$ of {harmonic cochains of level $\CK$} is the set of all maps $c\colon \CT^\ori_{1,\CK} \to N$ such that for all $g\in\GL_2(\BA_F^\infty)$ we have
\[c([e^*,g]_\CK)=-c([e,g]_\CK) \hbox{ for all }e\in \CT^\ori_1\hbox{ and } \sum_{t(e)=v}c([e,g]_\CK)=0\hbox{ for all }v\in \CT_0.\]
It is a left $F[\GL_2(F)]$-module by defining for $\gamma\in\GL_2(F)$ and $c\in C_\har^\ad(N,\CK)$ the map
\[ \gamma\circ c\colon \CT^\ori_{1,\CK} \to N , [e,g]_\CK \mapsto \gamma(c(\gamma^{-1} [e,g]_\CK)).\]

By strong approximation the determinant map induces a bijection of double cosets 
\begin{equation}\label{eq:StrongApp}
\det\colon \GL_2(F)\backslash\GL_2(\BA_F^\infty)/\CK\longto\GL_1(F)\backslash\GL_1(\BA_F^\infty)/\det(\CK)=:\Cl_\CK.
\end{equation}
The right hand side is a finite extension of the class group of $A$, and hence a finite group. Hence there is a tuple $(g_c)_{c\in \Cl_\CK}$ in $\GL_2(\BA_F^\infty)$, such that $\GL_2(\BA_F^\infty)$ is equal to the disjoint union $\coprod_{c\in \Cl_\CK}\GL_2(F)g_c\CK$. The groups $\Gamma_{g_c}:=g_c\CK g_c^{-1}\cap \GL_2(F)$ are congruence subgroups and one can construct a natural isomorphism between spaces of invariants
\begin{equation}\label{eq-LocalGlobal}
 C^\ad_\har(N,\CK)^{\GL_2(F)} \cong \bigoplus_{c\in \Cl_\CK} C_\har(N)^{\Gamma_{g_c}}.
\end{equation}
In particular, the space $C^\ad_\har(N,\CK)^{\GL_2(F)}$ of $\GL_2(F)$-invariant $N$-valued adelic harmonic cochains has finite $F$-dimension whenever $N$ has finite $F$-dimension.

The space $C^\ad_\har(N,\CK)^{\GL_2(F)} $ carries an action by Hecke operators. Let $y\in\GL_2(\BA_F^\infty)$. Then $\CK y\CK\subset \GL_2(\BA_F^\infty)$ is compact and hence a finite disjoint union $\CK y\CK=\coprod_i y_i\CK$ of right translates of the open subgroup $\CK$. Define the Hecke operator ${}_{|\CK y\CK}$ for $y$ by mapping $c\in C^\ad_\har(N,\CK)$ to
\[ c{}_{|\CK y\CK} \colon \CT^\ori_{1,\CK} \to N , [e,g]_\CK \mapsto \sum_{i} c([e,gy_i]). \]
The Hecke action clearly commutes with the left $\GL_2(F)$-action and thus preserves $C^\ad_\har(N,\CK)^{\GL_2(F)}$.  It is also clear that the operator ${}_{|\CK y\CK} $ is independent of the chosen $y_i$. The operators ${}_{|\CK y\CK} $ do in general not preserve the `local' factors in \eqref{eq-LocalGlobal}, i.e., the direct summands on the right.

Let $S_{k,l}(\CK)$ denote the space of adelic cusp forms of weight $k$ and type $l$ for $\CK$ from \cite{Boeckle-EiSh}. It is the natural generalization to the adelic setting of the local definition of $S_{k,l}(\Gamma)$ used by Goss, Gekeler et al. The following is obtained as in \cite[Sects.~5,~6]{Boeckle-EiSh}, building on \cite{Teitelbaum-Poisson} after taking into account the different normalizations regarding the action on $V_{k,l}$.
\begin{Thm}
\label{thm:Sk-CHar-Adel}
For any compact open subgroup $\CK\subset \GL_2(\BA_F^\infty)$ there is an isomorphism
\[S_{k,l}(\CK)\stackrel\simeq\longto C^\ad_\har(V_{k,l},\CK)^{\GL_2(F)} \otimes_F\BC_\infty,\]
which is equivariant for all Hecke operators ${}_{|\CK y\CK}$ for $y\in\GL_2(\BA_F^\infty)$.
\end{Thm}
The map is given by an explicit residue map which we will recall in \autoref{Sec-HyperDer}.

\medskip

Next we recall the adelic Steinberg module. We define it as the kernel in the short exact sequence
\begin{equation}\label{eqn:StAd-Defi}
0\to \St_\CK\to \BZ[\BP^1(F)\times \GL_2(\BA_F^\infty)/\CK]\stackrel{\deg_\CK}\to \BZ[\GL_2(\BA_F^\infty)/\CK]\to 0,
\end{equation}
with $\deg_\CK=\deg\otimes\id_{\GL_2(\BA_F^\infty)/\CK}$. The sequence carries a left $\BZ[\GL_2(F)]$-action induced from the left actions of $\GL_2(F)$ on $\BP^1(F)$ as given before \eqref{eq:Def-Steinberg} and on $\GL_2(\BA_F^\infty)/\CK$ by left multiplication.

 For $y\in\GL_2(\BA_F^\infty)$ and $\CK y\CK=\coprod_i y_i\CK$ from above, we define a Hecke action on $\St_\CK$ as follows: first define
\[
_{|\CK y\CK}\colon \BZ[\BP^1(F)\times \GL_2(\BA_F^\infty)/\CK]\to \BZ[\BP^1(F)\times \GL_2(\BA_F^\infty)/\CK], (P,g\CK)\mapsto \sum_i(P,gy_i\CK),\]
\[_{|\CK y\CK}\colon \BZ[\GL_2(\BA_F^\infty)/\CK]\to \BZ[\GL_2(\BA_F^\infty)/\CK], [g\CK]\mapsto \sum_i [gy_i\CK].
\]
It is clear that $\deg_\CK$ is equivariant for these definitions of ${}_{|\CK y\CK}$, and that the maps are independent of the chosen representatives $y_i$. Moreover the Hecke operation commutes with the $\GL_2(F)$-action since one acts from the right and the other from the left. Hence ${}_{|\CK y\CK}$ lies in $\End_{\BZ[\GL_2(F)]}(\St_\CK)$.

The advantage of the above action is that it works for any compact open subgroup $\CK$ of $\GL_2(\BA_F^\infty)$. To compare it with \cite[Sect.~6]{Boeckle-EiSh}, we need a definition.
\begin{Def}\label{def:KisSmall}
We call $\CK$ {\em small} if for all $g\in\GL_2(\BA_F^\infty)$ the congruence subgroup $\Gamma_g:=g\CK g^{-1}\cap \GL_2(F)$ is $p'$-torsion free.
\end{Def}
We refer the reader to \autoref{lem:KisSmall} for some criteria for $\CK$ to be small.

The construction of Hecke operators on $\St_\CK$ in  \cite[Sect.~6]{Boeckle-EiSh} was only given for small $\CK$. The construction was adapted to the adelization of the isomorphism  from \autoref{thm:CHar-St-local}. To continue, the following result is important.
\begin{Prop}\label{prop:HeckeAgreesForKsmall}
Suppose that $\CK$ is small. Then for any $y\in\GL_2(\BA_F^\infty)$, the operator ${}_{|\CK y\CK}$ constructed here and the operator ${}_{|\CK y\CK}$ constructed in  \cite[Sect.~6]{Boeckle-EiSh} agree.
\end{Prop}
For the proof we first need some preparations similar to \cite[\S~6.4]{Boeckle-EiSh}, that basically amount to an adelic version of \autoref{lem:FactsOnSt}. Suppose that $\CK$ is small. A simplex $[s,g]_\CK$ of the union of trees $\CT_\CK:=\CT\times \GL_2(\BA_F^\infty)/\CK$ is called $\CK$-stable if $\Stab_{\GL_2(F)}([s,g]_\CK)$ is trivial. Define sets of $\CK$-stable simplices $\CT_{0,\CK}^{\st}$ and $\CT_{1,\CK}^{\ori,\st}$ and observe that by definition $\GL_2(F)$ acts freely on these. One has a natural boundary map $\partial^\st_\CK\colon \BZ[\CT_{1,\CK}^{\ori,\st}]\to \BZ[\CT_{0,\CK}^{\st}]$ given by mapping $[e,g]_\CK$ to $[t(e),g]^\st_\CK-[o(e),g]^\st_\CK$, with the convention that for any simplex $[s,g]_\CK$ of $\CT_\CK$ we set  $[s,g]_\CK^\st=[s,g]_\CK$, if $[s,g]_\CK$ is $\CK$-stable, and $[s,g]_\CK^\st=0$, otherwise, as in the definition of $\partial_\Gamma$. One also has the usual boundary map $\partial_\CK\colon \BZ[\CT_{1,\CK}^{\ori}]\to \BZ[\CT_{0,\CK}]$ from algebraic topology. Let $\CT_{1,\CK}^{\ori,\ust} =\CT_{1,\CK}^{\ori}  \setminus \CT_{1,\CK}^{\ori,\st} $ and $\CT_{0,\CK}^{\ust} =\CT_{0,\CK}  \setminus \CT_{0,\CK}^{\st} $. This defines a subgraph of $\CT_\CK$, and the restriction of $\partial_\CK$ gives the standard boundary map $\BZ[\CT_{1,\CK}^{\ori,\ust}]\to \BZ[\CT^\ust_{0,\CK}]$. One has a $\BZ/(2)$-action $[e,g]_\CK\mapsto - [e^*,g]_\CK$ on the modules $\BZ[\CT_{1,\CK}^{\ori}]$ and $\BZ[\CT_{1,\CK}^{\ori,\st}]$ that commutes with $\partial_\CK$ and $\partial^\st_\CK$. For $*\in\{\st,\ust,\emptyset\}$, the modules of $\BZ/(2)$-coinvariants are denoted by $\BZ[\overline\CT_{1,\CK}^{*}]$, and we put a bar on the induced boundary maps. The $\BZ/(2)$-action also preserves $\BZ[\CT_{1,\CK}^*\cap (\CT_1\times\{g\CK\} )]$, and we write $ \BZ[\overline\CT_{1,\CK}^*\cap ( \overline\CT_1\times\{g\CK\} )]$ for the module of coinvariants. We observe that $\BZ[\overline\CT_{1,\CK}^{\st}]$ is again a free $\GL_2(F)$-module. 

Now for each $g\in\GL_2(F)$ consider the following commutative diagram with exact rows
\[\xymatrix@R-.5pc{
0\ar[r]& \BZ[\overline\CT_{1,\CK}^{\ust}\cap ( \overline\CT_1\times\{g\CK\} )]\ar[r] \ar[d]^{\bar\partial_\CK}&\BZ[ \overline\CT_{1,\CK}\cap ( \overline\CT_1\times\{g\CK\} )]\ar[r] \ar[d]^{\bar\partial_\CK}&\BZ[\overline\CT_{1,\CK}^{\st}\cap ( \overline\CT_1 \times\{g\CK\} )]\ar[r] \ar[d]^{\bar\partial_\CK^\st}&0\\
0\ar[r]& \BZ[\CT_{0,\CK}^{\ust}\cap ( \CT_0\times\{g\CK\} )]\ar[r] &\BZ[\CT_{0,\CK}\cap ( \CT_0\times\{g\CK\} )]\ar[r] &\BZ[\CT_{0,\CK}^{\st}\cap ( \CT_0\times\{g\CK\} )]\ar[r] &0\\
}\]
and where the horizontal maps are the natural choices. To verify the properties of the diagram and its well-definedness the reader should convince herself that the following identifications hold. One has  $\overline\CT_{i,\CK}^{\st}\cap ( \overline\CT_1\times\{g\CK\} )=\overline\CT^{\Gamma_{g}\dash\st}_i\times\{g\CK\}$, for $i=0,1$ and then the left vertical arrow is equal to $\bar\partial_{\Gamma_g}\times {\id_{\{g\CK\}}}$. This makes accessible all results of \autoref{lem:FactsOnSt} to the above diagram. The cokernel of the right vertical arrow is $0$, as are the kernels of the middle and left vertical maps, the kernel of the right vertical arrow is $\kernel(\bar\partial_{\Gamma_g}) \times\{g\CK\}$, the cokernel of the middle vertical arrow is $\BZ\times\{g\CK\}$ and that of the left one is $\BZ[\BP^1(F)]\times\{g\CK\}$. The map $b^{-1}_{\Gamma_g}$ is the connecting homomorphism in the Snake Lemma. See also \cite[II.2.9]{Serre-Trees} and in particular the proof of Lemma 13 therein. 

Choose now for each class $c\in \Cl_\CK$ from \eqref{eq:StrongApp} a representative by $g_c\in\GL_2(\BA_F^\infty)$, and choose representatives for the sets $\GL_2(F)/\Gamma_{g_c}$, $c\in\Cl_\CK$. Taking the direct sum $\bigoplus_{c\in \Cl_\CK} \bigoplus_{\gamma \in \GL_2(F)/\Gamma_{g_c}} $ over the above diagrams with $g$ being replaced by $\gamma g_c$, we obtain the above diagram with the expressions $\cap ( \overline\CT_1\times\{g\CK\} )$ removed in the top row and with $\cap ( \CT_0\times\{g\CK\} )$ removed in the bottom row, i.e., the diagram in the adelic setting! From \autoref{lem:FactsOnSt} and $\GL_2(F)/\Gamma_{g_c}\cong \GL_2(F)\backslash\GL_2(F)g_c\CK/\CK$, we deduce
\begin{Lem}\label{lem:FactsOnStAd}
Suppose that $\CK\subset \GL_2(\BA_F^\infty)$ is small. Then we have:
\begin{compactenum}
\item \label{lem:StFromBoe}
The sequence $0\to \kernel \bar\partial_\CK^\st\to \BZ[\overline\CT_{1,\CK}^{\st}]\stackrel{\bar\partial_\CK^\st}\to \BZ[\CT_{0,\CK}^{\st}]\to 0$ is exact.
\item The $\BZ[\GL_2(F)]$-modules $\BZ[\overline\CT_{1,\CK}^{\st}]$ and $\BZ[\CT_{0,\CK}^{\st}]$ are free and finitely generated.
\item There is a natural isomorphism $ b_\CK \colon \St_\CK\to \kernel(\bar\partial_\CK^\st)$ of $\BZ[\GL_2(F)]$-modules.
\item 
$\St_\CK$ is a finitely generated projective $\BZ[\GL_2(F)]$-module.
\end{compactenum}
\end{Lem}
Let us also describe $b_\CK$. Let $P,P'$ be distinct elements of $\BP^1(F)$ and $g\in\GL_2(\BA_F^\infty)$. Then as indicated below \autoref{lem:FactsOnSt} in the local situation, there is a geodesic path  $\wp_{P\to P',g\CK}$ from $(P,g\CK)$ to $(P',g\CK)$. Let $([e_i,g]_\CK)_{i\in\BZ}$ denote the sequence of successive edges forming $\wp_{P\to P',g\CK}$. Then \eqref{eq-bGamma-1} interpreted in the adelic situation gives
\begin{equation}\label{eq-bCK-1}
0\neq b_\CK((P',g\CK)-(P,g\CK))=\sum_{i\in\BZ} [e_i,g]_\CK^\st. 
\end{equation}

\begin{proof}[Proof of \autoref{prop:HeckeAgreesForKsmall}]
The Hecke action introduced here on $\St_\CK$ is defined via compatible actions on the middle and right terms of the complex \eqref{eqn:StAd-Defi}. The Hecke action in  \cite[Sect.~6]{Boeckle-EiSh} is defined via a compatible action on the middle and right term of the complex in \autoref{lem:FactsOnStAd}\ref{lem:StFromBoe}. The key for the comparison is the isomorphism $b_\CK$. Note that the elements $(P',g\CK)-(P,g\CK)$ generate $\St_\CK$ over $\BZ$, when $g\CK$ traverses $\GL_2(\BA_F^\infty)/\CK$ and $P,P'$ traverse all distinct pairs of points of $\BP^1(F)$.

Fix $g\CK$ and distinct points $P',P$ from $\BP^1(F)$, and represent $\wp_{P\to P',g\CK}$ by its sequence of successive edges $([e_i,g]_\CK)_{i\in\BZ}$. For a stable edge $[e,g]_\CK$, in  \cite[\S~6]{Boeckle-EiSh} its course $\crs([e,g]_\CK)$, a geodesic path of edges, is defined as follows: $\crs([e,g]_\CK)$ always contains the stable edge $[e,g]_\CK$. If $[t(e),g]_\CK$ is stable, the path ends at $[t(e),g]_\CK$; otherwise there is a unique boundary point $(Q',g\CK)$ such that 
\[\Stab_{\GL_2(F)}([t(e),g]_\CK)= \Stab_{\Gamma_g} (t(e)) \subset \Stab_{\Gamma_g} ( Q')= \Stab_{\GL_2(F)}( (Q',g\CK)),\]
and in that case $\crs([e,g]_\CK)$ contains the infinite geodesic path from $[t(e),g]_\CK$ to $(Q',g\CK)$. If $[o(e),g]_\CK$ is stable, the path begins at $[o(e),g]_\CK$; otherwise there is a unique boundary point $(Q,g\CK)$ such that $\Stab_{\GL_2(F)}([o(e),g]_\CK)  \subset  \Stab_{\GL_2(F)}( (Q,g\CK))$, and in that case $\crs([e,g]_\CK)$ contains the infinite geodesic path from $(Q,g\CK)$ to $[o(e),g]_\CK$. An important observation concerning $\wp_{P\to P',g\CK}$ is that in the formal completion of $\BZ[\overline\CT_{1,\CK}]$, one has the equality (of infinite sums) 
\[
\sum_{i\in\BZ} \Big(\sum_{[e',g]_\CK \in \crs ([e_i,g]_\CK^\st)} [e',g]_\CK\Big)
=
\sum_{i\in\BZ} [e_i,g]_\CK. 
\]
On the left, there are only finitely many $i$ for which the inner sum is non-zero; there are always inner sums that are infinite, and so we do need to consider this in the formal completion. Also, there may be infinite cancellations on the left; such a cancellation occurs if $[e_i,g]_\CK$ and $[e_{i'},g]_\CK$ are stable for some $i<i'$, and if $[t(e_i),g]_\CK$, $[o(e_{i'}),g]_\CK$ and all edges in between are unstable, because then $\crs([e_i,g]_\CK)$ and $\crs([e_{i'},g]_\CK)$ share all edges except for $[e_{i''},g]_\CK$, $i<i''<i'$, but with opposite orientation.
 
Let $y\in\GL_2(\BA_F^\infty)$ and write $\CK y\CK=\coprod_{j\in J} y_j\CK$. 
Now  \cite[\S~6]{Boeckle-EiSh} defines 
\[[e_i,g]_\CK^\st|_{\CK y\CK} =  \sum_{[e',g]_\CK \in \crs ([e_i,g]_\CK^\st)} \sum_{j\in J} [e',gy_j]^\st_\CK;\] 
here the right hand side is finite, because the sum is over stable edges only. It follows for the Hecke action as defined in \cite[\S~6]{Boeckle-EiSh} applied to $b_\CK((P',g\CK)-(P,g\CK))\stackrel{\eqref{eq-bCK-1}}=\sum_{i\in\BZ} [e_i,g]_\CK^\st$ that one has
\[ b_\CK\big((P',g\CK)-(P,g\CK)\big)|_{\CK y\CK}=\sum_{i\in \BZ} \sum_{j\in J} [e_i,gy_j]^\st_\CK.\] 
The Hecke action defined here on $\St_\CK=\kernel (\deg_\CK)$ is given directly by
\[ ((P',g\CK)-(P,g\CK)) |_{\CK y\CK}=\sum_{j\in J} \big((P',gy_j\CK)-(P,gy_j\CK)\big).\]
To see that the actions agree, we need to show that $b_\CK((P',gy_j\CK)-(P,gy_j\CK))= \sum_{i\in \BZ} [e_i,gy_j]^\st_\CK$ for all $j\in J$. The point is that the geodesic $\wp_{P\to P',gy_j\CK}$ is the sequence of successive edges $([e_i,gy_j]_\CK)_{i\in\BZ}$, and altering the  label $g\CK$ into $gy_j\CK$ maps a geodesic to a geodesic. This completes the proof.
\end{proof}

One can now generalize the map $\iota_{\Gamma,N}$ from \eqref{eq:Def-CharToSt} to the adelic setting. This gives the $F$-linear homomorphism $\iota_{\CK,N}\colon C^\ad_\har(N,\CK)^{\GL_2(F)}\to \St_\CK\otimes_{\GL_2(F)} N$,  defined by 
\begin{equation}\label{eq-AdelicTeitelbaumMap}
\iota_{\CK,N}(c) = \sum_{[e,g]_\CK \in\GL_2(F)\backslash \CT_{\CK,1}^{\ori,\st}/\{\pm1\}}
[e,g]_\CK \otimes c([e,g]_\CK)  \in \BZ[\overline\CT_{\CK,1}^{\ori,\st}]\otimes_{\GL_2(F)}N
\end{equation}
for any $F[\GL_2(F)]$-module $N$ of finite $F$-dimension. Because of \autoref{prop:HeckeAgreesForKsmall} we can quote from  \cite[\S~6]{Boeckle-EiSh} the following result where the Hecke action on $\St_\CK$ is the one defined here.

\begin{Thm}[{\cite[\S~6.4]{Boeckle-EiSh}}]
\label{thm:CHar-St-global-Small}
If $\CK$ is small, then $\iota_{\CK,N}$ defines Hecke-equivariant isomorphisms
\[C^\ad_\har(N,\CK)^{\GL_2(F)} \stackrel\simeq\longto \St_\CK\otimes_{\GL_2(F)}N.\]
\end{Thm}
\begin{proof}
Being an isomorphism as $F$-vector spaces is deduced from \autoref{thm:CHar-St-local} and \eqref{eq-LocalGlobal}. The Hecke-equivariance is seen as in \cite[\S~6.4]{Boeckle-EiSh}.
\end{proof}
We shall remove the requirement that $\CK$ is small  in \autoref{thm:CHar-St-global} and also show that the functor $\St_\CK\otimes_{\GL_2(F)} \ublank$ is exact on $F[\GL_2(F)]$-modules.

\subsection{Extensions of \autoref{thm:CHar-St-local} and \autoref{thm:CHar-St-global-Small}}
\label{Subsec2.3}
In this subsection we extend \autoref{thm:CHar-St-local} and \autoref{thm:CHar-St-global-Small} in the way promised after their formulation using some more general cohomological results proved in \autoref{Section3}.
One motivation is that it should be clarified that the above results hold for all congruence subgroups $\Gamma$ of $\GL_2(F)$ or all compact open subgroups $\CK$ of $\GL_2(\BA_F^\infty)$, respectively. Another is that the comparison to \cite{Bosser-Pellarin-HyperDiff} in \autoref{Sec-HyperDer}, as well as our computations for $A=\BF_q[t]$ in \autoref{MaedaSection} towards a Maeda type conjecture are for $\Gamma\in\{\SL_2(A),\GL_2(A)\}$ which are not $p'$-torsion free.

Let first $\Gamma$ be a congruence subgroup of $\GL_2(F)$, say $\Aut_A(P,\Fn)\subset \Gamma\subset \Aut_A(P)$ for some projective rank $2$ $A$-submodule $P$ of $F^2$ and some non-zero ideal $\Fn$ of $A$. Let $\Fp$ be any non-zero prime ideal of $A$. Then the homomorphism $\Aut_A(P)\to\Aut_{A_\Fp}(P\otimes_AA_\Fp)\cong\GL_2(A_\Fp)$ is injective. The kernel of $\Aut_{A_\Fp}(P\otimes_AA_\Fp)\to \Aut_{A/\Fp}(P/\Fp P)$ is an open pro-$p$ subgroup of $\GL_2(A_\Fp)$. Let $\Gamma_\Fp$ denote the intersection of this kernel with $\Gamma$; then clearly $\Gamma_\Fp$ is $p'$-torsion free. By construction, it is normal in $\Gamma$, and since $\Aut_A(P,\Fn\Fp)\subset \Gamma_\Fp$, it is a congruence subgroup.

Fix $\Fp$ as above and let $N$ be a $F[\GL_2(F)]$-module. Consider the isomorphism $\iota_{\Gamma_\Fp,N}\colon C_\har(N)^{\Gamma_\Fp}\to \St\otimes_{\Gamma_\Fp} N$. There is an obvious action of $\bar\Gamma:=\Gamma/\Gamma_\Fp$ on $C_\har(N)^{\Gamma_\Fp}$ and an action on $\St\otimes_{\Gamma_\Fp} N$ defined by
\begin{equation}
 \bar\gamma\cdot  m\otimes_{\Gamma_\Fp}n = m\gamma^{-1}\otimes_{\Gamma_\Fp}\gamma n,
\end{equation}
for $m\in \St$, $n\in N$ and $\gamma\in\Gamma$. One verifies that $\iota_{\Gamma_\Fp,N}$ induces an isomorphism $ C_\har(N)^{\Gamma}\to (\St\otimes_{\Gamma_\Fp} N)^{\bar \Gamma}$, using the explicit formula \eqref{eq:Def-CharToSt}. The norm map considered in \autoref{prop:OnGroupHomology}(b) gives an isomorphism $ \St\otimes_{\Gamma} N \stackrel\simeq\to (\St\otimes_{\Gamma_\Fp} N)^{\bar \Gamma}$, and we write $\tilde\iota_{\Gamma,N}$ for the combined isomorphism $C_\har(N)^{\Gamma}\to \St\otimes_{\Gamma} N$. 
\begin{Prop}\label{prop-TwoIGammaN}
If $\Gamma$ is $p'$-torsion free, then $\tilde\iota_{\Gamma,N}=\iota_{\Gamma,N}$.
\end{Prop}
\begin{proof}
We need to show that the following diagram commutes:
\[
\xymatrix{
C_\har(N)^{\Gamma_\Fp}\ar[rr]_-{\iota_{\Gamma_\Fp,N}}^-\simeq &&\St\otimes_{\Gamma_\Fp}N\\
\ar@{^{ (}->}[u]
C_\har(N)^{\Gamma}\ar[r]_-{\iota_{\Gamma,N}}^-\simeq &\St\otimes_{\Gamma}N\cong (\St\otimes_{\Gamma_\Fp}N)_{\overline\Gamma}\ar[r]_-{\nu_{N,\overline\Gamma}}&(\St\otimes_{\Gamma_\Fp}N)^{\overline\Gamma}\ar@{^{ (}->}[u]\rlap{,}\\
}\]
where $\nu_{N,\overline\Gamma}$ denotes the norm map  $m\otimes_\Gamma v\mapsto\sum_{\bar\gamma\in\overline\Gamma} m\gamma^{-1}\otimes_{\Gamma_\Fp} \gamma v$ from \eqref{eq:NormMap}. Let $c$ be in $C_\har(N)^{\Gamma}$. Going up and right gives
\[\iota_{\Gamma_\Fp,N}(c)=\sum_{e\in\Gamma_\Fp\backslash\CT_1^{\ori,\Gamma_\Fp\dash\st}/\{\pm\}} e\otimes_{\Gamma_\Fp}c(e) .\]
Because $\iota_{\Gamma_\Fp,N}(c)$ lies in $\St\otimes_{\Gamma_\Fp} N\subset \BZ[\overline\CT_1^{\Gamma_\Fp\dash\st}] \otimes_{\Gamma_\Fp}N$, it remains invariant if we apply to it $\overline p_{1,\Gamma_\Fp\subset\Gamma}\otimes \id_N$ with $\overline p_{1,\Gamma_\Fp\subset\Gamma}$ from \autoref{Lem-ForTwoIGammaN}. Hence 
\begin{equation}\label{eq-DiagIotaCommutes}
\iota_{\Gamma_\Fp,N}(c)=\sum_{e\in\Gamma_\Fp\backslash\CT_1^{\ori,\Gamma\dash\st}/\{\pm\}} e\otimes_{\Gamma_\Fp}c(e).
\end{equation}
Observe that in the summation index, we have $\Gamma_\Fp$-orbits of $\Gamma$-stable elements. If we follow the bottom path, then
\[ c\mapsto \sum_{e\in\Gamma\backslash\CT_1^{\ori,\Gamma\dash\st}/\{\pm\}} e\otimes_{\Gamma}c(e) \mapsto
 \sum_{e\in\Gamma\backslash\CT_1^{\ori,\Gamma\dash\st}/\{\pm\}} 
 \sum_{\bar\gamma\in\overline\Gamma} 
 e\gamma^{-1}\otimes_{\Gamma_\Fp} \gamma c(e).
\]
Now observe that  $e\gamma^{-1}=\gamma e$, that $\gamma c(e)=c(\gamma e)$ and that $\Gamma=\cup_{\bar\gamma\in\overline\Gamma}\gamma\Gamma_\Fp$. It follows that the last expression is equal to 
\eqref{eq-DiagIotaCommutes}, which shows that the diagram commutes.
\end{proof}
\begin{Rem}
In the above construction of $\tilde\iota_{\Gamma,N}$, we can replace $\Fp$ by any proper ideal $\Fm\subset A$ that is non-zero, and work with $\Gamma_\Fm$ instead of $\Gamma_\Fp$, and it is straightforward to observe that \autoref{prop-TwoIGammaN} also holds in that generality. It then follows in particular that $\tilde\iota_{\Gamma,N}$ is independent of any choice, and that for $p'$-torsion free $\Gamma$  the map is the one defined by Teitelbaum. 

Thus from now on we simply write $\iota_{\Gamma,N}$ instead of $\tilde\iota_{\Gamma,N}$ for all congruence subgroups $\Gamma$ of $\GL_2(F)$.
\end{Rem}
Next we state an immediate consequence of \autoref{prop:OnGroupHomology}(c):
\begin{Thm}
\label{thm:CHar-St-local-Gen}
Let $\Gamma\subset \GL_2(F)$ be a congruence subgroup, let $A$ be a ring of characteristic $p>0$ and let $N$ be an $A[\Gamma]$-module. Then the following hold:
\begin{compactenum}
\item The map  $\iota_{\Gamma,N}$ defines an $F$-vector space isomorphism
\[ C_\har(N)^\Gamma\to \St\otimes_{\Gamma} N. \]
\item The functor $ \St\otimes_{\Gamma} \ublank\colon$ from $A[\Gamma]$- to $A$-modules is exact.
\end{compactenum}
In particular $\iota_{\Gamma,V_{k,l}}\colon C_\har(V_{k,l})^\Gamma\to \St\otimes_{\Gamma} V_{k,l}$ is an isomorphism for all~$k,l$.
\end{Thm}

We now pass to the adelic situation and use the notation from \ref{Subsec2.2}. We denote by $\CK''\subset \CK'\subset \CK$ compact open subgroups of $\GL_2(\BA_F^\infty)$ such that $\CK''$ is normal in $\CK$. We consider $G:=\GL_2(F)$ as a subgroup of $\GL_2(\BA_F^\infty)$ via the diagonal embedding $F\into \BA_F^\infty$. Thereby $S:=\GL_2(\BA_F^\infty)/\CK''$ carries a left action by $G$; it also carries an obvious right action by $H:=\CK/\CK''$; a short computation reveals that $S$ is free as a right $H$-set, so that condition \eqref{eq:AssOnStabHs} in \autoref{Subsec2.2} is satisfied. Let $M:=\St$ be the Steinberg module. Then
\[ M[S] = \bigoplus_{s\in S} M = {\St_{\CK''}}\]
is a  right $\BZ[G\times H]$-module by the formula in \eqref{eqn:DefGH-Mod}, we have $M[S/H]=M[S]_H=\St_\CK$, and $\St_{\CK''}$ is free $H$-module. In the present situation, using $\St_{\CK''}=\St\otimes_\BZ \BZ[S]$, the action is given by 
\begin{equation}\label{eq-GxH-action}
(m\otimes [g\CK''])\cdot(\gamma,h)=(\gamma^{-1}m)\otimes [\gamma g h^{-1} \CK'']
\end{equation}
 for $(\gamma,h)\in\GL_2(F)\times \CK/\CK''$. Note also that for any $g\in \GL_2(\BA_F^\infty)$ the group $\Stab_G(g\CK'')=G\cap g\CK''g^{-1}$ is a congruence subgroup of $\GL_2(F)$. 

We remind the reader that for any $\CK$ there exists a free $\wh A$-submodule $\Lambda$ of $(\BA^\infty_F)^2$ of rank $2$ that is invariant under the action of $\CK$, so that $\CK\subset\Aut_{\wh A}(\Lambda)$; this uses that for any free $\wh A$-submodule of rank $2$ of $(\BA_F^\infty)^2$, the stabilizer under the action of $\CK$ is open in $\CK$. We choose any such $\Lambda$ for the following lemma. The lemma gathers result on smallness of $\CK$ and on smallness of $\CK'/\CK''$ for $\St_{\CK''}$.

\begin{Lem}\label{lem:KisSmall}
For any proper non-zero ideal $\Fn\subset A$, define $\CK(\Fn):=\kernel (\Aut_{\wh A}(\Lambda) \to \Aut_{A/\Fn}(\Lambda/\Fn\Lambda))$. We assume that $\CK''=\CK\cap\CK(\Fn)$ for some such $\Fn$. Then the following hold.
\begin{compactenum}
\item
The group $\CK(\Fn)$ is small.
\item
If $\CK'/\CK''\subset\CK/\CK''$ is a $p$-Sylow subgroup, then $\CK'$ is small.
\item 
Any $p$-Sylow subgroup $H'\subset H$ is small for $\St_{\CK''}$ in the sense of \ref{Subsec2.2}, i.e. $\St_{\CK''}$ is projective as a $\Stab_G(sH')$-module for all  $s\in S$.
\end{compactenum}
\end{Lem}
\begin{proof}
The proof of (a) and (b) is a standard argument: Let  $g\in\GL_2(\BA_F^\infty)$. Using the first isomorphism theorem in group theory, one sees that $g\CK''g^{-1}\cap\GL_2(F)\subset g\CK'g^{-1}\cap\GL_2(F)$ is a normal subgroup of finite $p$-power index. Hence to prove (b), it suffices to prove that $\CK''$ is small. Since smallness is inherited by compact open subgroups, it suffices to prove (a), and only in the case when $\Fn$ is a maximal ideal of $A$. Denote by $g_\Fn$ the component of $g$ at $\Fn$ and by $\CK_{\Fn}$ the component of $\CK(\Fn)$ at $\Fn$. Then $\CK_\Fn$ is a compact open pro-$p$ subgroup of $\GL_2(F_\Fn)$ and so is $g_\Fn\CK_{\Fn}(g_\Fn)^{-1}$. Moreover we have $\GL_2(F)\cap g\CK(\Fn) g^{-1}\subset \GL_2(F)\cap g_\Fn\CK_{\Fn}(g_\Fn)^{-1}$ as a subgroup of $\GL_2(F)$. Since pro-$p$ groups do not contain elements of finite order prime to $p$, part (a) follows.

In light of (b), to prove (c) it suffices to show that if $\CK'$ is small, then $H'=\CK'/\CK''$ is small for $\St_{\CK''}$. However if $\CK'$ is small, then for any $g\in \GL_2(\BA_F^\infty)$ the group $\Gamma_g:=g\CK' g^{-1}\cap\GL_2(F)$ is $p'$-torsion free. And then by \autoref{lem:FactsOnSt}(f), the module $\St$ is a projective $\BZ[\Gamma_g]$-module. Upon observering that  $\Stab_G(sH')=g\CK' g^{-1}\cap\GL_2(F)=\Gamma_g$ for $s=g\CK''$, this proves (c).
\end{proof}

We do have a norm map for the adelic Steinberg module relative to $\CK''\unlhd\CK$:
\begin{Lem}
The following map is an isomorphism of $\BZ[G]$-modules
\begin{equation}\label{eq-NormMap}
 \St_\CK \cong (\St_{\CK''})_H \longto (\St_{\CK''})^H, [m\otimes [g\CK'']]_H \mapsto m\otimes \bigg(\sum_{h\in H} [gh\CK'']\bigg).
\end{equation}
\end{Lem}
\begin{proof}
One deduces from \eqref{eq-GxH-action} that the map is one of $\BZ[G]$-modules. To see the bijectivity, by \eqref{eqn:StAd-Defi} it is enough to see that
\[  \BZ[\GL_2(\BA_F^\infty)/\CK] \to  \BZ[\GL_2(\BA_F^\infty)/\CK'']^H, [g\CK]\mapsto \sum_{h\in H}[gh\CK'']  \]
is bijective. But this is clear, since the target space is constant on $H$ orbits and because $H=\CK/\CK''$.
\end{proof}

From \autoref{lem:KisSmall}, the remarks preceding it and Propositions \ref{prop:AdelicGroupHomol} and \ref{prop:AbstractHeckeAction}, we conclude.
\begin{Cor}\label{cor:AdelicStProps}
Let $A$ be a ring of characteristic $p>0$. Let  $N$ be an $A[\GL_2(F)]$-module. Suppose that $\CK''=\CK\cap\CK(\Fn)$ for a proper non-zero ideal $\Fn\subset A$. Let $\CH_\CK$ be the subalgebra of $\End_{\BZ[\GL_2(F)]} (\St_\CK)$ generated by the Hecke operators $|_{\CK y\CK}$, $y\in\GL_2(\BA_F^\infty)$. 
Then the following hold:
\begin{compactenum}
\item
\label{cor:AdelicStPropsa}
The module $\St_{\CK''}\otimes_{\GL_2(F)}N$ is cohomologically trivial as an $A[\CK/\CK'']$-module.
\item
\label{cor:AdelicStPropsb}
The norm map $\nu_{\CK,\CK'',N}\colon \St_{\CK}\otimes_{\GL_2(F)}N  \longto(  \St_{\CK''}\otimes_{\GL_2(F)}N)^{\CK/\CK''}$ induced from \eqref{eq-NormMap} is an isomorphism.
\item
\label{cor:AdelicStPropsc}
For any short exact sequence $0\to N'\to N\to N''\to 0$ of left $A[\GL_2(F)]$-modules, the sequence
\[ \xymatrix{
0\ar[r]& \St_\CK \otimes_{\GL_2(F)} N'\ar[r]&  \St_\CK \otimes_{\GL_2(F)} N\ar[r]& \St_\CK \otimes_{\GL_2(F)} N''\ar[r]&0\\
}
\]
is exact. 
\item 
\label{cor:AdelicStPropsd}
With the action from \autoref{prop:AbstractHeckeAction}, the modules $(\St_{\CK''}\otimes_{\GL_2(F)}N)^{\CK/\CK''}$ and $\St_{\CK}\otimes_{\GL_2(F)}N$ are $\CH_\CK$-modules, the norm map is $\CH_\CK$-equivariant, the map $N\mapsto\St_\CK \otimes_{\GL_2(F)} N$ from  $A[\GL_2(F)]$-modules to $\CH_\CK$-modules is functorial, and the sequence in (c) is one of $\CH_\CK$-modules.
\end{compactenum}
\end{Cor}

Let $\Fn\subset A$ be any proper non-zero ideal and let $\CK'':=\CK(\Fn)\cap \CK$. Let $N$ be a $F[\GL_2(F)]$-module of finite $F$-dimension. Using \eqref{eq-AdelicTeitelbaumMap} one verifies that the isomorphism $\iota_{\CK'',N}\colon C^\ad_\har(N,\CK'')^{\GL_2(F)} \to \St_{\CK''}\otimes_{\GL_2(F)}N$ from \autoref{thm:CHar-St-global-Small} is equivariant for the action of $H=\CK/\CK''$, and it therefore induces an isomorphism 
\[C^\ad_\har(N,\CK)^{\GL_2(F)} \to (\St_{\CK''}\otimes_{\GL_2(F)}N)^{\CK/\CK''}.\]
By \autoref{cor:AdelicStProps}(b) we have an isomorphism  $\St_{\CK}\otimes_{\GL_2(F)}N \to (\St_{\CK''}\otimes_{\GL_2(F)}N)^{\CK/\CK''}$ given by the norm map, and together we obtain an isomorphism 
\[\tilde\iota_{\CK,N}\colon C^\ad_\har(N,\CK)^{\GL_2(F)} \to \St_{\CK}\otimes_{\GL_2(F)}N.\]

\begin{Prop} \label{prop:IndepOfAdelicIota}
If $\CK$ is small, then $\tilde\iota_{\CK,N}=\iota_{\CK,N}$.
\end{Prop}
The proposition immediately implies that $\tilde\iota_{\CK,N}$ is independent of any choices. After having given its proof, we shall drop the tilde from the notation.

The proof of \autoref{prop:IndepOfAdelicIota} requires the following analog of \autoref{Lem-ForTwoIGammaN}:

\begin{Lem}\label{lem-ForAdelicTwoIotas}
Let $\CK\subset\GL_2(\BA_F^\infty)$ be a small compact open subgroup. Let $\CK''\subset \CK$ be a normal open subgroup, and set $H=\CK/\CK''$. Define $\CT^{\CK\dash\st}_{0,\CK''}\subset \CT_{0,\CK''}$ and $\CT^{\ori,\CK\dash\st}_{1,\CK''}\subset \CT^{\ori}_{1,\CK''}$ as the subsets of those simplices $[t,g]_{\CK''}$ such that $[t,g]_\CK$ is $\CK\dash$stable, and define $[t,g]_{\CK''}^{\CK\dash\st}$ to be $[t,g]_{\CK''}$ if $[t,g]_\CK$ is stable, and let it be zero otherwise. Then the following hold.
\begin{compactenum}
\item The sets $\CT^\st_{0,\CK''}$ , $\CT^{\CK\dash\st}_{0,\CK''}$, $\CT^{\ori,\st}_{1,\CK''}$ and $\CT^{\ori,\CK\dash\st}_{1,\CK''}$ are invariant under the action of $H$.
\item We have $\CT^{\CK\dash\st}_{0,\CK''}\subset \CT^\st_{0,\CK''}$ and $\CT^{\ori,\CK\dash\st}_{1,\CK''}\subset \CT^{\ori,\st}_{1,\CK''}$.
\item Let $p_{0,\CK''\subset\CK}\colon \BZ[\CT^{\st}_{0,\CK''}]\to\BZ[\CT^{\CK\dash\st}_{0,\CK''}]$ and $\bar p_{1,\CK''\subset\CK}\colon \BZ[\overline\CT^{\st}_{1,\CK''}]\to\BZ[\overline\CT^{\CK\dash\st}_{1,\CK''}]$ be defined by $[s,g]_{\CK''}^\st\mapsto[s,g]_{\CK''}^{\CK\dash\st}$, let $\partial_{\CK''}^{\CK\dash\st}$ be defined by $[e,g]_{\CK''}^{\CK\dash\st}\mapsto [t(e),g]_{\CK''}^{\CK\dash\st}-[o(e),g]_{\CK''}^{\CK\dash\st}$, and let $b^{\CK\dash\st}_{\CK''}$ be as in \eqref{eq-bCK-1} for $\CK''$ with $[e_i,g]_{\CK''}^\st$ replaced by $ [e_i,g]_{\CK''}^{\CK\dash\st}$. Then the following diagram commutes
\[\xymatrix@R-.4pc@C+1pc{
0\ar[r] & \St_{\CK''}\ar[r]^-{b_{\CK''}}& \BZ[\overline\CT_{1,\CK''}^{\st}]\ar[r]^-{\bar\partial_{\CK''}^\st}&\BZ[\CT_{0,\CK''}^{\st}]\ar[r]& 0\\
0\ar[r] & \ar@{=}[d] \ar@{^{ (}->}[u] \St_{\CK''}^H\ar[r]^-{b_{\CK''}}&  \ar[d]^{\overline p_{1,\CK''\subset\CK}}  \ar@{^{ (}->}[u]  \BZ[\overline\CT_{1,\CK''}^{\st}]^H\ar[r]^-{\bar\partial_{\CK''}^\st}&\BZ[\CT_{0,\CK''}^{\st}]^H \ar[r] \ar[d]^{p_{0,\CK''\subset\CK}}   \ar@{^{ (}->}[u] & 0\\
0\ar[r] & \St_{\CK''}^H\ar[r]^-{b^{\CK\dash\st}_{\CK''}}& \BZ[\overline\CT_{1,\CK''}^{\CK\dash\st}]^H\ar[r]^-{\bar\partial_{\CK''}^{\CK\dash\st}}&\BZ[\CT_{0,\CK''}^{\CK\dash\st}]^H \ar[r]& 0\\
0\ar[r] & \ar[u]^\simeq\St_{\CK}\ar[r]^-{b_{\CK}}&\ar[u]^\simeq \BZ[\overline\CT_{1,\CK}^{\st}]\ar[r]^-{\bar\partial_{\CK}^\st}&\ar[u]^\simeq\BZ[\CT_{0,\CK}^{\st}]\ar[r]& 0\rlap{.}\\
}\]
\end{compactenum}
\end{Lem}
\begin{proof}
Part (a) for $\CT^{\CK\dash\st}_{0,\CK''}$ and $\CT^{\ori,\CK\dash\st}_{1,\CK''}$ is clear by definition. In the other two cases, it follows from
\[\Stab_G([t,gh]_{\CK''})=\Stab_{gh\CK'' (gh)^{-1}\cap G}(t)=\Stab_{g\CK'' g^{-1}\cap G}(t)=\Stab_G([t,g]_{\CK''}),\]
where in the middle step we use that $\CK''$ is normal in $\CK$. Part (b) is clear from
\[\Stab_G([t,g]_{\CK''})= \Stab_{g\CK'' g^{-1}\cap G}(t)\subset \Stab_{g\CK g^{-1}\cap G}(t)=\Stab_G([t,g]_{\CK}).\]

In part (c), the commutativity of the diagrams formed by the first two and by the last two rows is trivial. The proof of the commutativity for the diagram formed by rows 2 and 3 is analogous to that of \autoref{Lem-ForTwoIGammaN}. The exactness of the second row uses the freeness of the right $H$-action on $\St_{\CK''}$.
It follows from the explicit expressions for $b_{\CK''}$  and $b^{\CK\dash\st}_{\CK''}$, and from our definitions.
\end{proof}
\begin{proof}[Proof of \autoref{prop:IndepOfAdelicIota}]
The proof is analogous to the proof of \autoref{prop-TwoIGammaN}. Let $c$ be in $ C^\ad_\har(N,\CK)^{\GL_2(F)}$. Let $\pi \colon \CT_{1,\CK''}^{\ori}\to \CT_{1,\CK}^{\ori},[e,g]_{\CK''}\mapsto [e,g]_\CK$. We need to show that 
\[\iota_{\CK'',N}( c\circ\pi) = \nu_{\CK,\CK'',N} \circ \iota_{\CK,N}(c)\]
as elements in $\St_{\CK''}\otimes_GN$. We have
\begin{equation}\label{eq-ExprForTop}
\iota_{\CK'',N}( c\circ\pi)= 
\sum_{[e,g]_{\CK''} \in\GL_2(F)\backslash \CT_{\CK'',1}^{\ori,\st}/\{\pm1\}}
[e,g]_{\CK''} \otimes c([e,g]_\CK) \in \BZ[\overline\CT_{1,\CK''}^{\st}]\otimes_GN.
\end{equation}
Because $c$ is invariant under $H$, the expression for $\iota_{\CK'',N}( c\circ\pi)$ lies in $\BZ[\overline\CT_{1,\CK''}^{\st}]^H\otimes_GN$, and also in the kernel of $\bar\partial_{\CK''}^\st \otimes\id_N$; cf.~the second row of the diagram in \autoref{lem-ForAdelicTwoIotas}(c). Next, $ \nu_{\CK,\CK'',N} \circ \iota_{\CK,N}(c)$ is equal to
\begin{equation}\label{eq-ExprForBottom}
\sum_{[e,g]_{\CK} \in\GL_2(F)\backslash \CT_{\CK,1}^{\ori,\st}/\{\pm1\}}
\sum_{h\in H}
[e,gh]_{\CK''} \otimes c([e,g]_\CK) =
\sum_{[e,g]_{\CK''} \in\GL_2(F)\backslash \CT_{\CK'',1}^{\ori,\CK\dash\st}/\{\pm1\}}
[e,g]_{\CK''} \otimes c([e,g]_\CK).
\end{equation}
The right hand side lies in $\BZ[\overline\CT_{1,\CK''}^{\CK\dash\st}]^H\otimes_GN$, and also in the kernel of $\bar\partial_{\CK''}^{\CK\dash\st} \otimes\id_N$; cf.~the third row of the diagram in \autoref{lem-ForAdelicTwoIotas}(c). Moreover $\overline p_{1,\CK''\subset\CK}$ maps the right expression in \eqref{eq-ExprForTop} to the right expression in \eqref{eq-ExprForBottom}. To conclude, we use that the diagram formed by rows 2 and 3 in \autoref{lem-ForAdelicTwoIotas}(c) commutes, and that $\St_\CK\cong(\St_{\CK''})^H$ is a projective and hence flat $\BZ[G]$-module, because $\CK$ is small.
\end{proof}

We need one more preparatory result. For this, let $\supp y:=\{ \Fp\in \Max(A)\mid y_\Fp\notin \GL_2(A_\Fp)\}$ be the support of $y\in  \GL_2(\BA_F^\infty)$, and let $\supp \CK:= \{ \Fp\in \Max(A)\mid \GL_2(A_\Fp)\not\subset \CK\}$ be the support of a compact open subgroup $\CK\subset  \GL_2(\BA_F^\infty)$; in the second definition we embed $\GL_2(F_\Fp)$ into $\GL_2(\BA_F^\infty)$ via the component at $\Fp$. Both sets are finite. One has the following lemma:

\begin{Lem}
\label{lem:HeckeUnderEmbed}
For $\CK\subset \GL_2(\BA_F^\infty)$ a compact open subgroup denote by $X_\CK$ one of the symbols $ S_{k,l}(\CK)$, $C^\ad_\har(N,\CK)^{\GL_2(F)}$, or $\St_\CK\otimes_{\GL_2(F)}N$. For an open normal subgroup $\CK''\subset\CK$ let $\iota_{\CK,\CK''}\colon X_{\CK}\to X_{\CK''}$ be the canonical inclusion, where in the last case we take $\iota_{\CK,\CK''}=\nu_{\CK,\CK'',N}$ from \autoref{cor:AdelicStProps}(b). Let $_{|\CK y\CK}$ in $\End(X_\CK)$ be the Hecke operator attached to some $y\in \GL_2(\BA_F^\infty)$. Then for any non-zero ideal $\Fn\subset A$ such that $(\supp y\cup\supp\CK)\cap\supp\Fn=\emptyset$ and with $\CK'':=\CK\cap\CK(\Fn)$ as in \autoref{lem:KisSmall}, the following diagram commutes
\[
\xymatrix{
X_{\CK}\ar[r]^{\iota_{\CK,\CK''}} \ar[d]_{{}_{|\CK y''\CK}}& X_{\CK''}\ar[d]^{{}_{|\CK'' y''\CK''}}\\
X_{\CK}\ar[r]^{\iota_{\CK,\CK''}} & X_{\CK''}\rlap{.}\\
}
\]
\end{Lem}
\begin{proof}
Set $S:=\supp y\cup\supp\CK$ and $F_S:=\prod_{\Fp\in S} F_\Fp$, and let $\pr_S\colon \GL_2(\BA_F^\infty)\to\GL_2(F_S)$ denote the canonical projection. Define $\CK_S:=\pr_S(\CK)$ and $y_S:=\pr_S(y)$, and let $H:=\prod_{\Fp\in\Max A\setminus S}\GL_2(A_\Fp)$ and  $H(\Fn):=\{ h\in H\mid h\equiv 1\pmod{\Fn} \hbox{ in }\GL_2(A/\Fn)\}$; we regard $\CK_S$, $H$ and $H(\Fn)$ as subgroups of $\GL_2(\BA_F^\infty)$ by extending tuples by $1$ at the missing components. We first show that 
\begin{equation}\label{eq:OnKS}
 \CK = \CK_S \times H.\footnote{A decomposition of the same  form also holds when replacing $S$ by $\supp \CK$.}
\end{equation}
Because $\CK$ is compact open for the adelic topology, there exists a finite set $S'\subset\Max(A)$ such that $\CK\supset  \prod_{\Fp\in\Max A\setminus S'} \GL_2(A_\Fp)$, and by enlarging $S'$, we may assume $S'\supset S$. Let $\Fp\in S'\setminus S$. Then, by the definition of $\supp\CK$ and because $S\supset\supp\CK$, we have $\CK\supset\GL_2(A_\Fp)$. Thus, by taking a finite product of subgroups, it follows that
\[
\CK\supset  \prod_{\Fp\in S'\setminus S} \GL_2(A_\Fp)\times  \prod_{\Fp\in\Max A\setminus S'} \GL_2(A_\Fp)=H.
\]
Because $H$ is maximal compact in $\prod_{\Fp\in\Max A\setminus S}\GL_2(F_\Fp)$, any $y'\in \CK$ can be written as a product $y'=y'_S\cdot h$ with $y'\in\CK_S$ and $h\in H$. By using that in fact $H\subset\CK$, we obtain \eqref{eq:OnKS}. Intersecting  \eqref{eq:OnKS} with $\CK(\Fn)$ also gives $ \CK'' = \CK_S \times H(\Fn)$; we refer to this formula as $(16'')$.

Choose now $y_{S,i}\in\GL_2(F_S)$ such that $\CK_S y_S\CK_S =\coprod_i y_{S,i}\CK_S$, using that the left hand side is compact and right translation invariant under the compact open subgroup $\CK_S\subset\GL_2(F_S)$. For $g\in \GL_2(\BA_F^\infty)$, in the following we also write $g=(g_S,g^S)$ with $g_S=\pr_S(g)$ and $g^S$ the tuple of those components of $g$ with index set $\Max A\setminus S$; we use this notation also for subsets whenever this is permitted.

Define $y_i\in \GL_2(\BA_F^\infty)$ as the extension of $y_{S,i}$ given by $y_i=(y_{S,i}, y^S)$. Then we have
\[ \CK'' y\CK''\stackrel{(16'')}= (\CK_S y_S\CK_S,H(\Fn) y^S H(\Fn))= (\coprod_i y_{S,i}\CK_S,y^S H(\Fn))= \coprod_i (y_{S,i},y^S)(\CK_S,H(\Fn))= \coprod_i y_i\CK'', \]
because $y^S\in H$ and $H(\Fn)$ is normal in $H$. The analogous computation for $\CK$ in place of $\CK''$ gives $\CK y''\CK =\coprod_i y_{i}\CK$. This shows that  the Hecke operators ${}_{|\CK'' y''\CK''}$ and ${}_{|\CK y''\CK}$ can be defined by the same expressions and now the commutativity is clear.
\end{proof}

\begin{Cor}
\label{thm:CHar-St-global}
Let $\CK$ be any compact open subgroup of $\GL_2(\BA_F^\infty)$. Then  the following hold:
\begin{compactenum}
\item 
\label{thm:CHar-St-globala}
The map $\iota_{\CK,N}$ defines a Hecke-equivariant isomorphism
\[C^\ad_\har(N,\CK)^{\GL_2(F)} \stackrel\simeq\longto \St_{\CK}\otimes_{\GL_2(F)}N.\]
This holds in particular for $N=V_{k,l}$.
\item 
\label{thm:CHar-St-globalb}
For any short exact sequence $0\to N'\to N\to N''\to 0$ of $F[\GL_2(F)]$-vector spaces that are finite dimensional over $F$ one has a short exact sequence of Hecke-modules
\[0\longto C^\ad_\har(N',\CK)^{\GL_2(F)}\longto C^\ad_\har(N,\CK)^{\GL_2(F)}\longto C^\ad_\har(N'',\CK)^{\GL_2(F)}\longto 0.\]
\end{compactenum}
\end{Cor}

\begin{proof}
Part (b) follows from part (a) and from \autoref{cor:AdelicStProps}(c) and (d). To prove (a), let $y\in\GL_2(\BA_F^\infty)$. Choose a non-zero proper ideal $\Fn\subset A$ such that $(\supp y\cup\supp\CK)\cap\supp\Fn=\emptyset$. Consider the commutative diagram
\[
\xymatrix@C+1pc{
C^\ad_\har(N,\CK)^{\GL_2(F)} \ar[r]^-{\iota_{\CK,N}} \ar[d]^{\iota_{\CK,\CK''}}&\St_{\CK}\otimes_{\GL_2(F)}N \ar[d]^{\iota_{\CK,\CK''}} \\
C^\ad_\har(N,\CK'')^{\GL_2(F)} \ar[r]^-{\iota_{\CK'',N}} &\St_{\CK''}\otimes_{\GL_2(F)}N\rlap{.}\\
}
\]
By \autoref{lem:HeckeUnderEmbed} the vertical maps are equivariant for the Hecke operator $|_{\CK y\CK}$. By \autoref{thm:CHar-St-global-Small} the lower horizontal map has the same property. Since the vertical maps are inclusions and the horizontal maps are isomorphisms, the result follows.
\end{proof}

\begin{Rem}
Suppose that $\Cl_\CK$ is the trivial group (which for instance happens for $A=\BF_q[t]$ and $\CK$ maximal compact). Let $\Gamma:= \GL_2(F)\cap\CK$.
Then one can define meaningful Hecke actions directly on $C_\har(V_{k,l})^\Gamma$ and on $\St \otimes_\Gamma V_{k,l}$, so that $\iota_{\Gamma,V_{k,l}}$ becomes a Hecke equivariant isomorphism; see \cite[\S~6]{Gekeler} for the case $S_{k,l}(\Gamma)$, or \cite[Ch.~6]{Boeckle-EiSh} (note that the actions are normalized differently, which we will discuss further at the end of \autoref{RelBosPel}). In this particular case, there is no gain of the adelic over the local situation. The Hecke actions agree and proofs from either perspective are of a similar difficulty. However, we should point out that, by definition, the Hecke action on $\St_\CK\otimes_{\GL_2(F)}V_{k,l}$ is induced from an action on $\St_\CK$. In the local situation this is not true anymore: The Hecke action involves the module $V_{k,l}$ directly. Thus to show the Hecke equivariance for morphisms under the functor in \autoref{thm:CHar-St-local-Gen}(b) the full $\GL_2(F)$-action on the coefficients needs to be considered.
\end{Rem}

\section{Representations of $\SL_2(F)$ and $\GL_2(F)$}
By the results of the previous section, in particular \autoref{thm:CHar-St-global}, we can use representation theory to gain insights into the Hecke-module structure of spaces of Drinfeld cusp forms. In this section, we first develop the relevant representation theory of $\SL_2(F)$ and $\GL_2(F)$ building heavily on \cite{Bonnafe}. Afterwards, we introduce certain explicit maps between symmetric powers which will have very natural interpretations when passing to Drinfeld cusp forms in \autoref{Section5}. Finally, we also briefly mention representations of the finite groups $\SL_2(\BF_q)$ and $\GL_2(\BF_q)$ which will become relevant in \autoref{Section6}.

\label{Section4}
\subsection{Symmetric powers and irreducible representations} 
\label{Subsec3.1}
This subsection is mostly standard material. For an excellent reference in the case of $\SL_2$ we refer to \cite{Bonnafe}. Denote by $F[X,Y]_k$ the subspace of homogeneous polynomials of degree $k$ in $F[X,Y]$.  We let the group $\GL_2(F)$ act on $F[X,Y]$ by 
\[\SMat{a}{b}{c}{d}\cdot X^iY^j = (ad-bc)^{-i-j}(dX-bY)^i(-cX+aY)^j.\]
This action clearly preserves the homogeneous components $F[X,Y]_k$. If one identifies  $F[X,Y]_1$ with the dual vector space of $F^2$ in the obvious way, then $\Sym^k ((F^2)^*)$ is naturally isomorphic to $F[X,Y]_k$ as a $\GL_2(F)$-module. For shorter notation, we often abbreviate $\Delta_k:=\Sym^k ((F^2)^*)$. For any $k\ge0$, we define $L_k\subset \Delta_k$ as the $F$-linear span of $\{X^iY^{k-i}\mid \binom ki\not\equiv0\pmod p\}$. In \autoref{Lem-BasicsOnLk} we shall show that $L_k$ is an irreducible subrepresentation of $\Delta_k$.

To see that $L_k$ is a subrepresentation of $\Delta_k$, and to introduce the Frobenius twist of a representation, note that 
\[\tau \colon\GL_2(F)\to \GL_2(F),\SMat{a}{b}{c}{d}\to\SMat{a^p}{b^p}{c^p}{d^p}\]
is a group homomorphism. Given any representation $V$ of $\GL_2(F)$ on an $F$-vector space, it follows that it $s$-fold Frobenius twist $V^{(s)}$, defined as the vector space $V$ together with the action $\GL_2(F)\times V\to V,(\gamma,v)\mapsto \tau^s(\gamma)v$, is again a $\GL_2(F)$-representation. One verifies that
\[F[X,Y]_k^{(s)}\to F[X,Y]_{p^sk},f(X,Y)\mapsto f(X^{p^s},Y^{p^s})\]
is a monomorphism of $\GL_2(F)$-representations. By a similar computation one also verifies the following: Write $k=k_0+k_1p+\ldots+k_sp^s$ in its base $p$ expansion, with $0\le k_i\le p-1$, where we insist that $k_s\neq0$ -- we also use the notation $k=(k_s\,\ldots\,k_1\,k_0)_p$ to denote the base $p$-expansion of $k$. Then one has a monomorphism of $\GL_2(F)$-representations
\begin{equation}\label{Eqn-Steinberg}
\otimes_{i=0}^s \Delta_{k_i}^{(i)}\to \Delta_k,f_0(X,Y)\otimes\ldots \otimes f_s(X,Y)\mapsto \prod_{i=0}^s f_i(X^{p^i},Y^{p^i})
\end{equation}
with image $L_k$. This is a form of the Steinberg Tensor Theorem specialized to the group $\GL_2$; for a direct proof see \cite[Lemma 8]{Pellarin-CharPReps}. It follows immediately that
\[\dim_F L_k = \prod_{i = 0}^s(k_i+1). \] 
In the sequel, for $m\in\BZ$ we denote by $\det^m$ the action of $\GL_2(F)\times F\to F,(\gamma,v)\mapsto \det\gamma^m v$.

Let us recall results on (highest) weights for $\GL_2$. Denote by $T(F)$ the subtorus $\big\{\SMat{a}00d\mid a,d\in F^\times\big\}$ of $\GL_2(F)$. It is commutative, and so the restriction of any $F$-linear algebraic representation $V$ of $\GL_2(F)$  to $T(F)$ decomposes as a direct sum $V=\oplus_{(n,m)\in\BZ^2}V(n,m)$ of $F$-subvector spaces where $\SMat{a}00d$ acts on $v\in V(n,m)$ as $a^nd^m$. Only finitely many $V(n,m)$ are non-zero, and the corresponding pairs $(n,m)$ are called the weights of $V$. From $\SMat0110\in\GL_2(F)$ one deduces that with $(n,m)$ also $(m,n)$ is a weight. A weight $(n,m)$ is called \emph{dominant} if $n\geq m$. We define a partial order on the weights of $V$ as follows:
\[
(n,m)\geq(n',m') \quad \text{if and only if} \quad (n-n',m-m')=l(1,-1) \text{ for some } l\geq 0.
\]
A weight $(n,m)$ is called \emph{highest weight} of $V$ if for all other weights $(n',m')$ of $V$ one has $(n,m)\geq(n',m')$. 
An important theorem for connected reductive groups, specialized to $\GL_2$, asserts that each irreducible representation has a unique highest weight and that for each dominant weight $(n,m)$ there is a unique irreducible representation of $\GL_2$ of that highest weight; see \cite[Part II, Prop.~2.4]{Jantzen-Representations}.

\begin{Lem}\label{Lem-BasicsOnLk}
Let $k\ge0$. The following hold:
\begin{compactenum}
\item The representation $L_k$ is irreducible. Its weights as a representation of $\GL_2(F)$ are the pairs $(-i,i-k)$ with $0\le i\le k$ and $\binom ki\not\equiv 0\pmod p$, and they all occur with multiplicity~$1$.
\item $L_k$ is the unique simple subrepresentation of $\Delta_k$. In particular, it is the socle of $\Delta_k$.
\item The irreducible representations of $\GL_2(F)$ up to isomorphism are the representations $L_{k}\otimes\det^m$ for integers $k,m\in\BZ$, $k\geq 0$.
\item One has $L^*_k\otimes \det^{-k}\cong L_k$.
\item One has $L_k=\Delta_k$ if and only if all but possibly the leading digit $k_s$ in the base $p$-expansion $k=(k_s\,\ldots\,k_1\,k_0)_p$ of $k$ are equal to $p-1$. \footnote{These are precisely D. Goss's so-called `$p$-Magic Numbers' which arise in connection with many objects in function field arithmetic, e.g. Goss's special $\Gamma$ and $\zeta$ values. One calls $k$ {\it even} when $k_s = p-1$ and {\it odd} otherwise. Here we see $p$-Magic Numbers characterized as precisely those positive integers $k$ for which ${k \choose i} \not\equiv 0 \pmod{p}$, for all $0 \leq i \leq k$. }
\end{compactenum}
\end{Lem}

We first recall some properties of binomial coefficients, that are also used in other places. For $m\in\BZ$ and $i\ge0$, they are defined by $\binom{m}i=\frac{m(m-1)\cdot\ldots\cdot(m-i+1)}{i!}\in\BZ$ with the convention that empty products equal $1$, so that $\binom{m}0=1$. To link negative to positive upper indices one has
\begin{equation}\label{eq-BinomReflection}
\binom{m}{i}=(-1)^i\binom{i-1-m}{i},
\end{equation}
which is easily verified. The defining formula also shows that $m\mapsto\binom{m}i$ extends to a uniformly continuous function $\binom{\cdot}i\colon \BZ_p\to\BZ_p$ for the $p$-adic topology. Moreover one has Lucas'~Theorem. 
\begin{Prop}[Lucas]\label{Prop-Lucas}
Let $i=i_0+i_1p+\ldots+i_dp^d$ and $m=\sum_{j\ge0} m_jp^j$ be the base $p$ expansions of $i\in\BN_0$ and $m\in\BZ_p$ with $i_j,m_j\in\{0,\ldots,p-1\}$. Then 
\[\binom{m}i\equiv \prod_{j=0}^d \binom{m_j}{i_j}\pmod p.\]
In particular, $\binom{m}i\pmod p$ is zero if and only if $m_j<i_j$ for some $j\in\{0,\ldots,d\}$, and moreover the map $m\mapsto \binom{m}i\pmod p$ is $p^{d+1}$ periodic.
\end{Prop}
\begin{proof}
For $m\in\BN_0$ this is classical, and can be found for instance in \cite{Granville}. For $m\in\BZ_p$, the formula follows from the uniform continuity of $\binom{\cdot}i$ and the congruence for all $m\in\BN_0$.
\end{proof}

\begin{proof}[Proof of \autoref{Lem-BasicsOnLk}:] We first prove (a) and (b). From the definition of $L_k$ it is straightforward to see that its weights are the pairs $(-i,i-k)$ with $0\le i\le k$ and $\binom ki\not\equiv 0\pmod p$, and for these $i$ each $L_k(-i,k-i)$ has dimension $1$. In particular, $L_k$ has highest weight $(0,-k)$. It is shown in  \cite[Thm.~10.1.8.(b)]{Bonnafe} that the remaining assertions are true for the restrictions of $L_k$ and $\Delta_k$ to $\SL_2(F)$. Since we know that $L_k$ is a subrepresentation of $\Delta_k$, knowing (a) and (b) for $\SL_2$ implies the analogous assertions for $\GL_2$.

Part (c) follows from (a) and the classification of irreducible representations in terms of highest weights stated above. For (d) note that the dual of an irreducible representation is again irreducible, since the canonical bidual-map is an isomorphism, and thus $L_k^*$ is irreducible. To determine its weights, let $\xi_i$ be the unique element of $L_k^*$ that is $1$ on $X^iY^{k-i}$ and zero on $X^{i'}Y^{k-i'}$ for $i'\neq i$. These elements form a basis of $L_k^*$ and 
\[(\SMat a00d \xi_i)(X^{i'}Y^{k-i'})=\xi_i(\SMat a00d^{-1} X^{i'}Y^{k-i'})=(ad)^k\xi_i(a^{i'-k}d^{-i'} X^{i'}Y^{k-i'})=a^{i}d^{k-i} \xi_i( X^{i'}Y^{k-i'}).\]
Twisting with $\det^{-k}$ preserves the irreducibility, but changes the weight $(i,k-i)$ to $(i-k,-i)$. Hence $L_k^*\otimes \det^{-k}$ has the same weights as $L_k$ and is thus isomorphic to it because both are irreducible. Part (e) is immediate from \autoref{Prop-Lucas} and the definition of $L_k$. One can also use (\ref{Eqn-Steinberg}) and a dimension count.
\end{proof}

For the algebraic group $G\in\{\SL_2,\GL_2\}$ let $K_0(G)$ denote the Grothendieck ring of algebraic representations of $G(F)$ on finite-dimensional $F$-vector spaces. The multiplication of the ring is given by the tensor product of representations. If $G=\SL_2$, then as an additive group $K_0(\SL_2)$ is the free module on the symbols $L_k$, $k\ge0$ (we regard $L_k$ as a module for $\SL_2\subset \GL_2$). By \cite[(10.1.14)]{Bonnafe} $\Delta_k$, $k\ge0$, is also a $\BZ$-basis of $K_0(\SL_2)$. If $G=\GL_2$, then from \autoref{Lem-BasicsOnLk} it easily follows that the $L_k\otimes\det^m$, $k,m\in\BZ$, $k\geq 0$ form a $\BZ$-basis of $K_0(\GL_2)$, and below we shall deduce rapidly from Bonnaf\'e's results that the symbols $\Delta_k\otimes\det^m$, $k\ge0$, $m\in\BZ$ also form a $\BZ$-basis of $K_0(\GL_2)$. It follows from \autoref{Lem-BasicsOnLk}(a) that when expressing $\Delta_k\otimes \det^m$ in $K_0(\GL_2)$ as a linear combination in the $L_{k'}\otimes\det^{m'}$ then only multiplicities $0$ and $1$ can occur, that $L_k\otimes\det^m$ occurs with multiplicity $1$, and that $L_{k'}\otimes\det^{m'}$ cannot occur in $\Delta_k\otimes\det^m$ if $k'>k$. 

To give a recursive formulas for the coefficients in $\{0,1\}$ that occur when expressing $\Delta_k\otimes \det^m$ in $K_0(\GL_2)$ as a linear combination in the $L_{k'}\otimes\det^{m'}$, we follow \cite[Subsec.~10.1.3]{Bonnafe}: Define recursively sets $\CE(k)$, $k\ge0$ as follows, where by $k_0$ we denote the lowest digit of $k$ in its base $p$ expansion, i.e., the residue of $k$ by division by $p$: 
\[\CE(k)=\left\{\begin{array}{ll}
\{0\}&\hbox{if $0\le k\le p-1$,}\\
p\CE\big(\frac{k-k_0}{p}\big)&\hbox{if $k\ge p$ and $k_0=p-1$,}\\
p\CE\big(\frac{k-k_0}{p}\big)\sqcup \big( k_0+1+p\CE(\frac{k-k_0-p}{p})\big)&\hbox{if $k\ge p$ and $k_0\leq p-2$.}
\end{array}
\right.\]
\begin{Prop}[{\cite[Proposition 10.1.18]{Bonnafe}}]\label{Prop-V_k-versus-L_k}
In $K_0(\SL_2)$ for $k\ge0$ one has
\[\Delta_k=\bigoplus_{k'\in k-2\CE(k)} L_{k'}.\]
\end{Prop} 
\begin{Cor}\label{Cor-V_k-versus-L_k}
In $K_0(\GL_2)$ for $k\ge0$ and $m\in\BZ$ one has
\[\Delta_k\otimes\det^m=\bigoplus_{k'\in k-2\CE(k)} L_{k'} \otimes\det^{(k-k')/2+m}.\]
\end{Cor}
\begin{proof}
To prove the result, by twisting with powers of $\det$ we may clearly assume $m=0$. By the basis property of the $L_k\otimes\det^m$, given $k$, there exist unique constants $e_{k',m'}\in\BN_0$ such that in $K_0(\GL_2)$ we have
\[\Delta_k=\bigoplus_{k',m'} (L_{k'} \otimes\det^{m'})^{\oplus e_{k',m'}}.\]
This equality can be restricted to $K_0(\SL_2)$. It follows from \autoref{Prop-V_k-versus-L_k} that 
\[ \sum_{m'\in\BZ} e_{k',m'} =\left\{\begin{array}{ll}
1,& \hbox{if }k' \in k-2\CE(k),\\
0,& \hbox{if }k' \not\in k-2\CE(k).
\end{array}
\right.\]
Hence for each $k'\in k-2\CE(k)$ there is a unique $m'$ for which $e_{k',m'}$ is non-zero, and equal to $1$, and for all other $k'$, all $e_{k',m'}$ are zero. It remains to determine the unique $m'$ for the former $k'$. For this note that $\SMat a00a$, $a\in F^\times$ acts on $\Delta_{k'}\otimes\det^{m'}$ as $a^{k'+2m'}$ and on $\Delta_k$ as $a^k$. This implies the claim on $m'$ in the corollary, and we are done.
\end{proof}

\subsection{Hyperderivatives on symmetric powers}
\label{Subsec3.2}
Motivated by results of Bosser--Pellarin on certain hyperderivatives between certain spaces of Drinfeld modular forms, we were seeking a representation theoretic description of these maps in via the residue map. This is how we arrived at the definitions in this section. The detailed relation to the work of Bosser--Pellarin will be explained \autoref{Sec-HyperDer}. 

We need the following result.
\begin{Lem}\label{Lem-OnBinomials}
Let $s\ge1$ and let $r=1+\lfloor \log_p s \rfloor \in\BN$.
\begin{compactenum} 
\item For all $a\in\BZ$ and $i\ge0$ we have $\binom{i+ap^{r}}s\equiv\binom{i}s\pmod p$. 
\item For $k\ge2$ the following conditions are equivalent:
\begin{compactenum}
\item For $l=1,\ldots,s$ one has $\binom{k+s-1}l\equiv0\pmod p$.
\item The lowest $r$ digits of $k+s-1$ in its base $p$ expansion are zero.
\item For $i=0,\ldots,k+2s-2$ and $j=k+2s-2-i$ we have $\binom{j}s\equiv(-1)^s\binom{i}s\pmod p$.
\end{compactenum}
\end{compactenum}
\end{Lem}

\begin{proof}
Part (a) is an immediate consequence of the last assertion of \autoref{Prop-Lucas}.

In (b), the equivalence of (i) and (ii) is immediate from Lucas' theorem \autoref{Prop-Lucas} applied with $m=k+s-1$. Next we prove (ii) $\Rightarrow$ (iii). By (ii), the number $k+s-1$ is divisible by $p^{r}$, and (iii) follows from
\begin{equation}\label{Eqn-BinomEqn}
\binom{j}s\equiv\binom{(k+s-1)+(s-1-i)}s\stackrel{\mathrm{(a)}}\equiv\binom{s-1-i}s \stackrel{\eqref{eq-BinomReflection}}\equiv(-1)^s\binom{i}s\pmod p.
\end{equation}
Suppose conversely that (iii) holds. Let $m$ be the remainder of $-(k+s-1)$ in $\{0,\ldots,p^r-1\}$ under division by $p^r$. Then for all $i=0,\ldots,k+2s-2$ and in particular for $i=0,\ldots,s-1$, we deduce
\[\binom{i+m}s  \stackrel{\mathrm{(a)}}\equiv \binom{i-(k+s-1)}s \stackrel{\eqref{eq-BinomReflection}}\equiv  (-1)^s \binom{(k+s-1)+(s-1-i)}s  \stackrel{(iii)}\equiv \binom{i}s\pmod p.\]
 Assume that $0<m<p^r$. Let $m=m_0+m_1p+\ldots+m_{r-1}p^{r-1}$ and $s=s_0+s_1p+\ldots+s_{r-1}p^{r-1}$ and define $i_t=\max\{0,s_t-m_t\}$ for $t=0,\ldots,r-1$, and $i=i_0+i_1p+\ldots+i_{r-1}p^{r-1}$. Because $m>0$ we have $0\le i \le s-1$ and by Lucas' theorem we deduce $\binom{i+m}s\not\equiv0\pmod p$ while $\binom is\equiv0\pmod p$, which is a contradiction. Hence we must have $m=0$ and hence (ii) holds.
\end{proof}

\begin{Prop}\label{Prop-DsAsMapOnReps} Let $k\geq 2$, $s\geq 1$, $m\in \BZ$.
Suppose  that $\binom{k+s-1}i=0\pmod p$ for $i=1,\ldots,s$. Define
\[ \mathcal{D}_s\colon \Delta_{k-2+2s}\otimes \det^{m+s}\to \Delta_{k-2}\otimes \det^{m} \]
on the basis $X^iY^j$ by $\mathcal{D}_s X^iY^j=(-1)^s\binom{i}s X^{i-s} Y^{j-s}$ and $F$-linearly extended. Then the following hold:
\begin{compactenum}
\item The map $\mathcal{D}_s$ is well-defined, i.e., $\mathcal{D}_s X^iY^j=0$ for $i<s$ or $j<s$.
\item The map $\mathcal{D}_s$ is $\GL_2(F)$-equivariant.
\item The kernel of $\mathcal{D}_s$ is the $F$-linear span of $B:=\{X^iY^j\mid i+j=k+2s-2,\binom{i}s=0\pmod p\}$.
\end{compactenum}
In particular, the $F$-linear hull $W$ of $B$ in $\Delta_{k-2+2s}\otimes \det^{m+s}$ is stable under $\GL_2(F)$.
\end{Prop}
\begin{proof}
To see (a) observe that if $i<s$, then $\binom{i}s=0\pmod p$, and thus $\mathcal{D}_s X^iY^j=0$ for $i<s$. Moreover if $j<s$, then $\binom{i}s=(-1)^s\binom{j}s=0\pmod p$ by \autoref{Lem-OnBinomials}, and it follows that $\mathcal{D}_s X^iY^j=0$ also for~$j<s$.

To prove (b) it suffices to verify it for matrices $\gamma\in\GL_2(F)$ of the form form (i) $\SMat{1}{b}{0}{1}$, $b\in F$, (ii) $\SMat{1}{0}{c}{1}$, $c\in F$, and (iii) diagonal matrices $\SMat{a}{0}{0}{1}$, $a\in F^\times$, since matrices of this form generate $\GL_2(F)$. We denote the actions of $\gamma$ on $\Delta_{k-2+2s}\otimes \det^{m-s}$ and $\Delta_{k-2}\otimes \det^{m}$ by $\bullet $ and $\cdot$, respectively. We begin with (iii) and so $\gamma=\SMat{a}{0}{0}{1}$.
\begin{eqnarray*}
\mathcal{D}_s(\gamma\bullet X^iY^j)&=&\mathcal{D}_s(a^{-k+2-2s+j+m+s}X^iY^j) \ = \ (-1)^s a^{-k+2+j-s+m}\binom{i}s X^{i-s}Y^{j-s} \\ & = & (-1)^s\binom{i}s \gamma\cdot  X^{i-s}Y^{j-s}
\ = \ \gamma\cdot \mathcal{D}_s(X^iY^j).
\end{eqnarray*}
Next let $ \gamma=\SMat{1}{b}{0}{1}$. Then the assertion follows from
\begin{eqnarray*}
\mathcal{D}_s(\gamma\bullet X^iY^j)&=&\mathcal{D}_s((X-bY)^iY^j) \ = \ \sum_{m=0}^i\binom{i}{m} (-b)^{i-m}\mathcal{D}_s(X^{m}Y^{i+j-m}) \\
&=&\sum_{m=0}^i\binom{i}{m}(-b)^{i-m}(-1)^s\binom{m}s X^{m-s}Y^{k+2s-2-m-s} \\
&=&(-1)^s\sum_{m=s}^i(-b)^{i-m}\binom{i}{m} \binom{m}s X^{m-s}Y^{k+s-2-m} \\
&=&(-1)^s\sum_{m=0}^{i-s}(-b)^{i-m-s}\binom{i}{m+s} \binom{m+s}s X^{m}Y^{k-2-m} \\
&=&(-1)^s\sum_{m=0}^{i-s}(-b)^{i-m-s}\binom{i}s \binom{i-s}m X^{m}Y^{k-2-m}, 
\end{eqnarray*}
and
\begin{eqnarray*}
\gamma\cdot \mathcal{D}_s(X^iY^j)&=&\gamma\cdot (-1)^s\binom{i}s X^{i-s}Y^{j-s} \ = \ (-1)^s\binom{i}s (X-bY)^{i-s}Y^{j-s}
\\
&=&(-1)^s\binom{i}s \sum_{m=0}^{i-s}\binom{i-s}{m} (-b)^{i-s-m} X^{m}Y^{i-s-m+j-s} \\
&=&(-1)^s\binom{i}s \sum_{m=0}^{i-s}\binom{i-s}{m} (-b)^{i-s-m} X^{m}Y^{k-2-m}. 
\end{eqnarray*}

Finally the case $ \gamma=\SMat{1}{0}{c}{1}$ is completely analogous to the previous case using \autoref{Lem-OnBinomials}(b).

Regarding (c), note first that $B$ lies clearly in the kernel of $\mathcal{D}_s$. However it is also obvious that set $\{X^{i-s}Y^{j-s}\mid \binom{i}s\neq 0\pmod p,0\le i,j, i+j=k+2s-2\}$ is $F$-linearly independent. This completes the proof of (c). Finally, the last sentence is immediate from~(c).
\end{proof}

\begin{Rem}
While our main motivation to study these particular maps is tied to the results in \autoref{Sec-HyperDer}, they are already of purely representation-theoretic interest. Note however that the image of $\mathcal{D}_s$ can very well be (and is) $0$ in many cases (for example when the source is already irreducible as a  $\GL_2(F)$-representation.
\end{Rem}
\begin{Rem}
For a fixed $k\geq 0$, let $B_k$ denote the intersection of all kernels of the hyperderivatives with source $\Delta_k$. By \autoref{Lem-BasicsOnLk}(b) we have that $L_k\subset B_k$ once $B_k$ is non-trivial. It is a natural question to ask if equality holds. The following example shows that this is in general not the case: Let $q=3$ and $k=40$. Then, there are hyperderivatives only for $s_1=2$, $s_2=5$ and $s_3=14$. Now, $X^{16}Y^{24}$ is in $B_{40}$, since $\binom{16}{2}\equiv 0 \pmod 3$, $\binom{16}{5}\equiv 0 \pmod 3$ and $\binom{16}{14}\equiv 0 \pmod 3$ by \autoref{Prop-Lucas}. However, since $\binom{40}{16}\equiv 0 \pmod 3$, $X^{16}Y^{24}$ is not in $L_{40}$.
\end{Rem}

\subsection{The Cartier operator}
\label{Cartier}
In the spirt of the previous subsection, we were looking for a representation theoretic description of the Frobenius operator on Drinfeld modular forms. Via the residue map we arrived at the following construction. Since contrary to the previous case, our new map is only defined over a perfect field, we work with $\BC_\infty$ for simplicity of notation, and we identify $\Delta_k\otimes_F\BC_\infty$ with $\BC_\infty[X,Y]_k$. 
Let $\sigma$ denote the Frobenius of the perfect field $\BC_\infty$. Then, for any $\BC_\infty$-vector space $M$, we denote by $\sigma_\ast M$ the $\BC_\infty$-vector space with underlying abelian group $M$ and scalar multiplication given by
\[
a\cdot m=\sigma(a) m=a^{p} m \quad \text{for } a\in\BC_\infty, m\in M.
\]
If $M$ is a $\BC_\infty[\GL_2(F)]$-module, so is $\sigma_\ast M$.
\begin{Prop}
\label{Prop-Cartier} Let $k\geq 2$ and $m\in \BZ$.
The map $\CC_p\colon\sigma_*((\Delta_{pk-2}\otimes\det^{pm-1})\otimes_F\BC_\infty)\rightarrow (\Delta_{k-2}\otimes \det^{m-1})\otimes_F\BC_\infty$ given by
\[
\CC_p(aX^{i-1}Y^{j-1})=\left\{\begin{aligned}
&a^{\frac{1}{p}}X^{\frac{i}{p}-1}Y^{\frac{j}{p}-1}, && \text{if } p\mid i,\\
&0, && \text{otherwise,}\\
\end{aligned}
\right.
\]
is $\GL_2(F)$-equivariant and surjective.
\end{Prop}
\begin{proof}
Firstly, the map is $\BC_\infty$-linear by construction and clearly surjective. We proceed as in the proof of \autoref{Prop-DsAsMapOnReps} and consider matrices $\gamma\in\GL_2(F)$ of the form (i) $\SMat{1}{b}{0}{1}$, $b\in F$, (ii) $\SMat{1}{0}{c}{1}$, $c\in F$, and (iii) diagonal matrices $\SMat{a}{0}{0}{1}$, $a\in F^\times$. 

We begin with (iii) and so $\gamma=\SMat{a}{0}{0}{1}$. Here, we only need to consider the case $p\mid i$ (and so $p\mid j$).
\begin{eqnarray*}
\mathcal{C}_p(\gamma\cdot X^{i-1}Y^{j-1})&=&\mathcal{C}_p(a^{-pk+j+pm}X^{i-1}Y^{j-1}) \ = \ 
a^{-k+\frac{j}{p}+m}X^{\frac{i}{p}-1}Y^{\frac{j}{p}-1}
\ = \ \gamma\cdot \mathcal{C}_p(X^{i-1}Y^{j-1})
\end{eqnarray*}
Next let $ \gamma=\SMat{1}{b}{0}{1}$. Then we have  
\begin{eqnarray*}
\mathcal{C}_p(\gamma\cdot X^{i-1}Y^{j-1})&=&\mathcal{C}_p((X-bY)^{i-1}Y^{j-1})\\
&=&\CC_p\left(\sum_{n=0}^{i-1}\binom{i-1}{n}(-b)^nX^{i-1-n}Y^{n+j-1}\right)\\
&=&\sum_{n=0}^{i-1}\binom{i-1}{n}(-b)^{\frac{n}{p}}\CC_p(X^{i-1-n}Y^{n+j-1}).
\end{eqnarray*}
In this sum, only the terms with $i-n= 0 \pmod p$ can possibly be non-zero by definition of $\CC_p$. Let $i=i_0+pi_1+\ldots+p^{d}i_d$ be the base $p$ expansion of $i$. Then, if $i-n= 0 \pmod p$, the base $p$ expansion of $n$ is of the form $n=i_0+pn_1+\ldots+p^{d}n_d$. We first consider the case $p\nmid i$, i.e. $i_0\neq 0$. Then, by Lucas's theorem we obtain
\[
\binom{i-1}{n}=\binom{i_0-1}{i_0}\prod_{j=1}^{d}\binom{i_j}{n_j}=0 \pmod p.
\]
Consequently, $\mathcal{C}_p(\gamma\cdot X^{i-1}Y^{j-1})=0=\gamma\cdot\mathcal{C}_p( X^{i-1}Y^{j-1})$. In the other case $p\mid i$, i.e. $i_0=0$, we have $i-1=p(\frac{i}{p}-1)+p-1$. Let $m=\frac{n}{p}\in\BZ$. Then by Lucas's theorem
\[
\binom{i-1}{n}=\binom{p-1}{0}\binom{\frac{i}{p}-1}{m}=\binom{\frac{i}{p}-1}{m} \pmod p.
\]
We obtain
\begin{eqnarray*}
\mathcal{C}_p(\gamma\cdot X^{i-1}Y^{j-1})&=&\sum_{m=0}^{\frac{i}{p}-1}\binom{\frac{i}{p}-1}{m}(-b)^{m}X^{\frac{i}{p}-1-m}Y^{m+\frac{j}{p}-1}
=\gamma\cdot X^{\frac{i}{p}-1}Y^{\frac{j}{p}-1}= \gamma\cdot \CC_p(X^{i-1}Y^{j-1}).
\end{eqnarray*}
The case $ \gamma=\SMat{1}{0}{c}{1}$ follows in a completely similar way, thus completing the proof.
\end{proof}

\begin{Rem}
For $m=0$ the map in \autoref{Prop-Cartier}  is in fact a special case of the so called \emph{Cartier operator}: Considering the natural action of $\GL_2(\BC_\infty)$ on the affine plane $\BA^2_{\BC_\infty}$, for the corresponding module of differentials $\Omega_{\BC_\infty[X,Y]/\BC_\infty}$ one has the $\GL_2(\BC_\infty)$-equivariant isomorphism 
\[\BC_\infty[X,Y]\otimes\det \longto \omega_{\BC_\infty[X,Y]/\BC_\infty}, g(X,Y)\mapsto g(X,Y)\mathrm{d}X\wedge\mathrm{d}Y\]
with $ \omega_{\BC_\infty[X,Y]/\BC_\infty}=\bigwedge^2 \Omega_{\BC_\infty[X,Y]/\BC_\infty}$. The Cartier operator is a very general map $\sigma_*\omega_{Z/K}\to \omega_{Z/K}$ for $K$ perfect of characteristic $p>0$ and $Z/K$ a smooth variety. In the case at hand, it is the $\BC_\infty[X,Y]$-linear map $ \sigma_*\omega_{\BC_\infty[X,Y]/\BC_\infty}\to  \omega_{\BC_\infty[X,Y]/\BC_\infty} $ characterized by having $\frac{\mathrm{d}X\wedge\mathrm{d}Y}{XY}$ as a fixed point. In this perspective, another possible proof of the $\GL_2(F)$-equivariance of the map $\CC_p$ is to use $\BC_\infty[X,Y]=\BC_\infty[\alpha X+\beta Y,\gamma X+\delta Y]$ for $\SMat{\alpha}{\beta}{\gamma}{\delta}\in\GL_2(\BC_\infty)$, i.e., to apply an automorphism of $\BA^2_{\BC_\infty}$.
\end{Rem}

\begin{Rem}
Using \autoref{Prop-Lucas} it is easy to show that $\mathcal{C}_p(\sigma_*((L_{pk-2}\otimes\det^{pm+1})\otimes_F\BC_\infty))=0$.
\end{Rem}

\begin{Lem}
\label{Lem-DualFrob}
Let $M$ be a $\BC_\infty[\GL_2(F)]$-module, finite dimensional over $\BC_\infty$. Then,
\begin{align*}
\phi\colon(\sigma_*M)^\ast&\rightarrow \sigma_\ast(M^\ast)\\
\lambda&\mapsto\sigma\circ\lambda
\end{align*}
is well-defined and a $\GL_2(F)$-equivariant isomorphism.
\end{Lem}
\begin{proof}
Let $\lambda\in (\sigma_\ast M)^\ast$, i.e. $\lambda$ is a $\BC_\infty$-linear map $\sigma_\ast M\rightarrow \BC_\infty$. In other words, $\lambda(\sigma(a)m)=a\lambda(m)$ for all $a\in\BC_\infty$ and $m\in M$. Thus, for all $a,b\in\BC_\infty$, $m\in M$, we have
\[
\phi(b\lambda)(am)=\sigma((b\lambda)(am))=\sigma(b)\sigma(\sigma^{-1}(a))\sigma(\lambda(m))=a\sigma(b)\phi(\lambda)(m),
\]
proving the well-definedness.  The $\GL_2(F)$-equivariance is obvious. It is easily checked that composition with $\sigma^{-1}$ gives the inverse map, proving that $\phi$ is an isomorphism.
\end{proof}

Using the notation from \autoref{Section2}, we apply \autoref{Prop-Cartier} in the case $m=l$ and \autoref{Lem-DualFrob} to $M=(\Delta_{pk-2}\otimes\det^{pl-1})\otimes_F\BC_\infty$ to obtain an injective $\GL_2(F)$-equivariant map 
\begin{equation}\label{eq-DefAlphap}
\alpha_p:=\phi\circ\CC_p^*\colon V_{k,l}\otimes_F\BC_\infty\rightarrow\sigma_* (V_{pk,pl}\otimes_F\BC_\infty).
\end{equation}
We will relate this map to the Frobenius map on Drinfeld modular forms in \autoref{Section5}.

\subsection{Irreducible representations of $\SL_2(\BF_q)$ and $\GL_2(\BF_q)$}
\label{Subsec3.3}
We conclude the section with the classification of the irreducible representations of $\BF_q[\SL_2(\BF_q)]$ and $\BF_q[\GL_2(\BF_q)]$ following Bonnaf\'e. This will be relevant in \autoref{Section6}. Define $\overline{L}_k$ as the $\BF_q$-span in $L_k$ of the monomials $X^iY^{k-i}$, $0\le i\le k$, $\binom ki\not\equiv0\pmod p$ and define $\overline{\Delta}_k$ in the same manner. Observe that by definition, $\overline{L}_k=\overline{\Delta}_k$ for $0\leq k\leq p-1$.
\begin{Prop}[{\cite[Thm.~10.1.8(e')]{Bonnafe}}]\label{Prop-IrredsSL2Fq}
Up to isomorphism, the irreducible representations of the ring $\BF_q[\SL_2(\BF_q)]$ are the modules $\overline{L}_k$, $0\le k\le q-1$. They are absolutely irreducible.
\end{Prop}
We write $\overline{\det}^{m}$ for the representation $\GL_2(\BF_q)\times\BF_q\to\BF_q,(\gamma,v)\mapsto (\det\gamma)^m v$.
\begin{Cor}\label{Cor-IrredsGL2Fq}
Up to isomorphism, the irreducible representations of the ring $\BF_q[\GL_2(\BF_q)]$ are the modules $\overline{L}_k\otimes_{\BF_q}\overline{\det}^{m}$, $0\le k\le q-1$, $0\le m\le q-2$.
\end{Cor}
\begin{proof} We shall use the following result that can be found in \cite[\S7, Thm.~1.9]{Weintraub}: Let $G$ be a finite group. Then the number of non-isomorphic irreducible representations of $G$ over $\overline \BF_p$ is equal to the number of $p$-regular conjugacy classes of $G$. 

By the theory of the rational canonical form, any semisimple element is either (i) scalar, or conjugate to a unique companion matrix of the form $\SMat 0 b 1 d$ with $b\neq0$ such that $x^2-dx-b\in\BF_q[x]$ is (ii) irreducible or (iii) a product of two distinct linear factors. The number of classes in (i) is $q-1$, that in (iii) is $1/2(q-1)(q-2)$, that in (ii) is $1/2\deg (x^{q^2}-x)/(x^q-x)=1/2(q^2-q)$. Therefore the number of semisimple (i.e., $p$-regular) conjugacy classes of $\GL_2(\BF_q)$ is
\[ q-1+\frac12(q-1)(q-2)+\frac12(q-1)q=(q-1)q.\]
By \autoref{Prop-IrredsSL2Fq} the representations $L(k,m):=\overline{L}_k\otimes_{\BF_q}\overline{\det}^{m}\otimes_{\BF_q}\overline\BF_p$ are irreducible, and their number is $q(q-1)$. Hence it suffices to show that they are pairwise non-isomorphic. So suppose there is an isomorphism of $\GL_2(\BF_q)$-representations $\phi\colon L(k,m)\rightarrow L(k',m')$. By \autoref{Prop-IrredsSL2Fq} we have $k=k'$. Thus, we can regard $\phi$ as an automorphism of the underlying $\overline\BF_p$-vector space; in particular it has a non-zero eigenvalue $\lambda\in \overline\BF_p$. Let $\mu_\lambda\colon  L(k,m)\rightarrow L(k,m')$ denote the multiplication by $\lambda$. Clearly, $\mu_\lambda$ is $\SL_2(\BF_q)$-equivariant. Thus, $\phi-\mu_\lambda$ is a morphism of $\SL_2(\BF_q)$-representations with non-zero kernel, hence by \autoref{Prop-IrredsSL2Fq}, $\phi=\mu_\lambda$. Now, we consider the elements $g=\SMat \xi 00 1$, $\xi\in\BF^\times_q$. We obtain 
\[
\lambda\xi^{m-k}X^k=\phi(g\cdot X^k)=g\cdot\phi(X^k)=\lambda\xi^{m'-k}X^k.
\]
Since $\lambda\neq 0$, this implies $m \equiv m'\pmod {q-1}$ and thus $m=m'$.
\end{proof}

\section{Hecke-equivariant maps between spaces of Drinfeld cusp forms}
\label{Section5}
In this section we apply the machinery developed in \autoref{Section2} to the maps constructed in \autoref{Subsec3.2} and \autoref{Cartier}. As a result, we obtain Hecke equivariant maps betweens spaces of (adelic) Drinfeld cusp forms of different weights. More generally, in \autoref{Prop-HeckeFil} we see that any filtration of the $\GL_2(F)$-module $V_{k,l}$ gives rise to a Hecke-stable filtration on Drinfeld cusp forms of weight $k$ and type $l$. Afterwards, we show that our constructions coincide with previous work of Bosser--Pellarin and the well-known Frobenius map in the special case $A=\BF_q[t]$.  We conclude the section with some computational examples.
\label{Sec-HyperDer}
\subsection{Maps induced from representation theory}
\label{SubSec5.1}
Recall that in \autoref{Prop-DsAsMapOnReps} we constructed certain hyperderivative maps on symmetric powers: Let $k\geq 2$, $s\geq 1$, $m\in \BZ$ and suppose  that $\binom{k+s-1}i=0\pmod p$ for $i=1,\ldots,s$. Then there is an explicit $\GL_2(F)$-equivariant map
\[ \mathcal{D}_s\colon \Delta_{k-2+2s}\otimes \det^{m+s}\to \Delta_{k-2}\otimes \det^{m}. \]
Upon noting that $V_{k,l}=(\Delta_{k-2}\otimes \det^{l-1})^\ast$, we obtain (by duality) a $\GL_2(F)$-equivariant map $\mathcal{D}^\ast_s\colon V_{k,l}\rightarrow V_{k+2s,l+s}$. Now, let $\CK$ be any compact open subgroup of $\GL_2(\BA_F^\infty)$. We can apply \autoref{cor:AdelicStProps}(d) to obtain a Hecke-equivariant map
\[
\mathcal{D}^\ast_s\colon \St_{\CK}\otimes_{\GL_2(F)}V_{k,l}\rightarrow \St_{\CK}\otimes_{\GL_2(F)}V_{k+2s,l+s}.
\]
Finally, by invoking the isomorphisms in \autoref{thm:CHar-St-global} and \autoref{thm:Sk-CHar-Adel}, we obtain the Hecke-equivariant map
\[
\mathcal{D}^\ast_s\colon S_{k,l}(\CK)\to S_{k+2s,l+s}(\CK).
\]
In particular, $\mathcal{D}^\ast_s(S_{k,l}(\CK))\subset S_{k+2s,l+s}(\CK)$ is a Hecke-stable subspace. 

We obtain the following consequence of \autoref{Prop-DsAsMapOnReps}.

\begin{Cor} \label{HCCFsubspacecor}
Suppose  that $\binom{k+s-1}i=0\pmod p$ for $i=1,\ldots,s$ and let $W$ be as in \autoref{Prop-DsAsMapOnReps}. Then, we have
\[
\mathcal{D}^\ast_s(\St_{\CK}\otimes_{\GL_2(F)}V_{k,l})\subset \St_{\CK}\otimes_{\GL_2(F)}W^\perp.
\]
\end{Cor}

\begin{Rem}
We directly obtain that for all $s \geq 1$ such that $\mathcal{D}_s$ is defined and such that $k+2s-2$ satisfies the condition on its base $p$-expansion in \autoref{Lem-BasicsOnLk}(e), we must have $\mathcal{D}^\ast_s(\St_{\CK}\otimes_{\GL_2(F)}V_{k,l})= \{0\}$. Indeed, for such weights $\Delta_{k+2s - 2} = L_{k+2s-2}$ is irreducible as a $\GL_2(F)$-module. So, suppose that there is an $s \geq 1$ such that $\mathcal{D}_s$ is defined. Since the target space will always have smaller dimension, $\mathcal{D}_s$ must have a kernel, and as  $\Delta_{k+2s - 2} = L_{k+2s-2}$ is irreducible, we deduce that $\mathcal{D}_s$ is the zero map. Thus, $W^\perp = \{0\}$, and hence, by \autoref{HCCFsubspacecor}, $\mathcal{D}^\ast_s(\St_{\CK}\otimes_{\GL_2(F)}V_{k,l})= \{0\}$ as desired.
\end{Rem}

At this point it is natural to ask how many Hecke-stable subspaces arise via hyperderivatives, i.e. how many hyperderivatives are defined for a fixed target space.

\begin{Prop} \label{HeckestableProp}
Let $k\geq 2$, and write $k-1 = \sum_{i = 0}^r k_i p^i$ in base $p$. For each $0 \leq j \leq r-1$, let $s_j = \sum_{i = 0}^j k_i p^i$. For each such $j$, the subspace
\[ \mathcal{D}^\ast_{s_j}S_{k-2s_j,l-s_j}(\CK) \subset S_{k,l}(\CK) \]
is a Hecke-stable subspace of $S_{k,l}(\CK)$. Furthermore, these are all of the Hecke-stable subspaces of $S_{k,l}(\CK)$ arising from the hyperderivatives defined in \autoref{Prop-DsAsMapOnReps}.
\end{Prop}
\begin{proof}
The Hecke-stability is immediate from the above. We only need to check that $\mathcal{D}_{s_j}$ is defined and that these are all possible choices of $s$.
By \autoref{Lem-OnBinomials}, we have ${k-1-s \choose i} \equiv 0 \pmod{p}$ for all $1 \leq i \leq s$ if and only if $s$ is of the form $s = s_j := \sum_{i = 0}^j k_i p^i$, for some $0 \leq j \leq r-1$, in particular, this shows that there are at most $\lfloor \log_p(k-1) \rfloor $ possible such $s$. 
\end{proof}

\begin{Rem}
It follows immediately from \autoref{HeckestableProp} that for weights of the form $k = p^n+1-c p^{n-1}$, with $0 \leq c < p-1$ and $n \geq 1$, there are no nonzero $s$ such that $\mathcal{D}^\ast_{s}S_{k-2s,l-s}(\CK) \subset S_{k,l}(\CK)$. Thus, there are no obstructions to the irreducibility of $S_{k,l}(\CK)$ as a Hecke module coming from the hyperderivatives considered in this paper. We will investigate a Maeda-type conjecture for Drinfeld cusp forms of weights of this form in \autoref{MaedaSection}. 

We point out that for $k = p^n + 1-cp^{n-1}$ the corresponding $V_{k,l}=(\Delta_{k-2}\otimes \det^{l+1-k})^\ast$ is irreducible as $k-2 = p^n - 1-cp^{n-1}$ is a magic number and thus $\Delta_{k-2} = L_{k-2}$; see \autoref{Lem-BasicsOnLk}(e). So we can even generalize the previous remark to the statement that there are no obstructions to the ireducibility of the Hecke module $S_{k,l}(\CK)$ coming from representation theory.
\end{Rem}
We can also apply the above technique to the Cartier operator constructed in \autoref{Cartier}: Recall that  we constructed a $\GL_2(F)$-equivariant injective  map 
\[
\alpha_p\colon V_{k,l}\otimes_F\BC_\infty\rightarrow\sigma_* (V_{pk,pl}\otimes_F\BC_\infty).
\] 
Again, we apply \autoref{cor:AdelicStProps} to obtain the injective Hecke-equivariant map
\[
\alpha_p\colon \St_{\CK}\otimes_{\GL_2(F)}(V_{k,l}\otimes_F\BC_\infty)\rightarrow \St_{\CK}\otimes_{\GL_2(F)}\sigma_*(V_{pk,pl}\otimes_F\BC_\infty).
\]
Now, since $\St_\CK$ is defined over $\BZ$, we obtain a canonical isomorphism
\[
\St_{\CK}\otimes_{\GL_2(F)}\sigma_*(V_{pk,pl}\otimes_F\BC_\infty)\simeq
\sigma_\ast(\St_{\CK}\otimes_{\GL_2(F)}(V_{pk,pl}\otimes_F\BC_\infty)).
\]
Thus, by invoking the isomorphisms in \autoref{thm:CHar-St-global} and \autoref{thm:Sk-CHar-Adel} we obtain a Hecke-equivariant map
\[
\alpha_p\colon S_{k,l}(\CK)\rightarrow \sigma_*S_{pk,pl}(\CK).
\]

While the above maps are so far the only explicit source of Hecke-stable subspaces, we are able to state the following more general result.
\begin{Prop}\label{Prop-HeckeFil}~
\begin{compactenum}
\item Any filtration of the $\GL_2(F)$-module $V_{k,l}$ induces Hecke-stable filtrations of $C^\ad_\har(V_{k,l},\CK)^{\GL_2(F)}$ and $S_{k,l}(\CK)$.
\item If a subrepresentation $M$ of $V_{k,l}$ is isomorphic to a subrepresentation $M'$ of $V_{k',l'}$,then any Hecke-eigensystem in $S_{k,l}(\CK)$ corresponding to $M$ also appears in $S_{k',l'}(\CK)$
\end{compactenum}
\end{Prop}

\begin{proof}
Both parts are immediate from \autoref{thm:CHar-St-global}(b).
\end{proof}

\subsection{Hyperderivatives and the work of Bosser--Pellarin}
\label{RelBosPel}
From here on we assume $A=\BF_q[t]$. Denote by $\Omega$ the Drinfeld symmetric space $\BP^1(\BC_\infty)\setminus\BP^1(F_\infty)$ with its natural structure as a rigid space, e.g.~\cite[\S~5]{Gekeler}. Let $f\colon\Omega\to \BC_\infty$ be a locally analytic function and $z\in\Omega$. Then the hyperderivatives $D_s f$ at $z$ are defined by the formula
\[ f(z+\eps)=\sum_{s\ge 0} (D_sf)(z)\eps^s \]
for $\eps\in\BC_\infty$ with $|\eps|$ sufficiently small; see \cite[Def.~2.3]{Uchino-Satoh} or \cite[\S~3.1]{Bosser-Pellarin-HyperDiff}.

Here are some formal properties of hyperderivatives:
\begin{compactenum}
\item The functions $D_sf$ are locally analytic on $\Omega$, and rigid analytic, if $f$ is so.
\item If $f$ has a Laurent series expansion $f=\sum_{n\in\BZ} a_n(z-b)^n$ converging on an annulus $r<|z-b|<R$, then $D_sf$ has the Laurent series expansion
\[D_sf= \sum_{n\in\BZ} a_n\binom{n}{s}(z-b)^{n-s} = \sum_{n\in\BZ^{(s)}} a_{n+s}\binom{n+s}{s}(z-b)^{n} ,\]
where $\BZ^{(s)}=\BZ\setminus\{-s,-s+1,\ldots,-1\}$; use $(z+\eps)^i=z^i\big(1+\frac \eps z)^i=z^i\sum_{j\ge0 }\binom{i}j\big(\frac\eps z\big)^j$.
\item $D_s(fg)=\sum_{r=0}^s D_r f\cdot D_{s-r}g$.
\item $D_i D_j=D_j D_i = \binom{i+j}i D_{i+j}$.
\end{compactenum}

\begin{Lem}\label{Lem-ResidueAndHyperDeriv}
Consider the annulus $A:=\{z\in\BC_{\infty}\mid r<|z-c|<R\}$ for $c\in\BC_\infty$ and $0<r<R$ rational numbers. For rigid analytic functions $f,g\colon A\to \BC_\infty$ one has
\[\Res_A D_s f \cdot g\mathrm{d}z=(-1)^s \Res_A f \cdot D_sg\mathrm{d}z.\]
\end{Lem}
\begin{proof}
By explicit computation: Write $f=\sum_{n\in\BZ} a_n(z-c)^n$ and $g=\sum_{n\in\BZ} b_n(z-c)^n$. Then by (b) above we have
\begin{eqnarray*}
\Res_A D_s f \cdot g\mathrm{d}z&=&\sum_{n\in\BZ^{(s)}} a_{n+s}\binom{n+s}{s}b_{-1-n}\\
&=&\sum_{n\ge 0} a_{n+s}\binom{n+s}{s}b_{-1-n} + \sum_{n\ge0} a_{-1-n}\binom{-1-n}{s}b_{n+s}\\
&=&\sum_{n\ge 0} a_{n+s}\binom{n+s}{s}b_{-1-n} + (-1)^s\sum_{n\ge0} a_{-1-n}\binom{n+s}{s}b_{n+s}\\
&=&\sum_{n\ge 0} \binom{n+s}{s}\big(a_{n+s}b_{-1-n} + (-1)^s a_{-1-n}b_{n+s}\big).
\end{eqnarray*}
Clearly if we exchange $a_?$ and $b_?$, then the sign changes by $(-1)^s$, and this proves the assertion.
\end{proof}

From \cite[Thm.~4.1 and formula (45)]{Bosser-Pellarin-DQMF}, we quote the following result:
\begin{Prop}[Bosser-Pellarin]
Let $A=\BF_q[t]$, $k\ge1$ and $l\in\BZ/(q-1)$. Suppose that $s\ge1$ satisfies $\binom{k+s-1}i=0\pmod p$ for $i=1,\ldots,s$. Then, up to a character twist, $D_s$ defines is a linear Hecke-equivariant map of spaces of modular forms 
\[D_s\colon M_{k,l}(\GL_2(A))\to M_{k+2s,l+s}(\GL_2(A)).\]
The map $D_s$ preserves the subspaces of cusp forms.
\end{Prop}

Note that for $A=\BF_q[t]$ and $\CK=\GL_2(\hat{A})$, we have that $\mathrm{Cl}_\CK$ is the trivial group and $C^\ad_\har(N,\CK)^{\GL_2(F)}\cong C_\har(N)^\Gamma$ for $\Gamma=\GL_2(A)$. Thus, upon observing that this isomorphism is functorial in $N$, we can reduce the results of the previous subsection to the local situation.

Denote the residue map $S_{k,l}(\Gamma)\to C_\har(V_{k,l})^\Gamma$ by $f\mapsto c_{k,l}(f)$, and recall that $c_{k,l}$ is defined by
\begin{equation}
(c_{k,l}(f))(e)(X^iY^j)=\Res_e f z^i \mathrm{d}z :=\Res_{A_e} f z^i \mathrm{d}z,
\end{equation}
where $e$ is any oriented edge of the Bruhat-Tits tree, $A_e$ denotes the associated annulus, see \cite[Preliminaries]{Teitelbaum-Poisson}, and $i+j=k-2$

\begin{Cor}\label{Cor-BP}
Suppose that $\binom{k+s-1}i=0\pmod p$ for $i=1,\ldots,s$. Then
\[c_{k+2s,l+s}(D_s f)=\mathcal{D}^\ast_s c_{k,l}(f)\]
\end{Cor}

\begin{proof} Let $i,j\geq 0$ such that $i+j=k-2+2s$. We have
\begin{eqnarray*}
(c_{k+2s,\ell+s}(D_s f))(e)(X^iY^j)&=&\Res_e (D_sf z^i \mathrm{d}z)\\
&\stackrel{\ref{Lem-ResidueAndHyperDeriv}}=&(-1)^s\Res_e (f D_s z^i \mathrm{d}z)\\
&=&(-1)^s\binom{i}s \Res_e (f z^{i-s} \mathrm{d}z).
\end{eqnarray*}
The assertion follows now formally from the definition of $c_{k,\ell}(f)(e)(X^{i-s}Y^{j-s})$ and the definition of $\mathcal{D}_s$ in \autoref{Prop-DsAsMapOnReps}.
\end{proof}

Note that while the map $\mathcal{D}^\ast_s$ is Hecke-equivariant by construction (in the adelic situation), it is shown in  \cite[Section 4.1.1]{Bosser-Pellarin-DQMF}, that $D_s$ is Hecke-equivariant only up to twist by an explicit character. The reason for this phenomenon is that the Hecke-actions in the adelic and local situation use different normalizations. Let $\mathfrak{p}A$ denote a maximal ideal in $A$ with monic generator $\mathfrak{p}$ and denote by $T_\mathfrak{p}'$ and $T_\mathfrak{p}$ the associated Hecke-operators on $S_{k,l}(\CK)$ and $S_{k,l}(\Gamma)$ as defined in \cite[Ch.~6]{Boeckle-EiSh} and \cite[\S~6]{Gekeler}. Then one has $T_\mathfrak{p}'=\mathfrak{p}^{l-k}T_\mathfrak{p}$ by \cite[Example 6.13]{Boeckle-EiSh}. In particular, in the situation in \autoref{Cor-BP}, one has
\begin{align*}
c_{k+2s,l+s}(D_s T_\mathfrak{p}f)&=\mathcal{D}^\ast_sc_{k,l}(T_\mathfrak{p}f)=\mathfrak{p}^{k-l}\mathcal{D}^\ast_sc_{k,l}(T_\mathfrak{p}'f)\\
&=\mathfrak{p}^{k-l}T_\mathfrak{p}'\mathcal{D}^\ast_sc_{k,l}(f)=\mathfrak{p}^{k-l}T_\mathfrak{p}'\mathcal{D}^\ast_sc_{k,l}(f)\\
&=\mathfrak{p}^{k-l}T_\mathfrak{p}'c_{k+2s,l+s}(D_sf)\\
&= \mathfrak{p}^{k-l}\mathfrak{p}^{l+s-k-2s}T_\mathfrak{p}c_{k+2s,l+s}(D_sf).
\end{align*}
In other words, the diagram

\begin{displaymath}
    \xymatrix{ S_{k,l}(\Gamma) \ar[d]_{\mathfrak{p}^sT_\mathfrak{p}}\ar[r]^-{D_s} & S_{k+2s,l+s}(\Gamma) \ar[d]^{T_\mathfrak{p}} \\
               S_{k,l}(\Gamma) \ar[r]_-{D_s}& S_{k+2s,l+s}(\Gamma)  }
\end{displaymath}
commutes and if $f$ is an eigenform for $T_\mathfrak{p}$ with eigenvalues $a_\mathfrak{p}$, then $D_s f$ is an eigenform for $T_\mathfrak{p}$ with eigenvalue $a_\mathfrak{p}\mathfrak{p}^s$, which recovers \cite[Lemma 4.6]{Bosser-Pellarin-DQMF}. By the same reasoning, we can reformulate \autoref{Prop-HeckeFil}(b) in this special situation as follows.
\begin{Cor}
If a subrepresentation $M$ of $V_{k,l}$ is isomorphic to a subrepresentation $M'$ of $V_{k',l'}$, then  for any Hecke-eigensystem $(a_\mathfrak{p})_\mathfrak{p}$ in $S_{k,l}(\Gamma)$ corresponding to $M$ there is a Hecke eigensystem $(a_\mathfrak{p}\mathfrak{p}^{k'-k-(l'-l)})_\mathfrak{p}$ in $S_{k',l'}(\Gamma)$.
\end{Cor}

\subsection{The Frobenius on Drinfeld cusp forms}
We keep the notation from the previous subsection, except we allow general congruence subgroups $\Gamma$. The injective map $\tau_p\colon S_{k,l}(\Gamma)\rightarrow S_{pk,pl}(\Gamma), f\mapsto f^p$ is Frobenius-linear and behaves nicely with respect to the Hecke-operators:  If $f\in S_{k,l}(\Gamma)$ is a Hecke eigenform with Hecke eigensystem $(a_\mathfrak{p})_\mathfrak{p}$, then the form $\tau_p(f)$ is again a Hecke eigenform with eigensystem $(a_\mathfrak{p}^p)_\mathfrak{p}$. If we denote the Frobenius morphism on the perfect field $\BC_\infty$ by $\sigma$, we can reformulate the above to the following statement:
\begin{Prop}
\label{Prop-frob}
The map $\tau_p\colon S_{k,l}(\Gamma)\rightarrow \sigma_\ast S_{pk,pl}(\Gamma)$ is Hecke-equivariant. 
\end{Prop}

It what follows, we give a representation theoretic explanation for this phenomenon using the map $\alpha_p$ constructed in \eqref{eq-DefAlphap} in \autoref{Cartier}.
\begin{Prop}\label{Prop-FrDrin}
We have
\[
c_{pk,pl}(\tau_p(f))(e)=\alpha_p(c_{k,l}(f)(e)).
\]
\end{Prop}
\begin{proof}
Assume that the annulus corresponding to $e$ is given by $\{z\in\BC_\infty\mid r<|z-c|<R\}$ for $c\in\BC_\infty$ and $0<r<R$. Write $f=\sum_{n\in\BZ}a_n(z-c)^n$. Then $\tau_p(f)=\sum_{n\in\BZ}a_n^p(z-c)^{pn}$. Thus, if we write
\[
z^{i-1}f=(z-c+c)^{i-1} \sum_{n\in\BZ}a_n(z-c)^n=\left(\sum_{m=1}^{i}\binom{i-1}{m-1}c^{i-m}(z-c)^{m-1}\right)\left(\sum_{n\in\BZ}a_n(z-c)^n\right),
\]
we obtain $\Res_e(fz^{i-1}\mathrm{d}z)=\sum_{m=1}^{i}\binom{i-1}{m-1}c^{i-m}a_{-m}$. Consequently, we have
\begin{eqnarray*}
\alpha_p(c_{k,l}(f)(e))(X^{i-1}Y^{j-1})&=&\sigma(c_{k,l}(f)(e)(\CC_p(X^{i-1}Y^{j-1}))\\
&=&
\left\{
\begin{aligned}
&\left(\sum_{m=1}^{i/p}\binom{i/p-1}{m-1}c^{i/p-m}a_{-m}\right)^p
&& \text{if } p\mid i,\\
& 0, && \text{otherwise.}\\
\end{aligned}
\right.
\end{eqnarray*}
On the other hand, we have
\[
\begin{aligned}
c_{pk,pl}(\tau_p(f))(e)(X^{i-1}Y^{j-1})&=\Res_e(f^pz^{i-1}\mathrm{d} z)\\
&=\sum_{\substack{m=1\\ p\mid m}}^{i}\binom{i-1}{m-1}c^{i-m}a^p_{-m/p}=\sum_{n=1}^{\lfloor i/p \rfloor}\binom{i-1}{pn-1}c^{i-pn}a^p_{-n}.
\end{aligned}
\]
Since $pn-1=p(n-1)+p-1$, we see that $\binom{i-1}{pn-1}=0$ if $p\nmid i$ by  \autoref{Prop-Lucas}, which shows that in this case
\[
c_{pk,pl}(\tau_p(f))(e)(X^{i-1}Y^{j-1})=0.
\]
Assume now that $p\mid i$. Then we have, again by \autoref{Prop-Lucas}, 
\[
\binom{i-1}{pn-1}=\binom{p(i/p-1) +p-1}{p(n-1)+p-1}=\binom{i/p-1}{n-1},
\]
and thus,
\[
c_{pk,pl}(\tau_p(f))(e)(X^{i-1}Y^{j-1})=\sum_{n=1}^{i/p}\binom{i/p-1}{n-1}c^{p(i/p-n)}a^p_{-n}=\left(\sum_{n=1}^{i/p}\binom{i/p-1}{n-1}c^{i/p-n}a_{-n}\right)^p,
\]
which completes the proof; note that binomials mod $p$ lie in $\BF_p$ and are thus fixed by $\sigma$.
\end{proof}

Thus, the Hecke-equivariant map constructed in \autoref{SubSec5.1},
\[
\alpha_p\colon S_{k,l}(\CK)\rightarrow \sigma_*S_{pk,pl}(\CK),
\]
is just $\tau_p$ when switching to the local situation. As the final step to obtain a representation theoretic proof for \autoref{Prop-frob}, note that contrary to the situation in \autoref{RelBosPel}, the different normalizations of the local and adelic Hecke operators are not visible in this case as
\begin{align*}
c_{pk,pl}(\tau_p (T_\mathfrak{p}f))&=\alpha_p(c_{k,l}(T_\mathfrak{p}f))=\alpha_p(\mathfrak{p}^{k-l}c_{k,l}(T_\mathfrak{p}'f))\\
&=\mathfrak{p}^{pk-pl}\alpha_p(T_\mathfrak{p}'c_{k,l}(f))=\mathfrak{p}^{pk-pl}T_\mathfrak{p}'\alpha_p(c_{k,l}(f))\\
&=T_\mathfrak{p}\alpha_p(c_{k,l}(f))=c_{pk,pl}(T_\mathfrak{p}\tau_p(f)).
\end{align*}

\subsection{Examples}
In the following examples, we always fix $A=\BF_q[t]$ and $\Gamma=\GL_2(A)$.
\subsubsection*{Non-vanishing hyperderivatives via Petrov's special family}
In his thesis, A. Petrov introduced the following family ($n \geq 1$) of single cuspidal Hecke eigenforms $f_n \in S_{2+n(q-1),1}(\Gamma)$ with $A$-expansions,
\[ f_n = \sum_{\substack{a \in A \\ \text{monic}}} a^{1+n(q-1)} u_a,\]
where for each monic $a \in A$ one defines $u_a(z) := u(az)$. The importance of Petrov's family is that it gives an explicit description of all cuspforms in level 1 which vanish exactly to the order 1 at the cusp at infinity. It is not difficult to show that the analytic hyperderivatives defined in this section transform Petrov's forms in the following way
\[ D_sf_n = (-\widetilde{\pi})^s \sum_{\substack{a \in A \\ \text{monic}}} a^{s+1+n(q-1)} G_{s+1}(u_a),\]
where $G_{s+1}$ is the $s+1$-th Goss polynomial for the lattice $\widetilde{\pi} A$ as in \cite{Gekeler}; see e.g. \cite[(2)]{Petrov-Hyperderivatives} for this calculation. It follows then from \cite[Theorem 2.2]{Petrov-Aexpansion} that all of the forms $D_s f_n$ for $n,s \geq 1$ are non-zero functions, and thus whenever $D_s$ preserves modularity, the forms $D_s f_n$ give examples of non-zero cuspidal Hecke eigenforms in the image of the analytic hyperdifferential operators.

\subsubsection*{Chains of hyperderivative maps and interactions with the Frobenius}
We have observed computationally that it is possible to have chains of hyperderivative maps between spaces of various weights, but sometimes no direct map. So, for example for $q = 3$ we have chains of hyperderivatives
\[ S_{62,0}(\Gamma) \xrightarrow{D_{2}}{} S_{66,0}(\Gamma) \xrightarrow{D_{16}}{} S_{98,0}(\Gamma) \xrightarrow{D_2}{} S_{102,0}(\Gamma) 
\quad \text{and a direct map}\quad  S_{62,0}(\Gamma) \xrightarrow{D_{20}}{} S_{102,0}(\Gamma)\]
but \emph{no} direct maps
\[S_{62,0}(\Gamma) \xrightarrow{D_{18}}{} S_{98,0}(\Gamma) \quad \text{and} \quad  S_{66,0}(\Gamma) \xrightarrow{D_{18}}{} S_{102,0}(\Gamma). \]
Notice that one has $D_{2}D_{16} = D_{16}D_2 = {18 \choose 2}D_{18} \equiv 0$ in characteristic 3, and so while $D_{18}$ does not preserve modularity (i.e. $D_{18}$ does not meet the requirements of Bosser-Pellarin's result) on either of the spaces $S_{62,0}(\Gamma)$ or $S_{66,0}(\Gamma)$ when $q = 3$, the compositions $D_2 D_{16}$ and $D_{16}D_2$ still do, being just the zero map in all cases. As all of these maps are between spaces of type 0, Petrov's examples from the previous section do not apply.  There are also Frobenius maps
 \[ S_{22,0}(\Gamma) \xrightarrow{\tau_3}{} S_{66,0}(\Gamma) \quad \text{and}\quad S_{34,0}(\Gamma) \xrightarrow{\tau_3}{} S_{102,0}(\Gamma). \]
We have computed the action of the Hecke operator $T_t$ on the image of all of these maps using \texttt{Sage} which we summarize in the following diagram.

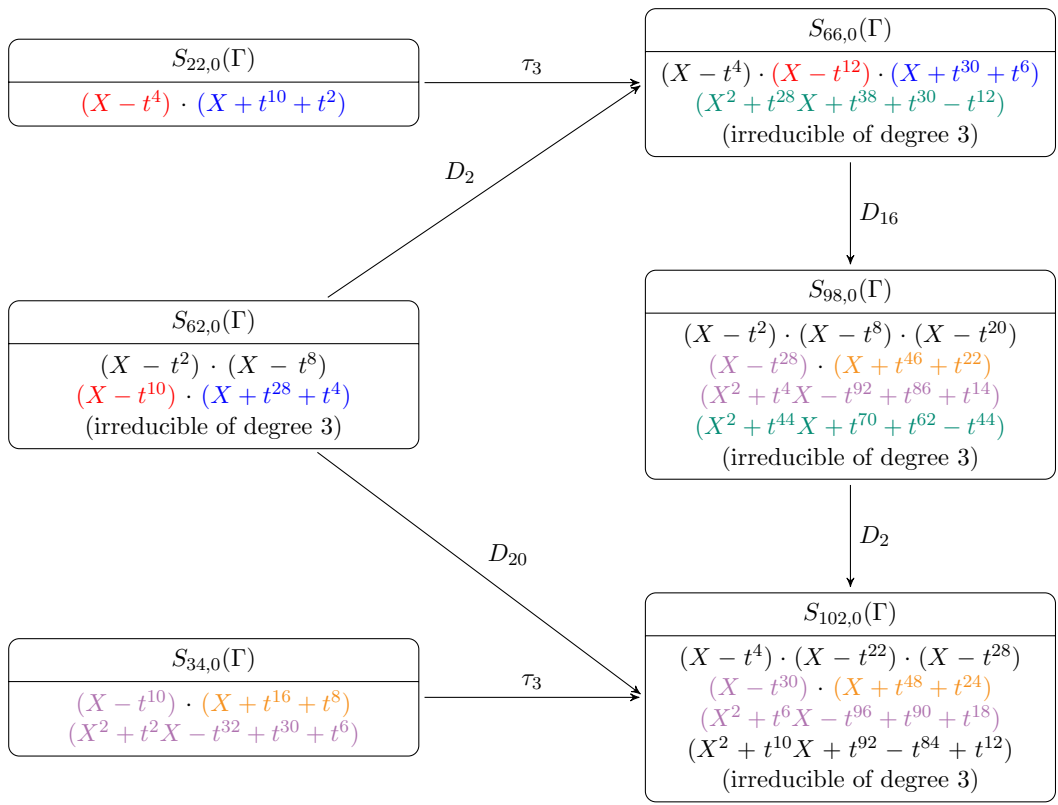
\begin{center}
\begin{tikzpicture}[node distance=3cm, auto,every node/.style={scale=0.85}]
 \node[punkt, shape=rectangle split, rectangle split parts=2](M102){$S_{102,0}(\Gamma)$\nodepart{two}
$(X -t^4)\cdot(X -t^{22})\cdot(X -t^{28})$\\
$\textcolor{Orchid}{(X -t^{30})}\cdot \textcolor{BurntOrange}{(X + t^{48} + t^{24})}$\\
$\textcolor{Orchid}{(X^2 + t^6X -t^{96} + t^{90} + t^{18})}$\\
$(X^2 + t^{10}X + t^{92} -t^{84} + t^{12})$ \\ $(\text{irreducible of degree 3})$};

 \node[punkt, shape=rectangle split, rectangle split parts=2, left=3cm of M102](M34){$S_{34,0}(\Gamma)$\nodepart{two}
$\textcolor{Orchid}{(X -t^{10})}\cdot\textcolor{BurntOrange}{(X + t^{16} + t^8)}$\\
$\textcolor{Orchid}{(X^2 + t^2X -t^{32} + t^{30} + t^6)}$}
    edge[pil] node[auto] {$\tau_3$} (M102.west);

 \node[punkt, shape=rectangle split, rectangle split parts=2, above=1.5cm of M102](M98){$S_{98,0}(\Gamma)$\nodepart{two}
$(X -t^2)\cdot(X -t^8)\cdot(X -t^{20})$\\
$\textcolor{Orchid}{(X -t^{28})}\cdot\textcolor{BurntOrange}{(X + t^{46} + t^{22})}$ \\ $\textcolor{Orchid}{(X^2 + t^4X -t^{92} + t^{86} + t^{14})}$\\
$\textcolor{PineGreen}{(X^2 + t^{44}X + t^{70} + t^{62} -t^{44})}$ \\ $(\text{irreducible of degree 3})$}
edge[pil] node[auto] {$D_{2}$} (M102.north);

 \node[punkt, shape=rectangle split, rectangle split parts=2, above=1.5cm of M98](M66){$S_{66,0}(\Gamma)$\nodepart{two}
$(X -t^4)\cdot\textcolor{red}{(X -t^{12})}\cdot\textcolor{blue}{(X + t^{30} + t^6)}$ \\ $\textcolor{PineGreen}{(X^2 + t^{28}X + t^{38} + t^{30} -t^{12})}$ \\ $(\text{irreducible of degree 3})$}
  edge[pil] node[auto] {$D_{16}$} (M98.north);

 \node[punkt, shape=rectangle split, rectangle split parts=2, left=3cm of M98](M62){$S_{62,0}(\Gamma)$\nodepart{two}
$(X -t^2)\cdot(X -t^8)$\\
$\textcolor{red}{(X -t^{10})}\cdot\textcolor{blue}{(X + t^{28} + t^4)}$ \\ $(\text{irreducible of degree 3})$}
   edge[pil] node[auto] {$D_2$} (M66.west)
   edge[pil] node[auto] {$D_{20}$} (M102.west);

 \node[punkt, shape=rectangle split, rectangle split parts=2, left=3cm of M66](M22){$S_{22,0}(\Gamma)$ \nodepart{two}  $\textcolor{red}{(X-t^4)}\cdot\textcolor{blue}{(X+t^{10}+t^2)}$}
    edge[pil] node[auto] {$\tau_3$} (M66.west);

\end{tikzpicture}
\end{center}
The lower part in each box displays the factorized characteristic polynomial of $T_t$ acting on the space of cusp forms displayed in the upper part of the box. The colors indicate which factors come from factors in lower weight via the maps above. We should point out that one such link is not displayed in the above diagram: The factor $\textcolor{BurntOrange}{(X+t^{48}+t^{24})}$ in weight $102$ in fact comes via $D_{20}$ from the factor $\textcolor{blue}{(X+t^{28}+t^4)}$ in weight $62$, i.e. should be colored blue. However, this coloring would be misleading, because this factor can't possibly explain the factor $\textcolor{BurntOrange}{(X+t^{46}+t^{22})}$ in weight $98$ as $D_{16} D_2=0$, nevertheless it behaves exactly as one would expect assuming the composition would be non-zero, which we find to be quite an interesting observation. We plan to investigate the interactions between the Frobenius and the hyperderivatives further in future work.

\section{Dimension formulas for $\SL_2(A)$-forms, $A=\BF_q[t]$}
\label{Section6}

Using representation theory for $\GL_2(F)$ we were able to gain in \autoref{Cor-V_k-versus-L_k} some understanding of the Jordan-H\"older constituents of $V_{k,l}$. By \autoref{Prop-HeckeFil}, any composition series of $V_{k,l}$ as an $\GL_2(F)$-module will induce a Hecke-stable filtration of $C_\har(V_{k,l})^\Gamma$ with subquotients isomorphic to $C_\har(L_{k'}\otimes\det^{m'})^\Gamma$ for suitable $(k',m')$. Thus, the number of non-zero such subquotients yields a lower bound on the number of factors of the characteristic polynomial for each Hecke operator $T_\mathfrak{p}$ acting on $S_{k,l}(\SL_2(A))$, and the dimension of each non-zero such subquotient yields an upper-bound on the degree of the factor of the characteristic polynomial for $T_\mathfrak{p}$ corresponding to the ``restriction of $T_\mathfrak{p}$ to this subquotient." We will see in \autoref{MaedaSection} that even in cases where there is just one non-zero irreducible subquotient, the characteristic polynomials of the Hecke operators can still factor further over $F$, and thus more work will be needed to understand this phenomenon.

Thus, in any effort to formulate an analog of the Maeda conjecture in the setting of Drinfeld modular forms, one is forced to grapple with the various modules $C_\har(L_{k'}\otimes\det^{m'})^\Gamma$, as above. 
The present section explains, how for $A=\BF_q[t]$ and $\Gamma=\SL_2(A)$, one can in principle compute the dimensions of the building blocks $C_\har(L_k)^{\SL_2(A)}$. The outcome, which we shall make explicit for $q\in\{2,3,5\}$, looks simple; see Propositions~\ref{prop:DimLkQ3}, \ref{prop:DimLkQ2} and \ref{prop:DimLkQ5}, but we have not found a simple way to prove this or to give closed formulas in general. Our explicit calculations do suggest that we have a rough estimate for the dimension given by
\[ \dim_F C_\har(L_k)^{\SL_2(A)} \approx \frac{\gcd(2,q^2-1)}{q^2-1}\dim_F L_k.\]
In \autoref{Prop-SteinbAndChar} we express $\dim C_\har(V)^{\SL_2(A)}$ for any $\SL_2(F)$-representation $V$ that is finite dimensional over $F$ in terms of a formula that only involves the action restricted to $\SL_2(\BF_q)$ and the Steinberg module $\overline\st$ for the latter finite group. For this we work out in detail some material from \cite{Teitelbaum-Poisson}. In \autoref{subsec:Reduzzi} we explain how from this and work of Reduzzi on $K_0$ of the representation category of $\BF_q[\SL_2(\BF_q)]$, in principle, one can derive explicit formulas for $\dim C_\har(V)^{\SL_2(A)}$. To obtain similar results for $\dim C_\har(V)^{\GL_2(A)}$, one would need to  extend \cite{Reduzzi} to $\GL_2(\BF_q)$.

\subsection{A model for $\SL_2(A)$-forms, $A=\BF_q[t]$}
\label{subsec:Model}

Throughout this section, let $A=\BF_q[t]$ and let $V$ denote a representation of $\SL_2(F)$ on a finite-dimensional $F$-vector space. Denote by $\ost:=\overline{L}_{q-1}$ the Steinberg module for $\BF_q[\GL_2(\BF_q)]$. It is a projective $\BF_q[\GL_2(\BF_q)]$-module by \cite[Lemma 10.2.4]{Bonnafe}. Let 
\begin{equation} \label{SB2def} \SB_2(\BF_q)=\{\SMat{a}bcd\in\SL_2(\BF_q)\mid c=0\}.
\end{equation}
Let $\CX:=\SL_2(\BF_q)/\SB_2(\BF_q)$. One can show that $\ost$ is the kernel of the $\BF_q[\SL_2(\BF_q)]$-linear map
\begin{equation}\label{eqn-Stbg-Def}
\BF_q[\CX]\longto \BF_q, \sum_{x\in \CX} a_x x\mapsto \sum_{x\in \CX} a_x,
\end{equation}
and observe that the map is split by $\BF_q\to\BF_q[\CX],a\mapsto a\sum_{x\in\CX}x$, because $\#\CX=q+1\equiv 1\pmod p$. Note that one has the $\BF_q[\SL_2(\BF_q)]$-linear isomorphism 
\[\BF_q[\CX]\stackrel\simeq\longto \Ind_{\SB_2(\BF_q)}^{\SL_2(\BF_q)}\BF_q=\BF_q[\SL_2(\BF_q)]\otimes_{\BF_q[\SB_2(\BF_q)]}\BF_q,g\SB_2(\BF_q)\mapsto g\otimes 1.\]

The following result is an analog of the isomorphism for type $1$ Drinfeld modular forms for $\GL_2(A)$ stated at the bottom of \cite[p.~507]{Teitelbaum-Poisson}. 
\begin{Prop}\label{Prop-SteinbAndChar}
One has an isomorphism of $F$-vector spaces
\[C_\har(V)^{\SL_2(A)}\cong\Hom_{\BF_q[\SL_2(\BF_q)]}(\overline\st,V).\]
\end{Prop}
The $F$-vector space structure on the right is induced from that on $V$, the $\SL_2(\BF_q)$-action on $V$ is induced from the map $\SL_2(\BF_q)\into\SL_2(A)$ that arises from the inclusion $\BF_q\into\BF_q[t]$.

\begin{proof}
The principal congruence subgroup of level $t$, denoted by $\Gamma(t)$, is normal in $\GL_2(A)$, and so by \autoref{lem:FactsOnSt}(d) and \autoref{rem:NormalSt} we have an exact sequence of $\BZ[\GL_2(A)]$-modules
\begin{equation}\label{Gammatex}
0\rightarrow \St\rightarrow \BZ[\overline\CT_1^{\Gamma(t)\dash\st}]\xrightarrow{\bar\partial_{\Gamma(t)}}\BZ[\CT_0^{\Gamma(t)\dash\st}]\rightarrow 0.
\end{equation}
Recall that the group $\GL_2(F_\infty)$ acts transitively on the tree $\CT$; the action is induced by the natural left action of $\GL_2(F_\infty)$ on $F_\infty^2$. The standard vertex $v_0$ of $\CT$ is given by the lattice class $v_0=[\mathcal{O}_\infty^2]$. The vertex $v_1=\left(\begin{smallmatrix}1 & 0\\
0 & \pi\\
\end{smallmatrix}\right) v_0=\left(\begin{smallmatrix}0 & 1\\
\pi & 0\\
\end{smallmatrix}\right) v_0$ is adjacent to $v_0$. The standard oriented edge $e_0$ is the oriented edge with $o(e_0)=v_0$ and $t(e_0)=v_1$. It is immediate from the definitions that one has
\[
\Stab_{\GL_2(F_\infty)}(v_0)=F_\infty^\ast \GL_2(\CO_\infty)\quad\text{and}\quad \Stab_{\GL_2(F_\infty)}(e_0)=F_\infty^\ast \CI,
\]
where $\Iw$ denotes the Iwahori subgroup of $\GL_2(\CO_\infty)$ of matrices whose reduction modulo $\pi$ is upper triangular.  This gives bijections $\GL_2(F_\infty)/F_\infty^*\GL_2(\CO_\infty)\to \CT_0$ and $\GL_2(F_\infty)/F_\infty^*\CI\to \CT^\ori_1$. Under this identification the map that associates to an oriented edge $e$ its origin $o(e)$ becomes the canonical map $\GL_2(F_\infty)/F_\infty^*\CI\rightarrow \GL_2(F_\infty)/F_\infty^*\GL_2(\CO_\infty)$. The inversion map $e\mapsto e^*$ on oriented edges $e$ is then induced by $g\mapsto g\left(\begin{smallmatrix}
0 & 1\\
\pi & 0\\
\end{smallmatrix}\right)$, and this gives a formula for $e\mapsto t(e)=o(e^*)$.

By \cite[Example~II.2.4.1]{Serre-Trees}, a fundamental domain for the action of $\GL_2(A)$ on $\CT$ is given by the subgraph with vertices 
$v_n=\left(\begin{smallmatrix}
1 & 0\\
0 & \pi^n\\
\end{smallmatrix}\right) v_0$, $n\in\BN_0$ and edges $e_n$ with $o(e_n)=v_n$ and $t(e_n)=v_{n+1}$ together with their opposites. By \cite[Section 2.5]{ccm}, we have that a vertex of the form $v_n$ is $\Gamma(t)$-stable if and only if $n=0$. Similarly, an edge of the form $e_n$ is $\Gamma(t)$-stable if and only if $n=0$. Let now $v\in\CT_0$ be an arbitrary vertex. Then we may find $g\in\GL_2(A)$ and $n\geq 0$ such that $v=gv_n$. Since $\Gamma(t)$ is normal in $\GL_2(A)$, it follows that $v\in\CT_0^{\Gamma(t)\dash\st}$ if and only if $v\in\GL_2(A)v_0$.  By the same argument, we obtain $\CT_1^{\ori,\Gamma(t)\dash\st}=\GL_2(A) e_{0}\cup\GL_2(A)e_{0}^\ast$.  From  \cite[II.1.6, Proposition~3]{Serre-Trees} we have
$\Stab_{\GL_2(A)}(v_0)=\GL_2(\BF_q)$ and $\Stab_{\GL_2(A)}(e_0)=\mathrm{B}_2(\BF_q)$. Thus, we obtain 
\[
\BZ[\CT_0^{\Gamma(t)\dash\st}]=\BZ[\GL_2(A)v_0]\cong \BZ[\GL_2(A)/\GL_2(\BF_q)]\cong \BZ[\SL_2(A)/\SL_2(\BF_q)].
\]
Similarly, by definition of $ \BZ[\overline\CT_1^{\Gamma(t)\dash\st}]$, see \autoref{lem:FactsOnSt}(b), we have
\[
\BZ[\overline\CT_1^{\Gamma(t)\dash\st}]\cong \BZ[\GL_2(A)e_0^\ast]\cong \BZ[\GL_2(A)/\mathrm{B}_2(\BF_q)]\cong \BZ[\SL_2(A)/\SB_2(\BF_q)].
\]
Combining this with \eqref{Gammatex}, we obtain an exact sequence of $\BZ[\SL_2(A)]$-modules,
\begin{equation}\label{Gammatex1}
0\rightarrow \St\rightarrow \BZ[\SL_2(A)/\SB_2(\BF_q)]\xrightarrow{\bar\partial_{\Gamma(t)}} \BZ[\SL_2(A)/\SL_2(\BF_q)]\rightarrow 0.
\end{equation}
Going back through the definitions, we see that the map $\bar\partial_{\Gamma(t)}$ is just given by the natural map $g\SB_2(\BF_q)\mapsto g\SL_2(\BF_q)$. Since \eqref{Gammatex} and thus also \eqref{Gammatex1} split as a sequence of $\BZ[\Gamma(t)]$-modules, we can tensor with $V$ over $\Gamma(t)$ to obtain an $\SL_2(A)/\Gamma(t)=\SL_2(\BF_q)$-equivariant exact sequence
\begin{equation}\label{Gammatex2}
0\rightarrow \St\otimes_{\Gamma(t)} V \rightarrow (\Ind^{\SL_2(A)}_{\SB_2(\BF_q)}\BZ)\otimes_{\Gamma(t)} V\xrightarrow{\bar\partial_{\Gamma(t)}} (\Ind^{\SL_2(A)}_{\SL_2(\BF_q)}\BZ)\otimes_{\Gamma(t)} V\rightarrow 0.
\end{equation}
To rewrite the middle and right term of \eqref{Gammatex2}, we let $G=\SL_2(A)$, $H=\SL_2(\BF_q)$, $H'=\SB_2(\BF_q)$ and $N=\Gamma(t)$, so that $G=N\rtimes H$. We claim that the following natural maps are $H$-equivariant isomorphisms:
\begin{itemize}
\item[(i)] $V|_H\rightarrow (\Ind_H^G\BZ)\otimes_NV$, $v\mapsto H\otimes v$.
\item[(ii)] $(\Ind_{H'}^H\BZ)\otimes_\BZ V|_H\rightarrow (\Ind_{H'}^G\BZ)\otimes_N V$, $hH'\otimes v\mapsto hH'\otimes v$.
\end{itemize}
For (i), this is easily verified directly: The injectivity is immediate as $\Ind_H^G\BZ$ is a free $\BZ[N]$-module. Similarly, the surjectivity follows immediately from $G=N\rtimes H$. The $H$-equivariance is obvious. The verification of (ii) can be done similarly upon observing that $\Ind_{H'}^G\BZ$ is a free $\BZ[N]$-module on the basis $H/H'$. By invoking these isomorphisms, we may rewrite \eqref{Gammatex2} as
\[
0\rightarrow \St\otimes_{\Gamma(t)} V \rightarrow (\Ind^{\SL_2(\BF_q)}_{\SB_2(\BF_q)}\BZ)\otimes_\BZ V|_{\SL_2(\BF_q)}\xrightarrow{\bar\partial_{\Gamma(t)}} V|_{\SL_2(\BF_q)}\rightarrow 0.
\]
By tracing back through the isomorphisms, one sees that $\bar\partial_{\Gamma(t)}$ is the map induced by \eqref{eqn-Stbg-Def}. We conclude that 
\begin{align*}
\St\otimes_{\Gamma(t)} V=\kernel\left( \Ind^{\SL_2(\BF_q)}_{\SB_2(\BF_q)}\BF_q\rightarrow \BF_q\right)\otimes_{\BF_q} V|_{\SL_2(\BF_q)}.
\end{align*}
Thus, by \eqref{eqn-Stbg-Def} we have $\St\otimes_{\Gamma(t)} V=\ost\otimes_{\BF_q}V|_{\SL_2(\BF_q)}$, and this yields
\[
C_{\har}(V)^{\SL_2(A)}\cong (\St\otimes_{\Gamma(t)} V)^{\SL_2(\BF_q)}=(\ost\otimes_{\BF_q}V)^{\SL_2(\BF_q)}.
\]
Now, by \autoref{Lem-BasicsOnLk}(d) with the action restricted to $\SL_2(\BF_q)\subset \GL_2(A)$, we have that $\ost=\overline L_{q-1}$ is self-dual and we deduce $C_{\har}(V)^{\SL_2(A)}\cong \Hom_{\BF_q[\SL_2(\BF_q)]}(\ost,V)$, proving the proposition.
\end{proof}

The following is an elementary result on $ \Hom_{\SL_2(\BF_q)}(\overline\st,V)$. Part (b) is essentially a special case of \cite[Cor.~1.3]{Kuhn-Mitchell}.
\begin{Prop} \label{cor:Hom-St-SS}\label{cor:Hom-StViaBorel}
Let $V$ be an $F[\SL_2(A)]$-module of finite $F$-dimension.
\begin{enumerate}
\item Denoting by $V^\ssi$ the semisimplification of $V$ considered as an $\SL_2(\BF_q)$-module, the $F$-dimension of $ \Hom_{\SL_2(\BF_q)}(\overline\st,V)$ is the multiplicity of $\overline L_{q-1}\otimes_{\BF_q}F$ in $V^\ssi$.
\item One has $\dim_F\Hom_{\SL_2(\BF_q)}(\overline\st,V)=\dim_F V^{\SB_2(\BF_q)} -\dim_F V^{\SL_2(\BF_q)} $.
\end{enumerate}
\end{Prop}
\begin{proof}
Part (a) follows from $\ost$ being a projective $\BF_q[\SL_2(\BF_q)]$-module. To see (b), apply the functor $\Hom_{\SL_2(\BF_q)}(\cdot,V)$ to the split exact sequence $0\to\ost\to\Ind_{\SB_2(\BF_q)}^{\SL_2(\BF_q)}\BF_q\to\BF_q\to0$ from \eqref{eqn-Stbg-Def}, to obtain the exact sequence 
\[ 0\to \Hom_{\SL_2(\BF_q)}(\BF_q,V)\to \Hom_{\SL_2(\BF_q)}(\Ind_{\SB_2(\BF_q)}^{\SL_2(\BF_q)}\BF_q,V)\to \Hom_{\SL_2(\BF_q)}(\ost,V) \to 0.\]
Applying Shapiro's Lemma to the middle term and computing dimensions proves (b).
\end{proof}
Because of \autoref{cor:Hom-St-SS}(a), we have $\dim  \Hom_{\SL_2(\BF_q)}(\overline\st,V_{k,l})=\dim  \Hom_{\SL_2(\BF_q)}(\overline\st,\Delta_{k-2})$. Now the latter dimension can be computed using \autoref{cor:Hom-StViaBorel}(b) and from the known structure of the $\SL_2(\BF_q)$- and $\SB_2(\BF_q)$-invariants of $\BF_q[X,Y]$: Define
\[f:=X^qY-XY^q, \ g:=\sum_{l=0}^q(X^{q-1})^l(Y^{q-1})^{q-l} \in \BF_q[X,Y].\]
\begin{Thm}[{\cite{Dickson}}]
For any field $\BF\supset \BF_q$ one has
\[\BF[X,Y]^{\SL_2(\BF_q)}=\BF[f,g]\]
and
\[\BF[X,Y]^{\SB_2(\BF_q)}=\BF[X,Y]\cap \BF[f,Y^{\pm (q-1)}]=\BF[f,(f/Y)^{(q-1)},Y^{(q-1)}].\]
\end{Thm}
\begin{proof}
The calculation of the $\SL_2(\BF_q)$-invariants is due to Dickson, see \cite{Dickson} or \cite[Thms.~2.1 and 2.2]{Kuhn-Mitchell}. We give the argument for the computation of the $\SB_2(\BF_q)$-invariants, and observe right away that the equality on the right is straightforward.

We fix $k\ge0$ and then identify $\BF[x,y]_k$ with $\BF[z]_{\le k}$, the polynomials in $z$ of degree at most $k$, by the map $f(X,Y)\mapsto f(z,1)$. Then the action $\SMat{a}bcd X^iY^j=(dX-bY)^i(-cX+aY)^j$ (with $i+j=k$ and $\SMat{a}{b}{c}{d}\in \SL_2(\BF_q)$) turns into 
\[\SMat{a}bcd z^i=\big(\frac{dz-b}{-cz+a}\big)^i(-cz+a)^k=(dz-b)^i(-cz+a)^j.\]
The action restricted to $\SB_2(\BF_q)$ extends to all of $\BF[z]$, and it is straightforward to verify that the invariants under the unipotent group $\SMat{1}{\BF_q}{0}1$ are given by the subring $\BF[z^q-z]$. Now note that $\SMat{a}00{a^{-1}}z^i=a^{k-2i}z^i$, and from this it follows that 
\[\SMat{a}00{a^{-1}}(z^q-z)^i=a^{k-2i}(z^q-z)^i.\] 
Hence the expression is invariant under all $a\in\BF_q^\times$ if and only if $2i\equiv k \pmod{q-1}$.

Let $i_0\ge0$ be minimal such that $2i_0\equiv k \pmod {q-1}$, 
let $q'=q-1$ if $q$ is even and $q'=\frac12(q-1)$ if $q$ is odd,
and let $h=(z^q-z)^{q'}$. Then we have
\[\BF[z]^{\SB_2(\BF_q)} =  (z^q-z)^{i_0}\BF[h].\]
Moreover $\BF[z]^{\SB_2(\BF_q)}_{\le k}= 0$ if $qi_0>k$, or if both $q$ and $k$ are odd, and $\BF[z]^{\SB_2(\BF_q)}_{\le k}= (z^q-z)^{i_0}\BF[h]_{ \le d_{k,q}}$ otherwise, where $d_{k,q}\ge0$ is the maximal integer such that $q(i_0 + d_{k,q} q') \le k$.

When converting back the answer to $\BF[X,Y]$ observe that
\[ (z^q-z)^{i_0}h^d = (X^qY-XY^q)^{i_0+dq'} \cdot (Y^{q'})^{\frac1{q'}(k-(q+1)i_0)-d(q+1)} , \]
where $X^qY - XY^q$ and $Y^{q'}$ are invariant under $\SB_2(\BF_q)$ and where $ k-(q+1)i_0\equiv k-2i_0\equiv 0\pmod {q-1}$, so that $\frac1{q'}(k-(q+1)i_0)\in\BZ$. In fact, if $q$ is odd, then $\frac1{q'}(k-(q+1)i_0)-d(q+1)$ is even. Note also that $\frac1{q'}(k-(q+1)i_0)-d(q+1)$ can be negative. This proves the claimed expression for $\BF_q[X,Y]^{\SB_2(\BF_q)}$.
\end{proof}

From the above, it follows that an $\BF$-basis of $\BF[X,Y]^{\SB_2(\BF_q)}$ is given by
\[\{ f^{i_0}\cdot (f/Y)^{(q-1)i_1}\cdot Y^{(q-1)i_2}\mid 0\le i_0\le q-2, i_1,i_2\in\BN_0\}.\]
Counting invariants in a fixed homogeneous degree, and using $\deg f=q+1$, the corresponding generating function for the dimensions of $\BF[X,Y]^{\SB_2(\BF_q)}$ is
\[ \sum_{i_0=0}^{q-2} (U^{q+1})^{i_0}\cdot \frac1{(1-U^{q-1})(1-U^{q(q-1)})}.\]
Subtracting from this the generating function for the dimensions of $\BF[X,Y]^{\SL_2(\BF_q)}$ gives
\[  \frac{1-U^{q^2-1}}{(1-U^{q+1})(1-U^{q-1})(1-U^{q(q-1)})} - 
 \frac1{(1-U^{q+1})(1-U^{q(q-1)})} 
=
 \frac{U^{q-1}}{(1-U^{q+1})(1-U^{q-1})} .\]

Noting that $V_{k+2,l}$ is related to $F[X,Y]_{k}$, the following is an immediate consequence.
\begin{Cor}\label{Cor:DimensionTeitel}
For any $l\in\BZ$ one has
\[\sum_{k\ge0} \dim C_\har(V_{k,l})^{\SL_2(A)} U^{k} 
=
 \frac{U^{q+1}}{(1-U^{q+1})(1-U^{q-1})}.\]
\end{Cor}
The above formula was obtained originally  via (rigid) analytic methods by Cornelissen building on work of Gekeler, see \cite[Prop.~4.3]{Cornelissen-Survey}. It only serves to demonstrate that via the above representation theoretic reformulation, one can also obtain this result; though perhaps via a more demanding route. For this reason Teitelbaum in \cite[p.~507ff.]{Teitelbaum-Poisson} promotes the viewpoint that one can use the analytic formulas, reinterpreted through representation theory to obtain formulas for spaces of invariants, or more concretely formulas for the multiplicity of $\ost$ in some $\SL_2(\BF_q)$-representations. 

In the case where the module is not $V_{k,l}$ but $L_k$ we see no analytic approach, since we have no explicit description of the corresponding quotient of $S_{k,l}(\SL_2(A))$. In the following subsection, we shall use representation theory to do the work for us. This is successful for some $q$, but we only see a general algorithm but no closed formula yet.

\begin{Rem}
The work of this subsection seems to be related to Gekeler's work \cite{Gekeler-Finite} on finite (Drinfeld) modular forms. 
\end{Rem}

\subsection{Dimension formulas via work of Reduzzi}
\label{subsec:Reduzzi}

Let $G:=\SL_2(\BF_q)$ with $q=p^e$. In \cite{Reduzzi-K0} Reduzzi proves the following result for the Grothendieck group $K_0(G)$ of the category of finite length representations of $\BF_q[G]$.
\begin{Thm}[{Reduzzi}]\label{Thm:Reduzzi}
Following the notation in  \autoref{Subsec3.3}, denote by $\overline{\Delta}_1$ the standard representation of $G$ on $\BF_q^2$ and let $T$ be an indeterminate over $\BZ$. The assignment $\overline{\Delta}_1\mapsto T$ induces a unique isomorphism of rings:
\[\phi\colon K_0(G)\stackrel\simeq\longto  \BZ[T]/(f^{[e]}(T)-T),\]
where $f^{[e]} = (f\circ \ldots\circ f)$ is the polynomial of $\BZ[T]$ of degree $q=p^e$ obtained by $e$-fold composition of the monic degree $p$ polynomial
\[f(T)=\sum_{j=0}^{\lfloor \frac{p}2\rfloor} (-1)^j\frac{p}{p-j} \binom{p-j}j T^{p-2j}.\]
Moreover one has the following:
\begin{compactenum}
\item The polynomial $f^{[e]}$ has the explicit expression $\sum_{j=0}^{\lfloor \frac{q}2\rfloor} (-1)^j\frac{q}{q-j} \binom{q-j}j T^{q-2j}$.
\item For $k\in\BZ_{\ge0}$ and $g_k(T):=\sum_{j=0}^{\lfloor \frac{k}2\rfloor} (-1)^j \binom{n-j}j T^{n-2j}\in\BZ[T]$ one has 
\[\phi(\overline{\Delta}_k)=g_k(T)\pmod{f^{[e]}(T)-T}.\]
\item One has $f^{[e]}(T)-T\equiv T^q-T \pmod p $ and $K_0(G)\otimes_\BZ\BF_p\cong \BF_p[T]/(T^q-T)$.
\item The $\BQ$-algebra $K_0(G)\otimes_\BZ\BQ$ is a product of number fields.
\item One has $K_0(G)\otimes_\BZ\BQ_q\cong\BQ_q^q$ for $\BQ_q$ the unramified extension of $\BQ_p$ of degree~$e$.
\end{compactenum}
\end{Thm}
\begin{Rem}
Experimentally, one also observes the following:
\begin{compactenum}
\item All factors of $K_0(G)\otimes_\BZ\BQ$ are totally real fields (Reduzzi).
\item The factors of $K_0(G)\otimes_\BZ\BQ$ are abelian extensions of $\BQ$; they only ramify at primes dividing~$q^2-1$.
\end{compactenum}
\end{Rem}
We now explain, how in principle \autoref{Thm:Reduzzi} allows one to give explicit expressions for the dimension of $\Hom_{\SL_2(\BF_q)}(\overline\st,L_k)$ and carry this through for $q=2,3,5$. For this, observe first that $L_k|_G=\overline{L}_k\otimes_{\BF_q}F$ and thus
\[
\dim_F \Hom_{\SL_2(\BF_q)}(\overline\st,L_k)= \dim_{\BF_q}\Hom_{\SL_2(\BF_q)}(\overline\st,\overline{L}_k).
\]
For simplicity of exposition, we assume $K_0(G)\otimes_\BZ\BR\cong\BR^q$, and we let $\sigma_0,\ldots,\sigma_{q-1}$ be the embeddings $K_0(G)\otimes_\BZ\BQ\to\BR$ -- we could also work with embeddings into $\BQ_q$, alternatively. 

A $\BZ$-basis of $K_0(G)$ is given by $\overline{\Delta}_k$, $k=0,\ldots,q-1$, or alternatively by $\overline{L}_k$, $k=0,\ldots,q-1$. Observe that the formulas of Bonnaf\'e in \autoref{Prop-V_k-versus-L_k} that allow one to go back an forth between the two bases continue to hold for $\BF_q[\SL_2(\BF_q)]$-representations and $k=0,\ldots,q-1$. Recall also that $\overline\st=\overline{L}_{q-1}$.

Note next that for an $\BF_q[G]$-representation $V$ one has $V^{(e)}\cong V$ for the $e$-th Frobenius twist. Let now $k$ be arbitrary in $\BZ_{\ge0}$ and write $k=\sum_{i\ge0} k_iq^i$ in its base $q$ expansion with $k_i\in\{0,\ldots,q-1\}$. Define furthermore $n_k(d):=\#\{i\ge0\mid k_i=d\}$ for $d=0,\ldots,q-1$. Then as an $\BF_q[G]$-representation one has
\[ \overline{L}_k\cong \bigotimes_{i\ge0} \overline{L}_{k_i}\cong  \bigotimes_{d=0}^{q-1} \overline{L}_d^{\otimes n_k(d)}.\]
To see the first isomorphism, observe that by the Steinberg Tensor Theorem, this holds for the base $p$ expansion, and if we insert twists. However after applying this theorem, we can group the tensor products in packages of $e$ consecutive digits and then use $V=V^{(e)}$ and run the Tensor Theorem backwards. The second isomorphism follows by simply regrouping the factors.

If we furthermore write any $d$ in its base $p$-expansion as $d=\sum_{r=0}^{e-1}d_rp^r$, then $\overline{L}_d=\bigotimes_{r=0}^{e-1}\overline{L}_{d_r}^{(r)}$.\footnote{Below the number $d$ varies. So we do regard $d\mapsto d_r$ as a function in $d$.} The next result is immediate from \cite[Proof of 3.1 (2nd paragr.), Lem.~3.2]{Reduzzi-K0} and induction.
\begin{Lem}
For $r,k\ge0$, one has 
\[\phi(\overline{\Delta}_k^{(r)})=g_k\circ f^{[r]}.\] 
In particular, for $0\le r\le e-1$ and $0\le k\le p-1$ one has $\phi(\overline{L}_k^{(r)})=g_k\circ f^{[r]}$.
\end{Lem}

Because $\overline\st$ is projective over $\BF_q[G]$, the dimension of $\Hom_{\BF_q[G]}(\overline\st, \overline{L}_k)$ is equal to the multiplicity of $\overline{L}_{q-1}=\overline{\Delta}_{q-1}$ when writing $\overline{L}_k$ in terms of the $\BZ$-basis $\overline{\Delta}_s$, $s=0,\ldots,q-1$ as $[\overline{L}_k]=\sum_{d=0}^{q-1} \lambda_k(d)[\overline{\Delta}_d]$. Let $\alpha_i$, $i=0,\ldots,q-1$, be the (pairwise distinct) roots of $f^{[e]}(T)-T$. Then $\underline{\sigma}:=(\sigma_d)_{d=0,\ldots,q-1}\colon K_0(G)\to\BR^q$ indicated above is given explicitly by
\[ \overline{\Delta}_k\mapsto (g_k(\alpha_j))_{j=0,\ldots,q-1}\in\BR^q. \]
Note that $\alpha\mapsto f(\alpha)$ induces a permutation of the roots $\alpha_i$ because $f^{[e]}=\id$ on these roots.

Under the map $(\sigma_i)_{i=0,\ldots,q-1}$ the equality $[\overline{L}_k]=\sum_{d=0}^{q-1} \lambda_k(d)[\overline{\Delta}_d]$ becomes
\begin{equation}\label{eqn:Lk-Eq}
\bigg(\prod_{d=0,\ldots,q-1} \Big(\prod_{r=0}^{e-1} g_{d_r}\big(f^{[r]}(\alpha_i)\big) \Big)^{n_k(d)}  \bigg)_{i=0,\ldots,q-1}=\sum_{d=0}^{q-1} \lambda_k(d) \Big( g_d(\alpha_i)\Big)_{i=0,\ldots,q-1}.
\end{equation}
The expression on the left may look complicated. But it simply says that after evaluating the basic polynomials for the $\overline{L}_k^{(r)}$ at the roots $\alpha_i$, computing $\overline{L}_k$ amounts to raising these values (separately) to suitable powers. A main point of Reduzzi's result is that it makes the a priori inexplicit product structure on $K_0(G)$ explicit after passing to the ring $\BR^q$.

To solve equation \eqref{eqn:Lk-Eq} one only needs linear algebra, since the images of the $\overline{\Delta}_k$, $k=0,\ldots,q-1$ form an $\BR$-basis of $\BR^q$. If the $\alpha_i$ are complicated, one has to think properly of how to solve the above system to the needed precision. This should not be too hard, since the tuple of $\lambda_k(d)$, $d=0,\ldots,q-1$ lies in $\BZ^d$, and we only need to approximate the $\lambda_k(d)$ to less than $1/2$. Also note that we are only interested in the coefficient of $\overline{\Delta}_{q-1}$. This element is invariant under Frobenius twist. Hence one could simply solve $\sum_{r=0}^{e-1}[\overline{L}_k^{(r)}]=e \lambda_k(q-1)[\overline{\Delta}_{q-1}] +\sum_{d=0}^{q-2} \lambda'_k(d)[\overline{\Delta}_d]$. Summing over the Frobenius twists may allow one to work over the trace field, i.e., the field that arises for $q=p$. We did not pursue this, since shall only present formulas for $q=2,3,5$. We give details for $q=3$, and only the solution in the other two cases.

For $q=p=3$ one has $f(T)=T^3-3T$, and hence $f^{[1]}(T)-T=T(T-2)(T+2)$. We set $\alpha_0=0$, $\alpha_1=2$, $\alpha_2=-2$. For the $\overline{L}_k=\overline{\Delta}_k$, $k=0,1,2$, we obtain $g_0(T)=1, g_1(T)=T, g_2(T)=T^2-1$ as polynomials over $\BZ[T]$. Write now a general $k\ge0$ as $\sum_{i\ge0} k_i3^i$ in its base $3$ expansion, and let $n_k(d)$ be the number of $i$ such that $k_i=d$, for $d=0,1,2$. Then $\underline\sigma{g_0}=(1,1,1)$, $\underline\sigma{g_1}=(0,2,-2)$ and $\underline\sigma{g_2}=(-1,3,3)$. For simplicity let $i_k=n_k(1)$ and $j_k=n_k(2)$. The linear system \eqref{eqn:Lk-Eq} now becomes
\[ (0,2,-2)^{i_k}\cdot(-1,3,3)^{j_k}= (\lambda_k(0),\lambda_k(1),\lambda_k(2))\cdot 
\left(\begin{array}{ccc}
1&1&1\\
0&2&-2\\
-1&3&3
\end{array}\right),\]
where the multiplication on the left is component-wise. This gives
\[ (0^{i_k}\cdot(-1)^{j_k},2^{i_k}3^{j_k},(-2)^{i_k}3^{j_k})\cdot\frac18 \left(\begin{array}{ccc}
6&0&-2\\
1&2&1\\
1&-2&1
\end{array}\right)= (\lambda_k(0),\lambda_k(1),\lambda_k(2))\]
with the convention $0^0=1$. This proves
\begin{Prop}\label{prop:DimLkQ3} If $q=p=3$ and $k\ge0$, then $\dim L_k=2^{i_k}3^{j_k}$ and 
\[\dim C_\har(L_k)^{\SL_2(A)}=\frac18\big(
(1+(-1)^{i_k}) \dim L_k -2\cdot 0^{i_k}(-1)^{j_k}
\big).\]
\end{Prop}

For $q=p=2$, similar but simpler arguments give:
\begin{Prop} \label{prop:DimLkQ2}
If $q=p=2$, then for all $k\ge0$ one has  $\dim L_k=2^{n_k(1)}$ and 
\[\dim C_\har(L_k)^{\SL_2(A)}=\frac13 \big( \dim L_k-(-1)^{n_k(1)} \big)
 \]
\end{Prop}

For $q=p=5$ a more elaborate argument gives
\begin{Prop} \label{prop:DimLkQ5}
If $q=p=5$ and $k\ge0$, then $\dim L_k=2^{n_k(1)}3^{n_k(2)}4^{n_k(3)}5^{n_k(4)}$ and 
\begin{eqnarray*}
\lefteqn{{\textstyle \dim C_\har(L_k)^{\SL_2(A)}}}\\
&=&{\textstyle \frac1{24}\big({6\cdot0^{n_k(1)+n_k(3)}(-1)^{n_k(2)}+(\dim L_k -4\cdot0^{n_k(2)}(-1)^{n_k(3)+n_k(4)})(1+(-1)^{n_k(1)+n_k(3)}) }\big).}
\end{eqnarray*}
\end{Prop}

\begin{Rem}
Dimension formulas for spaces of classical cusp forms of weight $k$ for $\SL_2(\BZ)$ can be written in a case-by-case way using congruences of $k-1=\dim_\BQ\Sym^{k-2}\BQ^2$ modulo $12$. The formulas in the above propositions can be written in a similar form using congruences of $\dim L_k$ modulo $(q^2-1)/\gcd(2,q^2-1)$ and of $k$ modulo $2$. However, our notation is more compact.

We think that the quantity $(q^2-1)/\gcd(2,q^2-1)$ should be interpreted as the index $[\PSL_2(A):\overline{\Gamma_1(t)}]$ where $\overline{\Gamma_1(t)}$ denotes the image of $\Gamma_1(t)$ in $\PSL_2(A)$, in analogy to the classical case. There the term $[\PSL_2(\BZ):\overline\Gamma]$ occurs as a factor in dimension formulas for classical modular forms of weight $k$ for congruence subgroups $\Gamma$ (up to an error $\CO(1)$), where now $\overline\Gamma$ denotes the image of $\Gamma$ in $\PSL_2(\BZ)$.
\end{Rem}

\begin{Rem}
An alternative approach to obtain the results in this subsection is to work out the Brauer characters of the $\overline{L}_k$ and use the orthogonality relations for such to compute the multiplicity of $\ost$ in $\overline{L}_k$. This allows one to recover the results of the present subsection. We plan to explore this further in future work.
\end{Rem}

For the sake of completeness, we note the following easy result for general odd $q$.

\begin{Prop} Let $q$ be odd and $k\geq 0$ odd. Then
\[
\dim C_\har(L_k)^{\SL_2(A)}=0.
\]
\end{Prop}
\begin{proof}
Since $L_k\subset \Delta_k$, we obtain a $\SL_2(A)$-equivariant surjection $V_{k+2,l}\rightarrow L_k$ (for any $l$). Thus, \autoref{thm:CHar-St-global}(b) implies
\[
\dim C_\har(L_k)^{\SL_2(A)}\leq \dim C_\har(V_{k+2,l})^{\SL_2(A)}.
\]
The latter is zero by \autoref{Cor:DimensionTeitel}.
\end{proof}

\begin{Rem}
The above proposition can also easily proved directly by considering the action of the matrix $\left(\begin{smallmatrix}
-1& 0\\ 0& -1
\end{smallmatrix}\right)\in \SL_2(A)$ directly on $C_\har(L_k)$ similar to the classical case.
\end{Rem}

\section{Investigations towards a Maeda-style conjecture for $q=3$ and special weights} \label{MaedaSection}
As indicated in \autoref{Section5}, the weights of the form $k=1+p^n-cp^{n-1}$ for $n\geq 0$ and $0\leq c <p-1$ are the most natural candidates for a simple formulation of a conjecture of Maeda-type as the relevant $\GL_2(F)$-representation $V_{k,l}$ is irreducible (for any $l$).  In particular, there are no non-zero hyperderivatives from lower weights. In this section, we want to investigate such a possible conjecture in the special case $q=3$. This is based on explicit computations using the computer algebra systems \texttt{Magma} and \texttt{Sage}. In the sequel, we consider $A=\BF_3[t]$ and $\Gamma=\GL_2(A)$.  We set $k_n=1+3^n$ and compute the characteristic polynomial $P_{t,n,l}\in F[X]$ of the Hecke operator $T_t$ associated to $t\in A$ acting on the spaces $S_{k_n,l}(\Gamma)$ for the two types $l=0,1$. We denote by $d_{n,l}$ the dimension of these spaces. Note that the weights of the form $1+3^n-3^{n-1}$ are not relevant for us as there are no non-zero Drinfeld cusp forms of odd weight for the group $\Gamma$.

\begin{table}[!htb]
\centering
\begin{tabular}{c|c|c|c|c|c}
$n$ & $k_n$ & $d_{n,0}$ & $d_{n,1}$ & Irreducible factors of $P_{t,n,0}$ & Irreducible factors of $P_{t,n,1}$\\
\hline
$2$ & $10$& $1$ & $1$ & $(X+t^4+t^2)$ & $(X-t)$\\
\hline
$3$ & $28$ & $3$ & $4$ & $(\text{irreducible of degree $3$})$ & \makecell{$(X-t)$\\
$(X+t^{13}-t^{11}-t^5)$\\$(X-t^{13}+t^{11}-t^5)$\\
$(X-t^{13}-t^{11}+t^5)$}\\
\hline
$4$ & $82$ & $10$ & $10$ &\makecell{$(X-t^{40}-t^{38}-t^{32}-t^{14})$\\ $(X+t^{40}+t^{38}+t^{32}-t^{14})$\\
 $(X+t^{40}+t^{38}-t^{32}+t^{14})$\\ $(X+t^{40}-t^{38}+t^{32}+t^{14})$\\ $(X-t^{40}+t^{38}+t^{32}+t^{14})$\\ $ (\text{irreducible of degree $5$})$}  & \makecell{$(X-t)$\\$(\text{irreducible of degree $9$})$}\\
\hline
$5$ & $244$ &$30$ & $31$&
\makecell{
$(\text{irreducible of degree $7$})$\\
$(\text{irreducible of degree $23$})$
}
& 
\makecell{
$(X-t)$\\
$(X+t^{121}+t^{119}+t^{113}+t^{95}+t^{41})$\\
$(X-t^{121}-t^{119}+t^{113}+t^{95}-t^{41})$\\
$(X-t^{121}+t^{119}-t^{113}+t^{95}-t^{41})$\\
$(X-t^{121}+t^{119}+t^{113}-t^{95}-t^{41})$\\
$(X-t^{121}-t^{119}-t^{113}+t^{95}+t^{41})$\\
$(X-t^{121}+t^{119}-t^{113}-t^{95}+t^{41})$\\
$(X-t^{121}-t^{119}+t^{113}-t^{95}+t^{41})$\\
$(X+t^{121}+t^{119}-t^{113}-t^{95}-t^{41})$\\
$(X+t^{121}-t^{119}+t^{113}-t^{95}-t^{41})$\\
$(X+t^{121}-t^{119}-t^{113}-t^{95}+t^{41})$\\
$(X+t^{121}-t^{119}-t^{113}+t^{95}-t^{41})$\\
$(\text{irreducible of degree $19$})$
}
\end{tabular}
\caption{Factorized Hecke polynomials}
\label{Table1}
\end{table}
Note that the factorizations in \autoref{Table1} take place in $F[X]$. Note also that \[\emph{all irreducible factors of degree $d>1$ that appear in \autoref{Table1} have Galois group $S_d$.}\]

Thus, we observe that aside the from the large number of linear factors appearing and the splitting into two large factors in weight $244$, we have a Maeda-style behavior. The factor $(X-t)$ showing up in all spaces of type $1$ corresponds to the unique cuspidal, but non double cuspidal form by \cite[Theorem 3.2]{Petrov-Aexpansion} and can therefore be omitted from the discussion. However, the other linear terms and their corresponding eigenvalues show remarkable symmetries, which we are going to explore further. In the sequel, we refer to these eigenvalues as \emph{special eigenvalues}.

The first key observation is that there are $n$ monomials appearing in the special eigenvalues in weight $k_n$ with non-zero coefficients and that their degrees are easily predictable as follows: If we divide a special eigenvalue in weight $k_n$ by $t^{k_n/2}$, the resulting polynomial in $t^{-1}$ has monomials $t^{-1},t^{-3},\dots,t^{-3^{n-1}}$, in particular it is  additive in $t^{-1}$. This leads to the natural question if the coefficients of the monomials (in this case just signs) are easily predictable as well. To make the situation more transparent, we use the short notation $(-1)^n_i$ for the set of all tuples in $(\BF_3^\times)^n$ with exactly $i$ entries being $-1$. Moreover, we denote the number of special eigenvalues in weight $k_n$ by $s_n$. Now, we can summarize the situation in \autoref{Table2}.

\begin{table}[!htb]
\centering
\begin{tabular}{c|c|c|c}
$n$ & $k_n$ & $s_n$ & Signs\\
\hline
$2$ & $10$ &$1$ & $(-1)^2_2$\\
\hline
$3$ & $28$ & $3$ & $(-1)^3_1$\\
\hline
$4$ & $82$ & $5$ & $(-1)^4_0,(-1)^4_3$\\
\hline
$5$ & $244$ &$11$ & $(-1)^5_2,(-1)^5_5$\\
\hline
$6$ &  $730$ & $21$ & $(-1)^6_1,(-1)^6_4$\\
\end{tabular}
\caption{Sign patterns for the special eigenvalues}
\label{Table2}
\end{table}

We first note that this suggests $s_n=(2^n-(-1)^n)/3$. However, we can say even more: The signs $(-1)^n_i$ in weight $k_n$ are precisely those with $0\leq i\leq n$ and $i\equiv 1-n \pmod {3}$. If we denote by $\Sigma_n\subset(\BF^\times_3)^n$ the union of these signs, by $S^0_{k_n,l}(\Gamma)$ the space of double cups forms in $S_{k_n,l}(\Gamma)$ and by $P^0_{t,n,l}$ the factor of $P_{t,n,l}$ corresponding to this space, we are able to formulate the following conjecture.
\begin{Conj}\label{Conj:RatEigen}
For all $n\geq 2$, with $k_n=1+3^n$, the elements in
\[
\left\{t^{k_n/2}\sum_{i=0}^{n-1}c_{n,i}t^{-3^i} \mid (
c_{n,i})_i\in \Sigma_n\right\}\subset A
\]
are eigenvalues for the Hecke operator associated to $t\in A$ acting on $S^0_{k_n,n}(\Gamma)$. The remaining irreducible factors $Q_i$ of $P^0_{t,n,n}$ are of degree larger than one with Galois group $S_{\deg Q_i}$. The factorization of $P^0_{t,n,n-1}$ only contains such polynomials $Q_i$.
\end{Conj}
\begin{Rem}
We would like to mention that we have also computed the matrix of $T_t$ in weight $730$. However, computing the characteristic polynomial of this matrix has proven to be not possible within reasonable running time using the build-in methods of standard computer algebra systems. Using reduction modulo large primes, we were nevertheless able to confirm that the only degree one factors in weight $730$ and type $0$ are those coming from the special eigenvalues.\end{Rem}

\begin{Rem}
At present, it is completely unclear how and why these special eigenvalues arise. While we are able to write down a purely combinatorial recipe predicting the eigenvalues, all attempts to get a grasp for the corresponding eigenforms have failed so far. We should note however that it is an easy exercise to check that at least the "sanity check" equality $\# \Sigma_n=(2^n-(-1)^n)/3$ holds. To formulate a full Maeda-type conjecture, we would need to predict the number of higher degree irreducible factors and their degrees. While our data suggests growth of the degrees with the weight, our data is too limited to formulate a precise conjecture. First results for $q=5$ indicate that this phenomenon could generalize. The case $q=2$ is in general more involved, as there are inseparability phenomena that seem to disappear for odd primes. It would also be of great interest to investigate the behaviour for more general weights $k$ when considering the representations $L_k$ (instead of $V_{k,l}$) as in the previous section. We plan to further investigate these points in future work. 
\end{Rem}
Note that the dimension of the ``special part" grows (conjecturally) much slower than the dimension of the full space of double cusp forms. Asymptotically, one would have
\[
s_n  \approx \frac{2^n}{3} \quad \text{and} \quad \dim S^0_{k_n,n}(\Gamma)\approx \frac{3^n}{16}.
\]

\begin{Rem}
Note that (classically) computational evidence suggests that it is sufficient to consider only one Hecke operator to observe the full Maeda-style behavior. In the function field case considered here one can say a bit more; one obtains the same results for degree one primes in $A$ different from $(t)$ rigorously by simple substitution.
\end{Rem}

\begin{Rem}
The special eigenform of weight $10$ in \autoref{Table1} has been studied in more detail by Shin Hattori in \cite[Section 4.2.2]{shin}. In particular, in his analysis of a $t$-adic family through this form he provides an explicit description as an element of $C_{\har}(V_{10,0})^{\Gamma_1(t)}\simeq V_{10,0}$.
\end{Rem}

\appendix

\renewcommand{\thesubsection}{\Alph{section}.\arabic{subsection}}
\setcounter{subsection}{0}

\section{Appendix}

\label{Section3}

The aim of this section is to provide some auxiliary cohomological results used in \autoref{Section2}. In \ref{Subsec2.1} we consider the local setting and in \ref{Subsec2.2} the adelic setting. 
\subsection{Tor exactness for certain modules over group rings}
\label{Subsec2.1}

For a finite group $G$ and a $\BZ[G]$-module $X$ let $X_G =X/\{gx-x\mid x\in X,g\in G\}$ and $X^G=\{x\in X\mid \forall g\in G: gx=x\}$ denote the modules of $G$-coinvariants and $G$-invariants of $X$, respectively. Recall (see \cite[VIII.1]{Serre-LocalFields}) that one has a well-defined {\em norm map}
\begin{equation}\label{eq:NormMap}
\nu_{X,G}\colon H_0(G,X)=X_G \longto H^0(G,X)= X^G, x \mapsto \sum_{g\in G}gx.
\end{equation}
Its kernel and cokernel are the Tate cohomology groups $\hat H^{-1}(G,X)$ and $\hat H^0(G,X)$, respectively. Recall also from \cite[IX.3]{Serre-LocalFields} that $X$ is called cohomologically trivial if $\hat H^i(H,X)=0$ for all $i\in\BZ$ and for all subgroups $H$ of~$G$; observe that $\nu_{X,G}$ is an isomorphism if $X$ is cohomologically trivial.

Let now $\Gamma$ be any group, let $\Gamma'\subset\Gamma$ be a normal subgroup of finite index and set $\bar\Gamma=\Gamma/\Gamma'$. For any subgroup $H\subset\bar \Gamma$, we define $\Gamma_{H}$ as the inverse image of $H$ under the canonical projection $\Gamma\to\bar\Gamma,\gamma\mapsto\bar\gamma$. Let $A$ be a commutative ring, and let $M$ be a right $\BZ[\Gamma]$-module and let $N$ is a left $A[\Gamma]$-module. It is straightforward to verify that $M\otimes_{\Gamma'} N$ is a left $A[\bar\Gamma]$-module via
\begin{equation}\label{eq:ModStructure}
 \bar\gamma\cdot  m\otimes_{\Gamma'}n = m\gamma^{-1}\otimes_{\Gamma'}\gamma n,
\end{equation}
for $m\in M$, $n\in N$ and $\gamma\in\Gamma$, and that the map 
\begin{equation}\label{eq:tensor-coinv}
\lambda_{M,N}\colon(M\otimes_{\Gamma'}N)_{\bar\Gamma}\to M\otimes_\Gamma N, [m\otimes_{\Gamma'} n]\mapsto m\otimes_\Gamma n
\end{equation}
is an $A$-module isomorphism. We also use the term {\em norm map} for the composition $\nu_{M\otimes_{\Gamma'}N,\bar\Gamma}\circ \lambda_{M,N}^{-1}\colon M\otimes_\Gamma N\to (M\otimes_{\Gamma'}N)^{\bar \Gamma}$. The aim of this section is the following technical result as well as the variant \autoref{prop:AdelicGroupHomol} in an adelic setting.
\begin{Prop}\label{prop:OnGroupHomology}
Suppose that $A$ is a ring of characteristic $p>0$ and that $M$ is projective as a $\BZ[\Gamma_H]$-module for any $p$-Sylow-subgroup $H$ of $\bar\Gamma$. Then the following hold.
\begin{compactenum}
\item
\label{prop:OnGroupHomologya}
The module $M\otimes_{\Gamma'}N$ is cohomologically trivial as an $A[\bar\Gamma]$-module.
\item
\label{prop:OnGroupHomologyb}
The norm map $M\otimes_\Gamma N\stackrel{{\scriptscriptstyle \eqref{eq:tensor-coinv}}}=(M\otimes_{\Gamma'}N)_{\bar \Gamma} \to (M\otimes_{\Gamma'} N)^{\bar\Gamma}$ is an isomorphism.
\item
\label{prop:OnGroupHomologyc}
For any short exact sequence $0\to N'\to N\to N''\to 0$ of left $A[\Gamma]$-modules the sequence
\[ \xymatrix{
0\ar[r]& M\otimes_{\Gamma} N'\ar[r]& M\otimes_{\Gamma} N\ar[r]&M\otimes_{\Gamma} N''\ar[r]&0\\
}
\]
is exact. 
\end{compactenum}
\end{Prop}
For the proof of the proposition we first need a lemma.
\begin{Lem}\label{Lem-BarGammaAction}
If $M$ is a projective $\BZ[\Gamma]$-module, then $M\otimes_{\Gamma'} N$ is a cohomologically~trivial $A[\bar\Gamma]$-module.
\end{Lem}
\begin{proof}
Observe first that we may assume that $M$ is a free $\BZ[\Gamma]$-module; this holds because cohomological triviality is inherited by direct summands and because any projective module is a direct summand of a free module. Next note that cohomological triviality is also preserved under filtered direct limits and under direct sums. Since our $M$ will be a limit of modules of the form $\BZ[\Gamma]^n$, it will thus suffice to prove the lemma assuming that $M=\BZ[\Gamma]$. For $M=\BZ[\Gamma]$ one verifies that that the map
\[ \BZ[\Gamma]\otimes_{\Gamma'}N\to A[\bar \Gamma]\otimes_AN,\gamma\otimes n\mapsto \bar\gamma^{-1}\otimes\gamma n\]
is an isomorphism of $A[\bar\Gamma]$-modules with inverse $\bar\gamma\otimes n\mapsto \gamma^{-1}\otimes\gamma n$ (the right hand expression is independent of the chosen representative $\gamma\in\Gamma$ of $\bar\gamma$); here $\bar\Gamma$ acts on the left module as defined in \eqref{eq:ModStructure}, and on the right module via its left action on $A[\bar\Gamma]$ (by left translation) and the trivial action on~$N$. Therefore $M\otimes_{\Gamma'}N$ is an induced $A[\bar\Gamma]$-module, and hence cohomologically trivial \hbox{by \cite[IX.3]{Serre-LocalFields}}. 
\end{proof}

\begin{proof}[Proof of \autoref{prop:OnGroupHomology}]
For the proof of part (a) we set $L:=M\otimes_{\Gamma'}N$. We apply \autoref{Lem-BarGammaAction} to $\Gamma_H\supset \Gamma'$ for any $p$-Sylow subgroup $H$ of $\bar\Gamma$. This shows that $L$ is cohomologically trivial as an $A[H]$-module for any $p$-Sylow subgroup $H$ of $\bar\Gamma$. Note next that $L$ is $p$-torsion since $\Char A=p$. Thus we may apply the criterion from \cite[IX.5~Thm.~8]{Serre-LocalFields} to deduce that $L$ is cohomologically trivial as an $A[\bar\Gamma]$-module.

Part (b) is immediate from part (a) and the observation on the norm map in the paragraph following the definition of the norm map in~\eqref{eq:NormMap}. For part (c) observe that the sequence
\[ \xymatrix{
0\ar[r]& M\otimes_{\Gamma'} N'\ar[r]& M\otimes_{\Gamma'} N\ar[r]&M\otimes_{\Gamma'} N''\ar[r]&0\\
}
\]
is exact, because $M$ is $\BZ[\Gamma']$-projective. We now apply part (b) to conclude.
\end{proof}

\subsection{Adelization of \autoref{prop:OnGroupHomology}}
\label{Subsec2.2}

In this section, we provide an abstract setting in which we formulate and prove an analog of \autoref{prop:OnGroupHomology}. This is used in \autoref{Subsec2.3} for the adelic Steinberg module. 

\medskip

Let $G$ and $H$ be groups with $H$ finite. Let $S$ be a bi-$G\times H$-set, by which we mean that $S$ carries a left action by $G$ and a right action by $H$ and that the two actions commute. Note that then $G\backslash S=\{Gs\mid s\in S\}$ carries a right $H$-action and that $S/H=\{sH\mid s\in S\}$ carries a left $G$-action. We make the assumption 
\begin{equation}\label{eq:AssOnStabHs}
\forall s\in S\colon \Stab_H(s)=\{e_H\}
\end{equation}
on $H$-stabilizers,  so that $S$ is a free right $H$-set. Let furthermore $M$ be a left $\BZ[G]$-module. We define the right $\BZ[G\times H]$-module $M[S]$ as follows: as a $\BZ$-module it is
\[ M[S]:=\bigoplus_{s\in S} M.\]
Elements are denoted as sequences $(m_s)_{s\in S}$ with $m_s\in M$, where it is understood that almost all $m_s=0$. The $G\times H$-right action is defined by 
\begin{equation}\label{eqn:DefGH-Mod}
 (m_s)_{s\in S} \cdot (g,h):=(g^{-1}\cdot m_{gsh^{-1}})
\end{equation}
for $(g,h)\in G\times H$ and $(m_s)\in M[S]$. Note that 
\begin{equation}
M[S]\otimes_H\BZ=M[S]_H=M[S/H]. 
\end{equation}
We call a subgroup $H'\subset H$ {\em small for $M[S]$} if $M$ is projective as a $\BZ[\Stab_G(sH')]$-module for all $s\in S$. Let $A$ be a ring of characteristic $p>0$. Let  $N$ be a left $A[G]$-module. 

The following result is an analog of \autoref{prop:OnGroupHomology}. 
\begin{Prop}\label{prop:AdelicGroupHomol}
If all $p$-Sylow subgroups of $H$ are small for $M[S]$, then the following hold:
\begin{compactenum}
\item
\label{prop:AdelicGroupHomola}
The module $M[S]\otimes_GN$ is cohomologically trivial as an $A[H]$-module.
\item
\label{prop:AdelicGroupHomolb}
The norm map $\nu\colon  M[S/H]\otimes_{G}N \longto (M[S]\otimes_GN)^H$ is an isomorphism.
\item
\label{prop:AdelicGroupHomolc}
For any short exact sequence $0\to N'\to N\to N''\to 0$ of left $A[G]$-modules, the sequence
\[ \xymatrix{
0\ar[r]& M[S/H]\otimes_G N'\ar[r]& M[S/H]\otimes_G N\ar[r]&M[S/H]\otimes_G N''\ar[r]&0\\
}
\]
is an exact sequence of $A$-modules. 
\end{compactenum}
\end{Prop}

\begin{Prop}
\label{prop:AbstractHeckeAction}
Let $T$ be an $A$-subalgebra of $\End_{\BZ[G\times H]}(M[S])$. Then the following hold:
\begin{compactenum}
\item On $(M[S]\otimes_G N)^H$ and $M[S/H]\otimes_G N$, an $A$-linear action of $T$ is given by $t\cdot(m\otimes n):=(tm)\otimes n$.
\item The norm map from \autoref{prop:AdelicGroupHomol}(b) is $T$-linear.
\item Let $N\to N'$ be a map of $A[G]$-modules. The induced maps $(M[S]\otimes_G N)^H\to (M[S]\otimes_G N')^H$ and $M[S/H]\otimes_G N \to M[S/H]\otimes_G N'$ are $T$-equivariant.
\item 
Suppose that $M[S]$ satisfies the condition from \autoref{prop:AdelicGroupHomol}.
Then for any short exact sequence $0\to N'\to N\to N''\to 0$ of left $A[G]$-modules, the sequence
\[ \xymatrix{
0\ar[r]& M[S/H]\otimes_G N'\ar[r]& M[S/H]\otimes_G N\ar[r]&M[S/H]\otimes_G N''\ar[r]&0\\
}
\]
is an exact sequence of $T$-modules.
\end{compactenum}
\end{Prop}
We first prove \autoref{prop:AbstractHeckeAction} using \autoref{prop:AdelicGroupHomol}.
\begin{proof}[Proof of \autoref{prop:AbstractHeckeAction}]
For part (a) we need to show that $\{ m\otimes g n-m g\otimes n\mid m\in M[S], n\in N, g\in G\}$ is preserved by multiplication by $t\in T$. This follows from 
\[t\cdot(m\otimes g n-m g\otimes n)=tm\otimes gn-t(m g)\otimes n\stackrel{t\in\End_{\BZ[G\times H]}(M[S])}=tm\otimes g n-(tm)g \otimes n.\]
Since $t\in\End_{\BZ[G\times H](M[S])}$, it also preserves taking $H$-invariants. The above computation is also valid when replacing $M[S]$ by $M[S/H]$ completing the proof of (a).
Parts (b) and (c) follow analogously. By (c), the maps in (d) are $T$-equivariant, the exactness follows \autoref{prop:AdelicGroupHomol}(c).
\end{proof}

In the remainder of this subsection, we give the proof of \autoref{prop:AdelicGroupHomol}. We define for $s\in S$ the groups
\[
H_s:=\Stab_H(Gs),\quad
\Gamma_s':=\Stab_G(s),\quad
\Gamma_s:=\Stab_G(sH).
\]
\begin{Lem}\label{lem:App-Lem1}
The following properties hold:
\begin{compactenum}
\item \label{lem:App-Lem1a} The subgroup $\Gamma_s'$ of $\Gamma_s$ is a normal subgroup.
\item \label{lem:App-Lem1b} The groups $\bar\Gamma_s:=\Gamma_s/\Gamma'_s$ and $H_s$ act simply transitively on $sH_s$ from the left and from the right, respectively.
\item \label{lem:App-Lem1c} For all $\bar\gamma\in\bar\Gamma_s$ there is a unique $h_{\bar\gamma}$ in $H_s$ such that $\bar\gamma s=sh_{\bar\gamma}^{-1}$, and the map
\[\bar\Gamma_s\to H_s,\bar\gamma\mapsto h_{\bar\gamma}\]
is an isomorphism of groups.
\item \label{lem:App-Lem1d} The map $\iota_s\colon \Gamma_s\to\Gamma_s\times H_s,\gamma\mapsto (\gamma,h_{\bar\gamma})$ is a group monomorphism and $\Stab_{\Gamma_s\times H_s}(s)=\iota_s(\Gamma_s)$.
\end{compactenum}
\end{Lem}
\begin{proof} (a) For $\gamma\in \Gamma_s$ let $h_\gamma\in H$ be such that $\gamma s=sh_\gamma$, and hence $\gamma^{-1}s=sh_\gamma^{-1}$. Then for 
$\gamma'\in\Gamma_s'$ we have 
\[\gamma\gamma'\gamma^{-1}s=\gamma\gamma'sh_\gamma^{-1}=\gamma sh_\gamma^{-1}=s,\]
the point being that $G$- and $H$- actions commute. This shows (a).

(b) The claim on $H_s$ is immediate from our hypothesis $\Stab_H(s)=\{e_H\}$. Also, it is clear that $\bar\Gamma_s$ acts simply transitively on the $\Gamma_s$-orbit of $s$ under $\Gamma_s$, since $\Gamma_s'$ is the stabilizer of this action. It remains to show that $\Gamma_s s=sH_s$. For ``$\subseteq$'', let $\gamma$ be in $\Gamma_s$ and $h_\gamma\in H$ as in the proof of  (a). Then $(Gs)h_\gamma=G\gamma s=G s$, and so $h_\gamma\in H_s$. For ``$\supseteq$'', let $h$ be in $H_s$ and $\gamma_h\in G$ so that $sh=\gamma_h s$. Then $\gamma_h sH=sh H=sH$, and so $\gamma_h\in \Gamma_s$, and this concludes the proof of~(b) .

(c) The existence and uniqueness of $h_{\bar\gamma}$ is immediate from (b). The fact that the map given is a group isomorphism also follows from (b), upon noting that $\bar\Gamma_s$ acts from the left and $H_s$ from the right, which explains the inversion in the definition of $h_{\bar\gamma}$.

(d) The first part is obvious, so that we only have to prove $\Stab_{\Gamma_s\times H_s}(s)=\iota_s(\Gamma_s)$. For ``$\subseteq$'' suppose that $(\gamma,h)\in\Gamma_s\times H_s$ fixes $s$, i.e, that $\gamma s h=s$. Then $\gamma s=s h^{-1}$, and so $h=h_{\bar \gamma}$, and hence $(\gamma,h)=\iota_s(\gamma)$. The proof of ``$\supseteq$''  is analogous, and we omit it. 
\end{proof}

Next note that since $GsH\subset S$ is a bi-$G\times H$-subset, we have an inclusion of $\BZ[G\times H]$-modules
\begin{equation}\label{eq:M[GsH]}
M[GsH]\subset M[S]
\end{equation}
\begin{Lem}
\label{lem:App-Lem2}
We have the following isomorphisms of right $\BZ[G\times H]$-modules
\begin{compactenum}
\item \label{lem:App-Lem2a}
$M[S]=\bigoplus_{i\in I} M[Gs_iH]$ for a set of representatives $\{s_i\}_{i\in I}$ in $S$ of $G\backslash S/H$.
\item \label{lem:App-Lem2b}
$M[GsH]\cong \Ind_{G\times H_s}^{G\times H}M[GsH_s]$
\item \label{lem:App-Lem2c}
$M[GsH_s]\cong \Ind_{\Gamma_s}^{G\times H_s} M$, where $\Gamma_s\to G\times H_s$ is given by $\iota_s$ from \autoref{lem:App-Lem1}.
\end{compactenum}
\end{Lem}
\begin{proof}
Part (a) is clear from the definitions and the partition $GSH=\coprod_{i\in I} Gs_iH$. For (b) and (c) recall that for a group $\CG$ with subgroup $\CH$, a right $\BZ[\CG]$-module $V$ is induced from a $\BZ[\CH]$-submodule $W\subset V$, if and only the map $\bigoplus_{g\in \CG/\CH} Wg\to V$, given by inclusion $Wg\subset V$ on each summand, is an isomorphism of $\BZ$-modules; see \cite[\S~3.3]{Serre-RepOfGroups}. 

To see (b), note now that $W:=M[GsH_s]$ is a $G\times H_s$-submodule of  $V:=M[GsH]$, and because $GsH$ is a disjoint union $\coprod_{\bar h \in H_s\backslash H} GsH_sh$, we have $\bigoplus_{h\in H_s\backslash H} M[GsH_s]h=M[GsH]$ as $\BZ$-modules, so that the criterion just stated applies.

To prove (c), denote by $V$ the $G\times H_s$-module $M[GsH_s]$ and write $W$ for the $\Gamma_s$-module $M[\{s\}]=M$; note that by \autoref{lem:App-Lem1}(d), $M[\{s\}]$ is indeed a $\Gamma_s$-module; it is also clear that $W$ is a $\Gamma_s$-submodule of $V$. Clearly $GsH_s$ is equal to the disjoint union $\bigcup_{(g,h)\in (G\times H_s)/\Stab_{G\times H_s}(s)} (g,h)\{s\}$, and the proof of (c) proceeds now as that for (b).
\end{proof}

\begin{proof}[Proof of \autoref{prop:AdelicGroupHomol}]
(a) Arguing as in the proof of \autoref{prop:OnGroupHomology}(a), it suffices to prove (a) for any $p$-Sylow subgroup of $H$; i.e., we may assume that $H$ is small. Since cohomological triviality passes from modules to their direct sum, by \autoref{lem:App-Lem2}(a) we may assume $S=GsH$ for some (or any) $s\in S$. By \autoref{lem:App-Lem2}(b) and~(c), we thus have $M[S]= \Ind_{\Gamma_s}^{G\times H} M$, where $\Gamma_s\to G\times H$ is given by $\iota_{s}$ from \autoref{lem:App-Lem1}. Since $\Ind_{H_s}^H$ maps cohomologically trivial $H_s$-modules to cohomologically trivial $H$-modules (see the beginning of \cite[IX.3]{Serre-LocalFields} ), we may assume $H=H_s$ when proving that $M[GsH]\otimes_G N$ is cohomologically trivial as an $H$-module. Arguing as in the proof of \autoref{Lem-BarGammaAction}, we may assume that $M\cong\BZ[\Gamma_s]$. In this case we have
\[\Ind_{\Gamma_s}^{G\times H} M\otimes_{\Gamma_s}N=\BZ[G\times H]\otimes_{G}N\cong\BZ[H]\otimes_\BZ N\cong \Ind_{\{e_H\}}^H N\]
as an $A[H]$-module; as an induced module, the latter is obviously cohomologically trivial, proving (a). 

Part (b) is immediate from part (a) and the observation after the definition of the norm map in~\eqref{eq:NormMap}. For part (c) we claim that $M[S]\otimes_G\ublank$ is an exact functor, i.e., that $M[S]$ is projective as a $\BZ[G]$-module. From this it follows that the sequence $ 0\to M[S]\otimes_G N'\to  M[S]\otimes_G N\to  M[S]\otimes_G N''\to 0$ is exact. By (a) it will remain so after taking $H$-invariants. From (b) the assertion in (c) will follow. To see the claim, note that $H$ is irrelevant, and so by our hypothesis we may assume that $H$ is small. Then the projectivity follows as in the proof of (a).
\end{proof}

\newpage
\bibliographystyle{hep}
\bibliography{FF-Maeda}
\vspace{\baselineskip}
\textsc{Gebhard B\"ockle}\par
Interdisziplin\"ares Zentrum f\"ur wissenschaftliches Rechnen\\
Ruprecht-Karls-Universit\"at Heidelberg\\
Im Neuenheimer Feld 205\\
69120 Heidelberg\\
Germany\par
\texttt{gebhard.boeckle@iwr.uni-heidelberg.de}
\\[2em]
\textsc{Peter Mathias Gr\"af}\par
Interdisziplin\"ares Zentrum f\"ur wissenschaftliches Rechnen\\
Ruprecht-Karls-Universit\"at Heidelberg\\
Im Neuenheimer Feld 205\\
69120 Heidelberg\\
Germany\par
\texttt{peter.graef@iwr.uni-heidelberg.de}
\\[2em]
\textsc{Rudolph Perkins}\par
Department of Mathematics\\
University of California, San Diego (UCSD)\\
9500 Gilman Drive \# 0112\\
La Jolla, CA  92093-0112\\
United States of America (USA)\par
\texttt{rperkins@ucsd.edu}

\end{document}